\documentclass[12pt,reqno]{amsart}

\usepackage{enumerate}
\usepackage{enumitem}
\usepackage{mathtools}
\usepackage{amsxtra}
\usepackage{amsmath}
\usepackage{amscd}
\usepackage{amssymb}
\usepackage{amsfonts}
\usepackage[all]{xy}
\usepackage{amsthm} 
\usepackage{color}
\usepackage{upgreek}
\usepackage{latexsym}
\usepackage[english]{babel} 
\usepackage{hyperref}
\usepackage{nicefrac}
\usepackage{euscript}
\usepackage[foot]{amsaddr}
\usepackage{amscd}
\usepackage[sort]{cite}           
\usepackage{fancyhdr}          
\usepackage{amsfonts}
\usepackage{color}
\usepackage{upgreek}
\usepackage{graphicx,xcolor}
	\usepackage{floatrow}
\usepackage[labelformat=simple,labelfont={}]{subcaption}
\usepackage{tikz}
	\usetikzlibrary {arrows.meta, decorations.markings}
\usepackage{hyperref}
\usepackage{stmaryrd}
\usepackage{nicefrac}
\usepackage{wasysym}
\usepackage{import}
\usepackage{xifthen}
\usepackage{pdfpages}
\usepackage{transparent}	
\usepackage{calc}

\usetikzlibrary{decorations.pathmorphing}

\xyoption{all}

\newsavebox{\mybox}
\newlength{\mydepth}
\newlength{\myheight}

\newenvironment{sidebar}%
{\begin{lrbox}{\mybox}\begin{minipage}{\textwidth}}%
{\end{minipage}\end{lrbox}%
 \settodepth{\mydepth}{\usebox{\mybox}}%
 \settoheight{\myheight}{\usebox{\mybox}}%
 \addtolength{\myheight}{\mydepth}%
 \noindent\makebox[0pt]{\hspace{-20pt}\rule[-\mydepth]{1pt}{\myheight}}%
 \usebox{\mybox}}

\theoremstyle{plain}
\newtheorem{theorem}[equation]{Theorem} 
\newtheorem{lemma}[equation]{Lemma}
\newtheorem{lemma-definition}[equation]{Lemma-Definition}
\newtheorem{definition-lemma}[equation]{Definition-Lemma}
\newtheorem{cor}[equation]{Corollary}

\newtheorem{prop}[equation]{Proposition}

\newtheorem{corollary}[equation]{Corollary}

\theoremstyle{definition}

\theoremstyle{remark}

\newtheorem{rem}[equation]{Remark}
\newtheorem{rems}[equation]{Remarks}

\newtheorem{ex}[equation]{Example}
\newtheorem{exas}[equation]{Examples}

\def\AAA{{\mathbb{A}}}

\def\C{\mathbb{C}}

\def\H{\mathbb{H}}

\def\N{\mathbb{N}}

\def\R{\mathbb{R}}

\def\Z{\mathbb{Z}}

\def\frakB{\mathfrak{B}}

\def\frakg{\mathfrak{G}}

\def\frakS{\mathfrak{S}}

\def\calC{\mathcal{C}}

\def\calE{\mathcal{E}}
\def\calF{\mathcal{F}}
\def\calG{\mathcal{G}}
\def\calH{\mathcal{H}}

\def\calK{\mathcal{K}}
\def\calL{\mathcal{L}}

\def\calN{\mathcal{N}}

\def\calP{\mathcal{P}}
\def\calQ{\mathcal{Q}}

\def\calS{\mathcal{S}}

\def\calV{\mathcal{V}}
\def\calW{\mathcal{W}}

\def\frakb{\mathfrak{b}}

\def\frakg{\mathfrak{g}}
\def\frakh{\mathfrak{h}}

\def\frakl{\mathfrak{l}}

\def\frakn{\mathfrak{n}}

\def\frakp{\mathfrak{p}}
\def\frakq{\mathfrak{q}}
\def\frakr{\mathfrak{r}}

\def\frakz{\mathfrak{z}}

\def\bfa{\mathbf{a}}

\def\bfk{\mathbf{k}}

\def\bfm{\mathbf{m}}
\def\bfn{\mathbf{n}}

\def\bfu{\mathbf{u}}

\def\bfI{\mathbf{I}}
\def\bfJ{\mathbf{J}}
\def\bfK{\mathbf{K}}
\def\bfL{\mathbf{L}}

\def\be{\begin{equation}}
\def\BFun{{\on{BFun}}}

\def\Br{{\on{Br}}}
\def\bs{\backslash}

\def\ch{{\on{char}\, }}

\def\Coker{{\on{Coker}}}
\def\Cox{{\on{Cox}}}
\def\coxcat{{\boldsymbol{\mathcal{C}}}}

\def\ee{\end{equation}}
\def\eps{{\varepsilon}}

\def\Fl{{\on{Fl}}}
\def\frakSB{{\frak{SB}}}
\def\Fun{{\on{Fun}}}

\def\GBialg{{\on{GBialg}}}

\def\gl{{\frakg\frakl}}

\def\hW{{W\backslash \frakh}}

\def\Id{{\on{Id}}}
\def\Ind{{\on{Ind}}}
\def\Irr{{\on{Irr}}}

\def\k{{\bfk}}
\def\Kar{{\on{Kar}}}

\def\Lie{{\on{Lie}}}

\def\lra{\longrightarrow}

\def\Mat{{\on{Mat}}}
\def\Mod{\hskip -0.05cm{\on {Mod}}}
\def\Mor{{\on{Mor}}}

\def\Ob{\operatorname{Ob}\nolimits}
\def\ol{\overline}
\def\on{\operatorname}

\def\Perv{\operatorname{Perv}\nolimits}
\def\phi{{\varphi}}

\def\QS{{\calQ\frakS}}

\def\Res{{\operatorname{Res}}}

\def\SBialg{{\on{SBialg}}}
\def\Span{{\on{Span}}}
\def\Sh{{\on{Sh}}}
\def\Stab{{\on{Stab}}}

\def\uInd{\ul{\Ind}}
\def\ul{\underline}
\def\uRes{\ul{\Res}}

\def\wt{\widetilde}

\def\={\,\, {\simeq}\, \,}
\def\1{{\bf 1}}

\def\Vect{\operatorname{Vect}\limits}

\def\top{{\operatorname{top}\nolimits}}
\def\alg{{\operatorname{alg}\nolimits}}

\def\top{{\operatorname{top}\nolimits}}
\def\top{{\operatorname{top}\nolimits}}

\def\-{{\operatorname{-}\!}}

\def\Im{\operatorname{Im}\nolimits}
\def\Ker{\operatorname{Ker}\nolimits}
\def\Hom{\operatorname{Hom}\nolimits}

\def\End{\operatorname{End}\nolimits}

\def\Ind{\operatorname{Ind}\nolimits}

\def\Inv{\operatorname{Inv}\nolimits}

\def\Spec{\operatorname{Spec}\nolimits}

\def\Fl{{\operatorname{Fl}\nolimits}}

\def\preceq{{\operatorname\preccurlyeq\nolimits}}

\numberwithin{itemcounter}{subsection}
\numberwithin{equation}{subsection}
\usepackage[toc,page]{appendix}

\usepackage{etoolbox}
\appto\appendix{\addtocontents{toc}{\protect\setcounter{tocdepth}{1}}}
\title[ Langlands formula and perverse sheaves]
{The Langlands formula and perverse sheaves}

\author{Mikhail Kapranov, Vadim Schechtman, Olivier Schiffmann, Jiangfan Yuan}

\begin{document}
 	
\today	
 
  \begin{abstract}
 	For a complex reductive Lie algebra $\frakg$ with Cartan subalgebra $\frakh$ and
 	Weyl group $W$ we consider the category $\Perv(\hW)$ of perverse sheaves on
 	$\hW$ smooth w.r.t. the natural stratification. We construct a category $\coxcat$
 	such that $\Perv(\hW)$ is identified with the category of  functors from $\coxcat$
 	to vector spaces. Objects of $\coxcat$ are labelled by standard parabolic
 	subalgebras in $\frakg$. It has morphisms analogous to the operations of
 	parabolic induction (Eisenstein series) and restriction (constant term) of
 	automorphic forms. In particular, the Langlands formula for the constant term
 	of an Eisenstein series  has a counterpart in the form of
 	an identity in $\coxcat$. We define $\coxcat$
 	as the  category of $W$-invariants  (in an appropriate sense) in the category $Q$
 	describing perverse sheaves on $\frakh$ smooth w.r.t. the root arrangement.
 	This matches, in an interesting way, the definition of $\hW$ itself as the
 	spectrum of the algebra of  $W$-invariants.
 \end{abstract}
 
\maketitle
 
\addtocounter{section}{-1}
 
\tableofcontents 
 
\vfill\eject

\section{Introduction}
	
\subsection{Overview of results} \label{subsec:sum-res}

The Langlands formula for the constant term of a (pseudo) Eisenstein series
\cite{Langlands-Eis, Langlands-LNM} \cite{Moeglin-W} is one of the cornerstones
of the theory of automorphic forms. It can be seen as a highly nontrivial extension to the level of functions
of the classical Mackey theorem which is a categorical result
 describing the restriction to a subgroup of a group representation
induced from another subgroup, see e.g.,  \cite{kirillov}. Both results involve summation over double cosets.  

\vskip .2cm

In this paper we relate the Langlands formula to the problem of describing perverse sheaves on 
the quotient $\hW$, where $\frakh$ is the Cartan subalgebra in (the Lie algebra $\frakg$ of)
 a complex reductive algebraic group $G$
and $W$ is the Weyl group.  More precisely, we fix  a field $\k$ of characteristic $0$ and
build a $\k$-linear category $\coxcat = \coxcat_G$ whose objects are
standard parabolics in $G$ and morphisms mimic the properties of the induction (Eisenstein series)
and restriction (constant term) operations, with a version of the Langlands formula serving as the
commutation rule between the two types of operations. We call $\coxcat$ the {\em perverse Coxeter
category} ({\em P-Coxeter category}, for short). 

Our first main result, Theorem \ref{thm:Perv-hW},  is that
\be\label{eq:Perv-hW-intro}
\Perv(\hW) \= \Fun(\coxcat, \Vect_\k). 
\ee
Here $\Perv(\hW) $ is the category of perverse sheaves of $\k$-vector spaces on $\hW$ smooth
with respect to the standard complex stratification $\calS_\frakh$, see \S \ref{subsec:rest-root} below. 
	\vskip .2cm
	
Recall that the Langlands formula involves intertwining operators which a priori satisfy the braid relations
but in the automorphic situation descend to an action of the Weyl group(oid) due to functional equations
for L-functions and Eisenstein series. In our definition of $\coxcat$ 
we retain only the braid relations.  The latter have a clear
topological meaning: they describe the fundamental
groupoids of various complex strata in $\hW$. So our $\coxcat$ can be seen as interpreting
``the Langlands formula without functional equations''. 

\vskip .2cm

In the case $G=GL_n$ we have $W=\frakS_n$, the symmetric group. We relate our category $\coxcat_{GL_n}$ to the universal braided monoidal
category $\frakB$ governing graded bialgebras $A=\bigoplus_{m\geq 0} A_m$ in braided $\k$-linear
 monoidal categories with $A_0=\1$ being the unit object  
\cite{KS-Prob}. That is, $\frakB$ is generated, as a braided category. by formal symbols $\bfa_n$, $n\geq 0$, $\bfa_0=\1$,
the component of the {\em universal graded bialgebra}
so that a general object is a product $\bfa_\bfn = \bfa_{n_1}\otimes\cdots\otimes\bfa_{n_k}$ for all possible
sequences $\bfn=(n_1,\cdots, n_k)$.  The space  $\Hom_\frakB (\bfa_\bfn, \bfa_\bfm)$
of all universal maps $\bfa_\bfn \to \bfa_\bfm$
obtained by composing multiplications, comultiplications and braidings, such compositions
being subject only to the axioms of a braided bialgebra and nothing else. The category $\frakB$ splits
into blocks $\frakB_n$  formed by objects $\bfa_\bfn$ with $\sum n_i=n$. Such objects
are in bijection with standard parabolics  $\frakp_\bfn$ in $GL_n$ and we show (Proposition \ref{prop:Inv*})  that 
$\bfa_\bfn\mapsto\frakp_\bfn$  defines a functor $\frakB_n\to\coxcat_{GL_n}$ which we conjecture to be
an equivalence. 
This interpretation reflects the phenomenon that the Eisenstein series and constant term 
operations for unramified automorphic forms on the groups $GL_n$ can be intepreted in terms of
 Hall (bi)algebras \cite {Kapr-Schiffman-Vasserot}.

 \vskip .2cm
 
 Since our motivation comes from representation theory, we work in the setup of reductive
 groups but our constructions (such as the category $\coxcat$) and results
 (such as Theorem  \ref{thm:Perv-hW}) can be formulated and proved (in the same way)
  in the slightly more general context of  arbitrary finite Coxeter groups.

\subsection{ The  perverse Coxeter category as the category of invariants} \label{subsec:PCox=Inv} 
The category $\coxcat$ is constructed starting
from another category $\calQ = \calQ_G$  which we call the {\em double incidence category}.
It goes back to \cite{K-S-realhyper}, whose main result can be formulated as
\be\label{eq:perv-h-intro}
\Perv(\frakh) \= \Fun(\calQ,\Vect_\k)
\ee
with $\Perv(\frakh)$ being the category of perverse sheaves on $\frakh$ smooth with respect
to the stratification by intersections of the root hyperplanes. Objects of $\calQ$ are all parabolic subalgebras
of $\frakg$ containing $\frakh$, i.e., they are labelled by the faces of the Coxeter complex of $G$,
see \S \ref{subsec:double-inc-cat}. 
The group $W$ acts on $\calQ$, and $\Ob(\coxcat)$ is identified with $\Ob(\calQ)/W$. Further,
\be\label{eq:cox-QW-intro}
\coxcat = \calQ^W
\ee
is obtained from $\calQ$ by an extension to linear categories of the familiar procedure
of passing to the subalgebra of invariants, see \S \ref{subsec:orb-cat-II}. For example, if we convert
$\coxcat$ and $\calQ$ into associative algebras $A(\coxcat)$ and $A(\calQ)$ (by taking endomorphisms
of the formal direct sum of all objects), then $A(\coxcat)=A(\calQ)^W$ is the usual invariant subalgebra. 
Among the relations in $\calQ$ we can distinguish a precursor of the Langlands formula for $\coxcat$, see
Proposition \ref{prop:proto-Lang}.

A general property of categories of invariants (Proposition \ref{prop:inv=kar} (b)) 
realizes $\calQ^W$ inside the Karoubi completion of  the semidirect product   $\calQ\rtimes W$. 
 The category $\calQ\rtimes W$ which describes $W$-equivariant perverse sheaves on $\frakh$
 has, in the case $G=GL_n$,  $W=\frakS_n$,  an appealing algebraic interpretation based on a version of the
 concept of a graded bialgebra, namely that of a  {\em set bialgebra}, or 
 a {\em Hopf monoid in species}  in the terminology of \cite{aguiar-species-book, aguiar:species-Hopf}. 
 Our second main result, Theorem \ref{thm:GL},
 realizes $\calQ_{GL_n}\rtimes \frakS_n$  as the $n$-th block of the braided monoidal category 
 generated by a universal set bialgebra. This corresponds to the fact that the exchange relation
 (compatibility of multiplication and comultiplication) for set bialgebras contains
 only one term  \eqref{Eq:setbialg} which matches the one-term ``proto-Langlands formula''  of Proposition
 \ref{prop:proto-Lang}.

\subsection{Invariants in commutative and noncommutative algebra}
The identifications \eqref{eq:Perv-hW-intro} \eqref{eq:perv-h-intro} and \eqref{eq:cox-QW-intro} 
mean that perverse sheaves on $\hW$  are identified with modules over  $W$-invariant subalgebra $A(\calQ)^W$
in the associative algebra $A(\calQ)$ whose modules are identified with perverse sheaves on $\frakh$. 
Note that quasicoherent sheaves on $\hW$, as on any quotient of an affine variety
by a finite group,  are identified with modules over $\C[\frakh]^W$,
the $W$-invariant subalgebra of the commutative algebra $\C[\frakh]$  governing
quasicoherent sheaves on $\frakh$. This match between two appearances of invariant theory
does not seem obvious on general grounds and deserves further study. 

\subsection {Relation to other work and futute directions} The earlier paper  \cite{K-S-hW} gave a more redundant description of 
$\Perv(\hW)$ by identifying it with   the category of  so-called
{\em mixed Bruhat sheaves} for $G$ (it depends only on $W$ and the set $S$ of its
Coxeter generators). That description proceeded by refining the complex stratification $\calS_\frakh$ of
$\frakh = \frakh_\R + i \frakh_\R$
to a real cell decomposition $\calS_{\frakh_\R} \times i\calS_{\frakh_\R}$, where $\calS_{\frakh_\R}$ is
the cell decomposition of the real part $\frakh_\R$ into faces of the root arrangement (the Coxeter complex). 
More precisely, mixed Bruhat sheaves are diagrams labelled by cells of the quotient cell decomposition
$W\backslash (\calS_{\frakh_\R} \times i\calS_{\frakh_\R})$ of $\hW$. 
At the same time, using $\calS_{\frakh_\R} \times i\calS_{\frakh_\R}$  itself leads to another,
more redundant than that of  \cite{K-S-realhyper},   description of $\Perv(\frakh)$, 
given in \cite {K-S-realhyper2}, 

Mixed Bruhat sheaves can be viewed as functors $\mathcal{MBS}_G\to\Vect_\k$ from some universal category
$\mathcal{MBS}_G$ given by generators and relations. 
In a series of unpublished notes  \cite{Tao1} - \cite{Tao3}  J. Tao established an equivalence of
  $\mathcal{MBS}_G$ with a more economic category
similar to our category $\coxcat$ but given explicitly by generators and relations. These relations 
(which include a version of the Langlands formula) do hold in
$\coxcat$ and we conjecture that our $\coxcat$ is equivalent to Tao's category and thus to $\mathcal{MBS}_G$. 
The statement in \S  \ref{subsec:sum-res} that $\coxcat_{GL_n}$   is equivalent to $\frakB_n$ is a particular case of this conjecture. 
We plan to address this in a future paper by detailed comparison of all the descriptions of $\Perv(\frakh)$ and
$\Perv(\hW)$ listed above.

\subsection{The structure of the paper} 
Section \ref{sec:notation} sets general conventions and fixes notation related to root systems.

Section \ref{sec:coxcom&hW} descibes the general properties of the root hyperplane arrangement in $\frakh$,
the fundamental groupoids of its complex strata (analogs of pure braid groupoids) and their descent to $\hW$ (analogs of
braid groupoids). We emphasize the known but little discussed subtlety in that the induced arrangements in
the flats need not, in general, be root arrangements of other root systems. 

In Section \ref{sec:examples} we work out the examples of 
types $A_2,B_2$ and $D_4$. 

In Section \ref{sec:percox} we recall the double incidence category $\calQ$ from  \cite{K-S-realhyper}
and analyze it  from the point of view of the Langlands formula, interpreting various relations as
``toy one-term prototypes'' of various identities involving Eisenstein series and constant terms.
Then we define the P-Coxeter category $\coxcat$ as the category of invariants and deduce
the full-fledged identities holding in it. 

Section \ref{sec:GLn} is devoted to the case $G=GL_n$, where $W=S_n$ is the symmetric group.
 After combinatorial preliminaries, we recall from  \cite{KS-Prob}
 the concept of graded bialgebras  and the
PROB (universal braided category) $\frakB$ governing them.  In addition, we consider related algebraic
structures which we call {\em set bialgebras}, by mild generalization (to the case of an arbitrary braided category)
the concept of {\em Hopf monoids in species} from \cite{aguiar-species-book, aguiar:species-Hopf}.
As for graded bialgebras, there is the PROB  $\frakSB$ governing them and we prove (Theorem \ref{thm:GL})
that $\frakSB_n$, its degree $n$ block, is equivalent to the semidirect product $\calQ_{GL_n}\rtimes S_n$. We then construct a braided monoidal functor $\frakB\to\frakSB$ which we conjecture to be
an equivalence. 

In Section \ref{sec:Perv-hW} we prove our main result, the identification \eqref{eq:Perv-hW-intro}, in
Theorem \ref{thm:Perv-hW}. The proof is based on analyzing the projection $p: [\hW] \to \hW$ of the
stacky quotient $[\hW]$ to the geometric quotient $\hW$. Perverse sheaves on $[\hW]$ are $W$-equivariant
perverse sheaves on $\frakh$ and are identified by \cite{K-S-realhyper} with functors $\calQ\rtimes W\to\Vect_\k$. 
We then proceed to use the adjoint pair $(p^{-1}, p_*)$ of the appropriate pullback and pushforward functors
to identify $\coxcat$  with  $\Perv(\hW)$ inside $\Perv[\hW]$ using dimensional d\'evissage. 

Finally, the Appendix collects general statements about categories of invariants and descent  
for sheaves and modules under quotienting of a space by a finite group or passing from a ring to the
invariant subring.

\subsection {Acknowledgements}

The work of M.K. was supported by the World Premier International Research Center Initiative (WPI
Initiative), MEXT, Japan and  by the JSPS KAKENHI grant 20H01794.
It was also supported by the US National Science Foundation under Grant No. DMS-1928930 and by the Alfred P. Sloan Foundation under grant G-2021-16778, while the author was in residence at the Simons Laufer Mathematical Sciences Institute (formerly MSRI) in Berkeley, California, during the Spring 2024 semester. He would like to thank  D. Kaledin and A. Polishchuk
for useful discussions  regarding this project. 

The work of O.S. was supported by PNRR Grant CF 44/14.11.2022. He thanks the Simion Stoilow Institute of Mathematics for a great research environment, as well as L. Hennecart, B. Hennion, C. Hohlweg and I. Marin for useful discussions.

We are grateful to J. Tao for sending us his unpublished notes \cite{Tao1} - \cite{Tao3}.

\vfill\eject

\section{Notations and conventions}\label{sec:notation}

 \begin{itemize}
 
  \item By $\N$ we denote the set $\{0,1,2, \cdots\}$ of non-negative integers. 
 
 \item All topological spaces in this paper are assumed to be locally homeomorphic to an open
  subset in
 a finite CW-complex. 
 
 \item $\k$ will denote a base field of characteristic $0$ (not necessarily equal to $\C$). 
 By $\Vect_\k$ we denote the category of $\k$-vector spaces. 
 
 \item  ``Braided category'' will stand for ``$\k$-linear braided monoidal category''.

\end{itemize} 

We further  fix the following:
\begin{itemize}
\item $G$ a complex reductive algebraic group,
$T$ a maximal torus of $G$, 
$r$ the semisimple rank of $G$.

\item $W=N(T)/Z(T)$ the associated Weyl group.

\item  $\frakg=\Lie(G), \frakh=\Lie(T)$, $\frakz\subset\frakh$ the center of $\frakg$. 

\item $\Delta, \Delta^\vee$ the corresponding root and coroot systems. For $\alpha\in\Delta$ we denote
by $\alpha^\vee\in\Delta^\vee$ the corresponding coroot. 

\item $\frakh_\R, \frakh^\vee_\R$ real forms of $\frakh, \frakh^\vee$ such that $\frakh_\R \supset \Span(\Delta^\vee), \frakh_\R^\vee \supset \Span(\Delta)$.

\item $\Pi=\{\alpha_1, \ldots, \alpha_r\} \subset \Delta$ a base, $B\subset G$ and $\frakb\subset \frakg$ the corresponding Borel subgroup and subalgebra.

\item $\Delta_+\subset\Delta$ the set of positive roots, i.e., roots of $\frakb$. 

\item $\mathcal{P}_{G,T}$ the set of parabolic subalgebras of $\frakg$ containing $\frakh$ (equivalently, the set of parabolic subgroups of $G$ containing $T$).

\item $\mathcal{P}_{G,B}$ the set of parabolic subalgebras containing $\frakb$ (equivalently, the set of parabolic subgroups of $G$ containing $B$).

	\item $\mathcal{L}_{G,T}$ the set of Levi subalgebras of $\frakg$ containing $\frakh$; 
there is a canonical projection $\pi: \mathcal{P}_{G,T} \to \mathcal{L}_{G,T}$ 
assigning to each $\frakp$ the unique $\frakh$-stable lift of $\frakp/\frakn(\frakp)$, 
where $\frakn(\frakp)$ is the nilpotent radical of $\frakp$.

\item For any $\frakp \in \calP_{G,T}$ we denote by $\Delta_\frakp \subset \Delta$ 
the set of roots occurring as $\frakh$-weights on $\frakp$. Thus $\Delta_+=\Delta_\frakb$. 

\item Likewise, for any $\frakl\in\calL_{G,T}$ we denote by $\Delta_\frakl\subset\Delta$ the  set of
roots occurring as $\frakh$-weights on  $\frakl$. This is also the root system of $\frakl$ as a reductive Lie algebra.

\item for any $I \subset \Pi$, $\frakp_I$ (or $P_I$) is the corresponding element of $\mathcal{P}_{G,B}$, so that 
$$\mathcal{P}_{G,B}=\{\frakp_I\;|\; I \subset \Pi\};$$
we also set $\frakl_I=\pi(\frakp_I)$. Elements of $\calP_{G,B}$ will be called 
{\em standard parabolic subalgebras}. There are $2^r$ of them. 

\item for any $I \subset \Pi$, $W_I$ the parabolic subgroup of $W$ generated by $s_{i}$ for $\alpha_i \in I$. 
It is also the Weyl group of (the reductive subgroup $L_I\subset G$ associated to) $\frakl_I$. 

\item
More generally, we denote by $W_\frakl \subset W$ the subgroup generated by all reflections
 $s_\alpha$ for $\alpha$ a root of $\frakl$. It is also the Weyl group of $\frakl$. 
 
 \end{itemize}

\vspace{.15in}

\vfill\eject

\section{ The Coxeter complex, the root arrangement and $\hW$}\label{sec:coxcom&hW}

\vspace{.1in}

\subsection{Combinatorics of parabolic subalgebras and the Coxeter complex} \label{subsec:comb-cox}

There is a natural action of $W$ on $\mathcal{P}_{G,T}$; every orbit meets $\mathcal{P}_{G,B}$ in exactly one point, i.e. $\mathcal{P}_{G,T}=\bigsqcup_I W \cdot \frakp_I$. Moreover, $\Stab_W(\frakp_I)=W_I$, yielding an identification
\[
\mathcal{P}_{G,T}=\bigsqcup_I \, W/W_I.
\]
If $\frakp \in W \cdot \frakp_I$ then we say that $\frakp$ is of \textit{type $I$}. Reverse inclusion defines a partial order $\preceq$ on $\calP_{G,T}$ : we write $\frakp \preceq \frakp'$ if $\frakp \supset \frakp'$. We put
\[
{C_\frakp}:=\{x \in \frakh_\R\;|\;  \alpha(x) \geq 0,\, \forall\,\alpha \in \Delta_\frakp \}.
\]
This is a closed polyhedral cone in $\frakh_\R$. Its relative interior $C_\frakp^\circ$ is described as
\[
C_\frakp^\circ \,=\,\bigl\{ x \in C_\frakp\;|\;  \alpha(x) > 0,\, \, \forall\,\alpha \in \Delta_\frakp \text{ s.t. } 
-\alpha \notin \Delta_\frakp\bigr\}. 
\]
In particular, for $\frakp=\frakb$ the standard Borel, we denote $C_+^\circ = C_\frakb^\circ$,
$C_+ = C_\frakb$ and call them the open and closed dominant Weyl chambers. 
We then have
\[
C_\frakp=\bigsqcup_{\frakp' \supset \frakp} C^\circ_{\frakp'}
\]
and a further disjoint decomposition
\be\label{eq:cox-complex}
\frakh_{\R}=\bigsqcup_{\frakp \in \calP_{G,T}} C^\circ_\frakp
\ee
of $\frakh$ itself into cones.
We will refer to the decomposition \eqref{eq:cox-complex} as the {\em Coxeter complex} and denote
it $\Cox$. Thus the minimal cone in $\Cox$ is $C_\frakg = \frakz_\R$.

\vspace{.15in}

\subsection{Combinatorics of Levi subalgebras and the root arrangement} 

To any root $\alpha \in \Delta$ is associated the {\em root hyperplane}
\[
H_{\alpha}=\R \alpha^\perp \subset \frakh_\R.
\]
Thus $H_\alpha=H_{-\alpha}$, so one can think of the root hyperplanes as labelled by $\Delta_+$. 
The arrangement of hyperplanes $\calH = \{H_\alpha|\, \alpha\in\Delta_+\}$, is known as the
{\em root arrangement}. The cones $C_\frakp^\circ$, $\frakp\in\calP_{G,B}$, are precisely the
{\em faces} of $\calH$, i.e., relatively open polyhedral cones into which $\calH$ partitions $\frakh$. 

Further, $\calH$ gives rise to the set of {\em flats}, i.e., linear subspaces in $\frakh_\R$ obtained as all possible intersections
of hyperplanes from $\calH$ (including  $\frakh_\R$ itself which is the intersection of the empty set of hyperplanes). 
Denote $\Fl(\calH)$ the poset of flats in $\calH$ ordered by inclusion. It is described as follows.

\begin{prop}\label{prop:flats-Hl}
Flats of $\calH$ are in order-reversing bijection with Levi subalgebras $\frakl\in\calL_{G,T}$: to each $\frakl$
we associate the flat $H_\frakl =\{x \in \frakh_\R\;|\;  \alpha(x) = 0,\, \forall\,\alpha \in \Delta_\frakl \}$. 
Thus $\Fl(\calH)\simeq \calL_{G,T}$.  In particular, the minimal flat is $H_\frakh = \frakg$. \qed
\end{prop}

The projection $\pi: \calP_{G,T} \to \calL_{G,T}$ from \S \ref{sec:notation} can be identified with
the map
\be\label{eq:pi-geometr}
\pi: \{\text{Faces of } \Cox\} \lra \Fl(\calH), \qquad F\mapsto \Span_\R (F). 
\ee

\subsection{Restriction of the root arrangement to a flat and its complexification}\label{subsec:rest-root}
For any $\frakl\in\calL_{G,T}$ we
denote by $\calH_\frakl$ the arrangement of hyperplanes in $H_\frakl$ formed by intersections
$H_\frakl \cap H_\alpha$, $H_\alpha \in\calH$, $H_\alpha \not\supset H_\frakl$. Note that
a hyperplane in $\calH_\frakl$ can sometimes be obtained as such an intersection in more than one way,
in which case we count it only once. The arrangement $\calH_\frakl$ cannot, in general,
be seen as a root arrangement for any reductive group, see \S~\ref{subsec:D4l=1} below for an example. 
But it is simplicial: each  face is a simplicial cone times $\frakz_\R$.
This is because the original arrangement $\calH$ is simplicial. As usual, open faces will be
called {\em (open) chambers} and their closures {\em closed chambers}.

We denote by $\calH_\C$ resp. $\calH_{\frakl, \C}$ the arrangement of complex hyperplanes in 
$\frakh$  resp. $H_\frakl$ obtained
by complexifying those from $\calH$ resp. $\calH_\frakl$. Thus $\calH_\C = \calH_{\frakh, \C}$. 
By the above, the union of the hyperplanes from $\calH_\frakl$ is $\bigcup_{\frakl'\supsetneq \frakl} 
H_{\frakl'}$. Removing the complexifications of these hyperplanes, we get an open subset
\[
H^\circ_{\frakl, \C} = H_{\frak l, \C} -  \bigcup_{\frakl'\supsetneq \frakl} 
H_{\frakl',\C}
\]
which we call the {\em generic stratum} of $H_{\frakl, \C}$. Its real part
\[
H^\circ_\frakl = H^\circ_{\frakl, \C} \cap H_\frakl = H_{\frak l} -  \bigcup_{\frakl'\supsetneq \frakl} 
H_{\frakl'}
\]
is the union of  open chambers in $H_\frakl$. 

The subvarieties $H^\circ_{\frakl, \C}$, $\frakl\in \calL_{G, T}$, form a complex stratification of 
$\frakh$ which we denote $\calS_\frakh$. The induced stratification of $\frakh_\R$
is essentially (up to splitting disconnected strata into connected components) the Coxeter
complex $\Cox$. 

\subsection{The fundamental groupoid of the generic stratum $H^\circ_{\frakl,\C}$} 
\label{subsec:pi-1-h}
The chambers in $H_\frakl$, i.e., the connected components of $H^\circ_{\frakl}$,
are contractible, being convex cones. Let us put one base point in each chamber
and denote  $\Pi_1(H^\circ_{\frakl,\C})$  the fundamental groupoid of $H^\circ_{\frakl,\C}$ with respect
to this set of base points. Because of contractibility of the chambers, the precise choice
of base points inside them is inessential. By \eqref{eq:pi-geometr}, 
$\Ob \Pi_1(H^\circ_{\frakl,\C})$,
i.e., the set of chambers in  $H^\circ_{\frakl}$, is identified with $\pi^{-1}(\frakl)$,
the set of parabolic algebras $\frakp\in \calP_{G,T}$ with Levi $\frakl$. 

\begin{prop}\label{prop:pi_1-linear}
$\Pi_1(H^\circ_{\frakl,\C})$ has the following description by generators and relations:
\begin{itemize}
\item[(0)] $\Ob \Pi_1(H^\circ_{\frakl,\C})=\pi^{-1}(\frakl)$.

\item[(1)]
For each pair $\frakp_1,\frakp_2\in \pi^{-1}(\frakl)$  of objects there is a generating morphism
$\tau_{\frakp_1}^{\frakp_2}: \frakp_1\to\frakp_2$, with $\tau_\frakp^\frakp=\Id$. They are
subject to the following relations:

\item[(2)]  Suppose $\frakp_1, \frakp_2, \frakp_3\in\pi^{-1}(\frakl)$ are such that there are
$x_i\in C^\circ_{\frakp_i}$ such that $x_2$ lies in the straight  line segment $[x_1, x_3]$.
Then $\tau_{\frakp_1}^{\frakp_3} = \tau_{\frakp_2}^{\frakp_3} \circ \tau_{\frakp_1}^{\frakp_2}$.
\end{itemize}

\end{prop}

\noindent{\sl Proof:} This is  \cite [Prop. 9.11]{K-S-realhyper} which is
an easy consequence of the explicit  cellular model for the homotopy 
type of $H^\circ_\frakl$ given by Salvetti \cite{Salvetti} which goes back to Deligne \cite{Deligne}
for simplicial arrangements.  

\qed

\begin{rems} The above generators and relations can be understood as follows. 

\vskip .2cm
(a) For two objects $\frakp_1, \frakp_2$, we say that $(x_1,x_2) \in C^\circ_{\frakp_1} \times C^\circ_{\frakp_2}$
is {\em generic} if the line segment $[x_1,x_2] \subset H_\frakl$ has only transverse intersections with the other cells $C^\circ_{\frakp'}$; in particular, $[x_1,x_2]$ only intersects codimension one walls. 
Let us choose a generic $(x_1, x_2)$.  The morphism $\tau_{\frakp_1}^{\frakp_2}$ is represented by a path in $H^\circ_{\frakl, \C}$
which coincides with $[x_1, x_2]$ away from small neighborhoods of the intersection points with
the walls and avoids the walls by going around them in the complex domain so that the imaginary
part of the value of the  equation of the wall has the same sign as the (real) value of this equation
on the vector $x_2-x_1$.  Thus all the $\tau_{\frakp_1}^{\frakp_2}$ are obtained by composition
from those in which  $\frakp_1, \frakp_2$ are adjacent, i.e., situated on opposite sides of a wall. 

\vskip .2cm
(b)  For a generic $(x_1, x_2)$ as above let us 
record the sequence of chambers crossed and write $G[x_1,x_2]=(\frakp_1=\frakq_1, \frakq_2, \ldots, \frakq_s=\frakp_2)$. After expressing all the $\tau_{\frakp_1}^{\frakp_2}$ through those with adjacent
$\frakp_1, \frakp_2$, the relations on the latter can be expressed as follows:~
for any two generic pairs $(x_1,x_2)$ and $(x'_1,x'_2)$ in $ C^\circ_{\frakp_1} \times C^\circ_{\frakp_2}$ we have
\begin{equation}\label{E:relations_fund_groupoid}
\tau_{\frakq_{s-1}}^{\frakp_2}  \cdots  \tau^{\frakq_2}_{\frakp_1}=\tau^{\frakp_2}_{\frakq'_{s'-1}} \cdots \tau^{\frakq'_2}_{\frakp_1} \in \Hom_{\Pi(H^\circ_{\frakl, \C})}(\frakp_1,\frakp_2)
\end{equation}
if $G[x_1,x_2]=(\frakq_1, \frakq_2, \ldots, \frakq_s), G[x'_1,x'_2]=(\frakq'_1, \frakq'_2, \ldots, \frakq'_{s'})$. 
Indeed, they should both be equal to  $\tau^{\frakp_2}_{\frakp_1}$.

\vskip .2cm

(c) Objects $\frakp$ correspond to $0$-cells (points), the elementary generators $\tau_{\frakp_1}^{\frakp_2}$ for adjacent $\frakp_1, \frakp_2$ correspond
to $1$-cells (intervals) and the relations
\eqref{E:relations_fund_groupoid} corresponds to   $2$-cells ($2s$-gons with various $s$) in Salvetti's
polyhedral model for $H_{\frakl,\C}^\circ$, see \cite{Salvetti}. 
\flushright{$\triangle$}
\end{rems}

\subsection{The quotients $\hW$ and $\hW_\R$}
The group $W$ acts on $\frakh$ and the quotient $\hW = \Spec\, \C[\frakh]^W$ is known to be isomorphic to
$\AAA^n$, $n=\dim\frakh$. The stratification $\calS_\frakh$ of $\frakh$ by open strata of the (complex)
flats of $\calH$
descends to a stratification of $\hW$ that we denote $\calS$. Thus strata of $\calS$ are labelled by the
set $W\backslash \calL_{G,T}$ of $W$-orbits on the set of Levis. We denote by $p: \frakh \to \hW$ the quotient map.

We also have the quotient $\hW_\R$ of the real part $\frakh_\R$ which we can identify  using the cell
decomposition of $\frakh_\R$ into the faces of the Coxeter complex.
The remarks at the beginning of \S \ref{subsec:comb-cox} imply that 
the (topological) quotient $\hW_\R$ is identified with $C_+$, the closed dominant Weyl chamber.
In particular, each face of $\calH$ is taken by $W$-action into a  face $C_I = C_{\frakp_I}\subset C_+$
for a unique subset $I\subset \Pi$ of simple roots.  
Let $2^\Pi$ be the set of such subsets. 

As any flat of $\calH$ is $\R$-linearly spanned by some face, we see that the set of flats $H_I = H_{\frakp_I}$ associated to  standard parabolics $\frakp_I$, $I\in 2^\Pi$, 
generates the set $\Fl(\calH)$ under $W$-action, in other words we have a surjective map
\be
\varpi: 2^\Pi \lra W\backslash \calL_{G,T}.
\ee
We denote by $S_I\subset \hW$ the stratum corresponding to $\varpi(I)$. 
Two standard parabolics $\frakp_I, \frakp_J$, $I,J\in 2^\Pi$, are called {\em associated}, if  
$S_I=S_J$, i.e.,  $\varpi(I)=\varpi(J)$.

\begin{ex}\label{ex:gl_n-strata}
Let $\frakg = \gl_n$, so $\frakh = \C^n$ and $W=S_n$. The quotient $\hW\simeq \C^n$ is identified with the space of monic polynomials
\[
f(x) = x^n+ a_1 x^{n-1} + \cdots + a_n, \quad a_i\in\C. 
\]
The strata in $\hW$ are labelled by unordered paritions $\alpha = \{\alpha_1, \cdots, \alpha_p\}$ of $n$. The stratum
$S_\alpha$ corrersponding to $\alpha$, consists of polynomials whose multiplicities of roots are given by the $\alpha_i$. 
In particular, the generic stratum $S_{1,\cdots ,1}$ consists of polynomials with $n$ distinct roots and its complement,
the  discriminantal hypersurface, is the closure of the stratum $S_{2,1, \cdots, 1}$. So unlike the case of the stratification
$\calS_\frakh$ of $\frakh$, the closures of the strata in $\calS$ can be singular. 

Standard parabolics in $\frakg$ are labelled by {\em compositions} or {\em ordered partitions}  
of $n$, i.e., sequences  $\bfn = (n_1,\cdots, n_p)$
s.t.  $\sum n_i=n$. That is, the parabolic  $\frakp_\bfn$ associated to $\bfn$ consists of block-upper-triangular matrices
with diagonal blocks of sizes $n_i$.  The map $\varpi$ sends an ordered partition into the corresponding unordered one, so that two standard
parabolics $\frakp_\bfn, \frakp_\bfm$ are associated iff $\bfm$ and $\bfn$ differ by a permutation. 
\flushright{$\triangle$}
\end{ex}

We now analyze the relation between strata in $\frakh$ and $\hW$ in more detail. 
For any $0 \leq l \leq r$,
we set 
\[
\calL_{G,T}^{(l)}=\bigcup_{I, |I|=l} W \cdot \frakl_I,
\quad \mathcal{H}^{(l)}=\{H_\frakl\;|\; \frakl \in \calL^{(l)}_{G,T}\}=\{H_{\frakl}\;|\; \mathrm{codim}_{\frakh_\R}(H_\frakl)=l\},
\quad \calH^{(1)}=\calH. 
\]
We will sometimes say that a Levi subalgebra $\frakl$ is of \textit{length} $l$ if $\frakl \in \calL^{(l)}_{G,T}$.

The group $W$ acts on $\calH$, preserving each $\calH^{(l)}$. The stabilizers $\Stab_W(\frakl)$ for $\frakl \in \calL_{G,T}$ have been determined by Howlett  \cite{Howlett}. For $\frakl=\frakl_I$ it contains $W_I \times W_{I^\perp}$ but it can sometimes be bigger, see e.g. the example of $D_4$ in Section~\ref{subsec:D4l=1}.

\begin{prop}\label{prop:W-l}
(a) $W_\frakl$, the Weyl group of $\frakl$ as a reductive Lie algebra, see \S \ref{sec:notation}, is a normal
subgroup in $\Stab_W(\frakl)$. 

\vskip .2cm 
(b) The quotient group $W(\frakl):= \Stab_W(\frakl)/W_\frakl$ acts  freely on $H^\circ _{\frakl, \C}$ and
$W_\frakl \backslash H^\circ _{\frakl, \C} = p(H^\circ _{\frakl, \C})$ is the $\calS$-stratum in $\hW$
corresponding to $H^\circ _{\frakl, \C}$. 

\end{prop}

\noindent{\sl Proof:} This is probably well-known, but we insert a proof for the reader's convenience. We may assume that $\frakl=\frakl_I$ for a subset $I \subset \Pi$. Let us denote by $\rho: \Stab_W(\frakl_I) \to \text{Aut}(H_{\frakl_I})$ the canonical map. We claim that $\Ker(\rho)=W_I$. Indeed, an element $w \in W$ belongs to $\Stab_W(H_{\frakl_I})$ if and only if it permutes the chambers of $H_{\frakl_I}$; on the other hand, the stabilizer of $C_{\frakp_I}$ in $W$ is well-known to be $W_I$. This shows (a). We now prove (b).
Let $x \in H^\circ_{\frakl}$ and assume that $w \in \Stab_W(\frakl)$ satisfies $w(x)=x$. Then $w$ stabilizes the chamber of $H_{\frakl}$ containing $x$. Conjugating by a suitable element of $W$ to bring that chamber into $C_+$ and arguing as in (a) we deduce that $w$ acts as $Id$ on $H_{\frakl}$, i.e. that $w \in W_I$, as wanted. 
\qed

\begin{ex}
For $l=1$, we have $\calH^{(1)}=\Delta/{\pm}$, hence $\calH^{(1)}/W$ has one or two elements depending on whether $\frakg$ is simply laced or not; in particular, the action of $W$ on $\calH^{(l)}$ is in general not transitive. For $\frakl_{\pm \alpha}:= \frakg_{\alpha} \oplus \frakh \oplus \frakg_{-\alpha}\in\calL^{(1)}_{G,T}$ we have $W(\frakl_{\pm \alpha})=\Stab_W(H_\alpha)/\langle s_\alpha\rangle$.
\flushright{$\triangle$}
\end{ex}

\subsection{The fundamental groupoid of a stratum in $\hW$} 
Let $\frakl\in\calL_{G,T}$, so $S_\frakl = W_\frakl \backslash H^\circ _{\frakl, \C}$ is the stratum
in $\hW$ corresponding to $\frakl$, see Proposition \ref{prop:W-l}. As in \S 
\ref{subsec:pi-1-h}, let $\Pi_1(H^\circ_{\frakl, \C})$ be the fundamental
groupoid of $H^\circ_{\frakl, \C}$ with the set of base points, one in each chamber.
We recall that the set of chambers is identified with $\pi^{-1}(\frakl)$ where
$\pi: \calP_{G,T}\to\calL_{G,T}$ is the natural projection. 
The set  of base points can be taken to be $W(\frakl)$-stable, so the group $W(\frakl)$ acts
on $\Pi_1(H^\circ_{\frakl, \C})$ in the strict sense, the action on the set of objects
being free. We define the {\em braid groupoid of} $\frakl$ to be the coinvariant
orbit category as defined in \S \ref{subsec:coinv-orbits}:
\begin{equation}\label{def:braidlevi}
\Br(\frakl,W):=\Pi_1(H_{\frakl,\C}^\circ)_{W(\frakl)}.
\end{equation}
Note that this only depends on the conjugacy class of $\frakl$, i.e. on the $W$-orbit of $\frakl$.
Observe that two parabolic subalgebras $\frakp,\frakp'$ corresponding to chambers of $H_{\frakl}$ are in the same $W(\frakl)$-orbit if and only if they are of the same type (i.e. if and only if they are $W$-conjugate); indeed, if $w \cdot\frakp=\frakp'$ then $w(C_\frakp)=C_{\frakp'}$ hence $w \in \Stab_W(H_\frakl)$ and the corresponding set of types forms a subset of $\Pi$.
This means that for each $\frakl$ the set $\Ob \Br(\frakl,W)$, i.e. the set of $W(\frakl)$-orbits on 
$\Ob \Pi_1(H_{\frakl,\C}^\circ) \simeq \pi^{-1}(\frakl)\subset \calP_{G,T}$ is identified with a  certain
class $A_\frakl$ of associated standard parabolics or, equivalently, with the
set $F_\frakl$ of corresponding faces of the closed dominant chamber $C_+ = C_\frakb$. 
Recalling that $C_+$ is identified with
$\hW_\R$, 
Example
\ref{ex:fund-grpd-coinv} implies the following.

\begin{prop}

(a) The union of the interiors of the faces from $F_\frakl$ is equal 
to the intersection $S_\frakl\cap(\hW_\R)$.

\vskip .1cm
(b) The groupoid $\Br(\frakl, W)$ is isomorphic to the fundamental groupoid of $S_\frakl$
with the set of based points obtained by taking one point inside each face 
from $F_\frakl$. 
\end{prop}

\begin{rems}
(a) By Proposition \ref{prop:pi_1-linear}, 
$\Br(\frakl,W)$ is equivalent to the groupoid whose set of objects is $W(\frakl) \backslash \pi^{-1}(\frakl)$, generating morphisms are given by $W(\frakl)$-orbits of maps $\tau_{\frakp_1}^{\frakp_2}$ as above, and relations are given by  relations (3) of the proposition or, equivalently, by
\eqref{E:relations_fund_groupoid}.

\medskip

(b) When $\frakl=\frakh$ we recover the well-known presentation of the classical braid  group
(groupoid with one object)
associated to the Weyl group $W$, i.e. $\Br(\frakh,W)\simeq \Br_W$.  In that case $W$ acts simply transitively on the set of chambers 
(so there only one orbit) and in the usual notation $\sigma_i=\tau_{w_1}^{w_2}$ whenever $s_i=w_2w_1^{-1}$.
\flushright{$\triangle$}
\end{rems}

\vfill\eject

\section{Examples : types $A_2,B_2$ and $D_4$}\label{sec:examples}

\vspace{.1in}

We now illustrate all the above notions in the following cases : $\frakg$ is a simple Lie algebra of type $A_2,B_2$ or $D_4$ and $\frakl=\frakl_\alpha$ corresponds to a simple root; $\frakg$ is of type $D_4$ and $\frakl$ corresponds to subset consisting of a pair of simple roots. 

\vspace{.2in}

\subsection{{Type $A_2$}} 
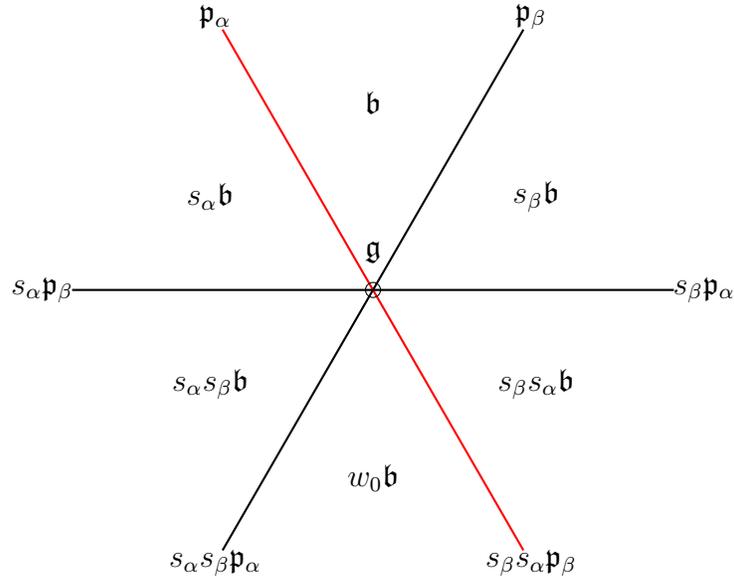
\begin{figure}[h]
\begin{tikzpicture} 
\draw[thick] (60:4)--(240:4);
\draw[thick,red] (120:4)--(300:4);
\draw[thick] (-4,0)--(4,0);
\draw (0,0) circle (0.1);
\node(a)at(0,2.5){$\frakb$};
\node(a)at(150:2.5){$s_\alpha\frakb$};
\node(a)at(30:2.5){$s_\beta\frakb$};
\node(a)at(-150:2.5){$s_\alpha s_\beta\frakb$};
\node(a)at(-30:2.5){$s_\beta s_\alpha\frakb$};
\node(a)at(0,-2.5){$w_0\frakb$};
\node(a1)at(4.4,0){$s_\beta \frakp_\alpha$};
\node(b1)at(-4.4,0){$s_\alpha\frakp_\beta$};
\node(c1)at(60:4.2){$\frakp_\beta$};
\node(d1)at(60:-4.2){$s_\alpha s_\beta \frakp_\alpha$};
\node(e1)at(120:4.2){$\frakp_\alpha$};
\node(f1)at(120:-4.2){$s_\beta s_\alpha \frakp_\beta$};
\node(o)at(0,0.5){$\frakg$};
\end{tikzpicture}
\caption{Type $A_2$}
\label{FigureA_2}
\end{figure}

Write $\Pi=\{\alpha,\beta\}$; the geometric realization of the Coxeter complex, with the parabolic subalgebras associated to each facet, is given in Figure~\ref{FigureA_2}.

There are $6$ chambers, $6$ facets of codimension $1$ (one for each ray), and one of codimension $2$ (the origin). The parabolics associated to the chambers are the Borel subalgebras; the parabolic associated to the cell at the origin is $\frakg$ itself. 
There are $3$ Levi subalgebras of length $1$, forming a single $W$-orbit. Let us consider $\frakl=\frakl_{\alpha}$. The subspace $H_{\frakl}$ corresponds to the red line in the above picture, which consists of three cells, namely $C^\circ_{\frakp_\alpha}, C^\circ_{s_\beta s_\alpha\frakp_{\beta}}$ and $C^\circ_{\frakg}$. We have
$\Stab_W(\frakl)=\{1,s_\alpha\}=W_{\frakl}$ hence $W(\frakl)=\{1\}$. It follows that $\Br(\frakl,W)$ is the fundamental groupoid of $H^\circ_{\frakl,\C}=\C^*$, which is equivalent to $\Z$. Note that there is no element of $W$ which maps $C^\circ_{\frakp_\alpha}$ to $C^\circ_{s_\beta s_\alpha\frakp_{\beta}}$, i.e. they belong to different conjugacy classes of parabolics. 
The combinatorial presentation of $\Br(\frakl,W)$ is
$$\xymatrix{\{\alpha\} \ar@/^/[r]^-{\{\alpha,\beta\}} &\{\beta\} \ar@/^/[l]^-{\{\alpha,\beta\}} }.$$

\vspace{.2in}

\subsection{{Type $B_2$}} write $\Pi=\{\alpha,\beta\}$ with $\alpha$ being the short root; the geometric realization of the Coxeter complex is as in Figure~\ref{FigureB_2} below.
\begin{figure}
\begin{tikzpicture}
	\draw[thick] (45:4)--(45:-4);
	\draw[thick] (-45:4)--(-45:-4);
	\draw[thick] (-4,0)--(4,0);
	\draw[thick,red] (0,4)--(0,-4);
	\draw (0,0) circle (0.1);
	\node(a)at(1,2.5){$\frakb$};
	\node(a)at(135:4.3){$s_\alpha\frakp_\beta$};
	\node(a)at(135:-4.3){$s_\beta s_\alpha\frakp_\beta$};
	\node(a)at(-90:4.2){$w_0\frakp_\alpha$};
	\node(a)at(45:-4.3){$w_0\frakp_\beta$};
	\node(a)at(22:2.5){$s_\beta\frakb$};
	\node(a)at(157:2.5){$s_\alpha s_\beta\frakb$};
	\node(a)at(112:2.5){$s_\alpha\frakb$};
	\node(a)at(112:-2.5){$s_\beta s_\alpha s_\beta \frakb$};
	\node(a)at(-22:2.5){$s_\beta s_\alpha\frakb$};
	\node(a)at(22:-2.5){$s_\alpha s_\beta s_\alpha\frakb$};
	\node(a)at(67:-2.5){$w_0\frakb$};
	\node(b1)at(4.4,0){$s_\beta\frakp_\alpha$};
	\node(b1)at(-4.6,0){$s_\alpha s_\beta\frakp_\alpha$};
	\node(c1)at(45:4.3){$\frakp_\beta$};
	\node(e1)at(90:4.2){$\frakp_\alpha$};
	\node(o)at(0.2,0.5){$\frakg$};
\end{tikzpicture}
\caption{Type $B_2$}
\label{FigureB_2}
\end{figure}
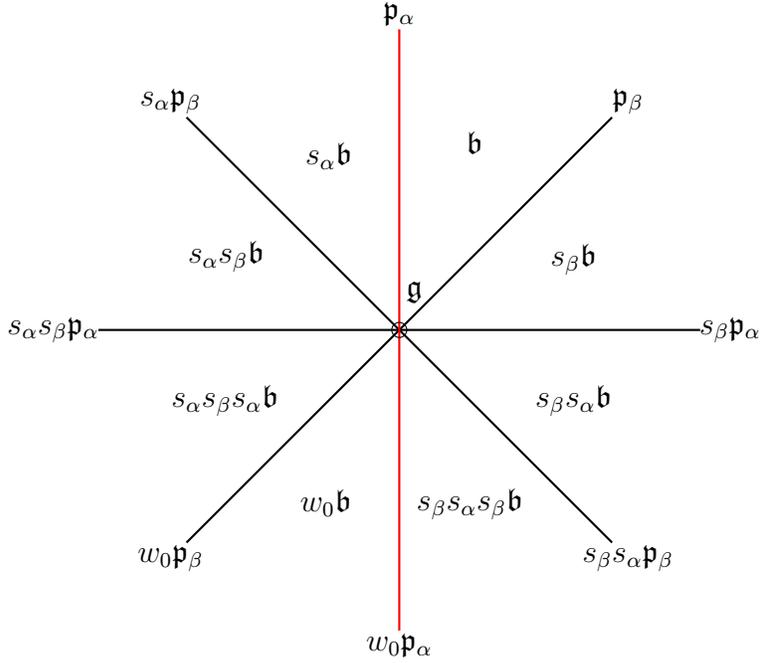
There are now two $W$-orbits of Levi subalgebras of length one : they can both be treated similarly. Let us consider $\frakl=\frakl_{\alpha}$. The subspace $H_{\frakl}$ corresponds to the red line in Figure~\ref{FigureB_2}, which consists of three cells, namely $C^\circ_{\frakp_\alpha}, w_0C^\circ_{\frakp_\alpha}$ and $C^\circ_{\frakg}$. We have $\Stab_W(\frakl)=\langle s_\alpha, w_0\rangle \simeq (\Z/2\Z)^2$ and $W_{\frakl}=\{1,s_\alpha\}$ hence $W(\frakl)\simeq \Z/2\Z$. It follows that $\Br(\frakl,W)$ is the fundamental groupoid of $H^\circ_{\frakl,\C}/(\Z/2\Z)=\C^*/(\Z/2\Z)$, which is again equivalent to $\Z$. 
The combinatorial presentation of $\Br(\frakl,W)$ is

\centerline{
\begin{tikzpicture}
	\node(a)at(0,0){$\{\alpha\}$};
	\node(a)at(-2.3,0){${\tiny{\{\alpha,\beta\}}}$};
	\draw[thick,->](-0.4,0.1)..controls(-2,0.9) and (-2,-0.9)..(-0.4,-0.1);
\end{tikzpicture}
}

\vspace{.2in}

\subsection{{Type $D_4$ and $l=1$}}\label{subsec:D4l=1}
 We write $\Pi=\{\alpha_1,\alpha_2,\alpha_3,\alpha_*\}$, where $\alpha_*$ is the trivalent vertex of the Dynkin diagram. The Coxeter complex has $192$ chambers (parametrized by the Weyl group $W=(\Z/2\Z)^3 \times \mathfrak{S}_4$, where $\mathfrak{S}_4$ is the permutation group on $4$ letters) and $384$ codimension $1$-facets. There are $12$ Levi subalgebras of length one, which are all conjugate as $\frakg$ is simply laced. Without loss of generality, we may assume that $\frakl=\frakl_{\alpha_*}$. The traces of the root hyperplanes on $\frakl$ form an arrangement with $7$ hyperplanes. One checks that $W(\frakl)\simeq (\Z/2\Z)^3$. Note that $\{\alpha_*\}^\perp$ is empty in this case so that we have a strict inclusion $W_I \times W_{I^\perp} \subset \Stab_W(H_\frakl)$. The group $W(\frakl)$ acts freely on the set of chambers of $H_\frakl$ which decomposes into $4$ orbits according to the types $I \subset \Pi$ of the corresponding parabolic subalgebras. 

\medskip

In order to describe the simplicial complex $H_\frakl$, we may intersect it with the unit sphere to obtain a triangulation of the $2$-sphere. The group $W(\frakl)$ is generated by the reflections along the three coordinate planes, thus it is enough to draw the triangulation of, say, the upper half-sphere, remembering that the entire complex is obtained by gluing a mirror image along the outer border. Up to homotopy, we obtain the complex described in Figure~\ref{FigureD_41} below. In addition to the three hyperplanes given by the coordinate axes, there are four more : we have indicated one of these hyperplanes in blue; the others are obtained by applying a $90$ degree rotations. We have indicated the type of each chamber, as well as the type of certain lower-dimensional cells. The types of the other cells are uniquely determined by the rules that
\begin{enumerate} 
\item if $C^\circ_\frakp \preceq C^\circ_{\frakp'}$ then the type of $\frakp'$ is contained in that of $\frakp$,
\item the types of the boundary cells of a cell of type $I$ run over all subsets of $\Pi$ containing $I$.
\end{enumerate}

\medskip

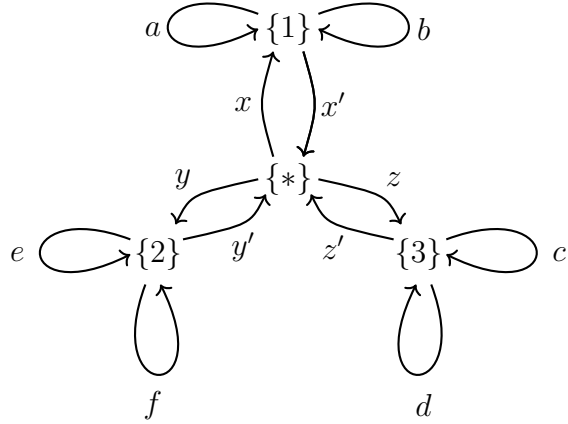
\begin{figure}
\begin{tikzpicture}
	\node(a)at(0,0){$\{*\}$};
	\node(a)at(90:2){$\{1\}$};
	\node(a)at(-30:2){$\{3\}$};
	\node(a)at(-150:2){$\{2\}$};	
	\node(a)at(-1.8,2){${a}$};
	\node(a)at(1.8,2){${b}$};
	\node(a)at(-0.6,1){$x$};
	\node(a)at(0.6,1){$x'$};
	\node(a)at(-1.4,0){$y$};
	\node(a)at(-0.6,-0.9){$y'$};
	\node(a)at(1.4,0){$z$};
	\node(a)at(0.6,-0.9){$z'$};
	\node(a)at(-1.8,-3){$f$};
	\node(a)at(1.8,-3){$d$};	
	\node(a)at(-3.6,-1){${e}$};
	\node(a)at(3.6,-1){${c}$};
	\draw[thick,->](-0.4,2.2)..controls(-2,2.8) and (-2,1.2)..(-0.4,2);
	\draw[thick,->](0.4,2.2)..controls(2,2.8) and (2,1.2)..(0.4,2);	
	\draw[thick,->](-0.2,0.3)..controls(-0.4,1)..(-0.2,1.7);		
	\draw[thick,->](0.2,1.7)..controls(0.4,1)..(0.2,0.3);
	\draw[thick,->](-2.1,-0.8)..controls(-3.7,-0.2) and (-3.7,-1.8)..(-2.1,-1);
	\draw[thick,->](-1.9,-1.4)..controls(-2.5,-3) and (-0.9,-3)..(-1.7,-1.4);
	\draw[thick,->](0.2,1.7)..controls(0.4,1)..(0.2,0.3);
	\draw[thick,->](2.1,-0.8)..controls(3.7,-0.2) and (3.7,-1.8)..(2.1,-1);
	\draw[thick,->](1.9,-1.4)..controls(2.5,-3) and (0.9,-3)..(1.7,-1.4);
	\draw[thick,->](-0.4,0)..controls(-1.3,-0.2)..(-1.5,-0.6);
	\draw[thick,<-](-0.3,-0.2)..controls(-0.5,-0.6)..(-1.4,-0.8);
	\draw[thick,->](0.4,0)..controls(1.3,-0.2)..(1.5,-0.6);
	\draw[thick,<-](0.3,-0.2)..controls(0.5,-0.6)..(1.4,-0.8);		
\end{tikzpicture}
\caption{Generators of the groupoid $Br(\frakl,W)$ in type $D_4$, $l=1$.}
\label{FigureBrD_41}
\end{figure}

We have also indicated the generators $x,y,z, x',y',z'$ and $a,b,c,d,e,f$ of the fundamental groupoid
$\Br(\frakl,W)$. The latter groupoid can be combinatorially depicted as in Figure~\ref{FigureBrD_41}, with the relations obtained by going around each codimension $2$ cell. This gives the folllowing set of relations~:
$$dc=cd, \qquad ef=fe, \qquad ab=ba,$$
$$z^{-1}dz=y^{-1}fy, \qquad x^{-1}ax=z^{-1}cz, \qquad y^{-1}ey=x^{-1}bx,$$
$$z'd(z')^{-1}=y'f(y')^{-1}, \qquad x'a(x')^{-1}=z'c(z')^{-1}, \qquad y'e(y')^{-1}=x'b(x')^{-1}.$$

In particular, setting
$$\gamma_1=x'x, \quad \gamma_2=y'y, \quad \gamma_3=z'z,\quad u_{13}=x^{-1}ax, \quad u_{12}=y^{-1}ey, \quad u_{23}=z^{-1}dz$$
we see that $\Br(\frakl,W)$ is equivalent to the group
$$Aut(*)=\langle \gamma_1, \gamma_2,\gamma_3, u_{12},u_{13},u_{23}\rangle / \{[u_{ij},u_{kl}]=0, [\gamma_i^{-1}\gamma_j,u_{ij}]=0\}.$$

\vspace{.1in}

\begin{figure}
	\begin{tikzpicture}
		\draw[thick] (-6,0)--(0,6)--(6,0)--(0,-6)--(-6,0);
		\draw[thick] (-3,0)--(0,3)--(3,0)--(0,-3)--(-3,0);
		\draw[thick] (-3,-3)--(-3,3)--(3,3);
		\draw[thick] (0,-3)--(3,-3)--(3,0);
		\draw[thick,blue] (-3,-3)--(0,-3)--(3,0)--(3,3);
		\draw[thick] (-6,0)--(6,0);
		\draw[thick] (0,6)--(0,-6);
		\draw[thick,red,->] (-1.7,2.7)--(-1.7,3.3);
		\draw[thick,red,->] (-1.3,2.9)--(-1.3,2.7);
		\draw[thick,red] (-1.3,3.3)--(-1.3,3.1);
		\draw[thick,red,->] (-2.7,1.7)--(-3.3,1.7);
		\draw[thick,red,->] (-2.9,1.3)--(-2.7,1.3);
		\draw[thick,red] (-3.3,1.3)--(-3.1,1.3);
		\draw[thick,red,->] (-1.8,1.6)--(-1.4,1.2);
		\draw[thick,red,<-] (-1.6,1.8)--(-1.45,1.65);
		\draw[thick,red] (-1.3,1.5)--(-1.2,1.4);
		\draw[thick,red,->] (-0.3,4)--(0.3,4);
		\draw[thick,red,->] (-1.4,4.3)--(-1.75,4.65);
		\draw[thick,red,->] (-0.3,1.2)--(0.3,1.2);
		\draw[thick,red,->] (-4.4,1.3)--(-4.75,1.65);
		\draw[thick,red,->] (-1,0.3)--(-1,-0.3);
		\draw[thick,red,->] (-4,0.3)--(-4,-0.3);
		\node(a)at(-1.9,2.7){${\color{red}{x}}$};
		\node(a)at(-1,3.4){${\color{red}{x'}}$};
		\node(a)at(-2.5,1.7){${\color{red}{y}}$};

		\node(a)at(-3.5,1.3){${\color{red}{y'}}$};
		\node(a)at(-2,1.6){${\color{red}{z}}$};
		\node(a)at(-0.95,1.5){${\color{red}{z'}}$};
		\node(a)at(-0.45,4){${\color{red}{a}}$};
		\node(a)at(-1.3,4.3){${\color{red}{b}}$};
		\node(a)at(-0.45,1.2){${\color{red}{c}}$};
		\node(a)at(-1,0.5){${\color{red}{d}}$};
		\node(a)at(-4.3,1.3){${\color{red}{e}}$};
		\node(a)at(-4,0.5){${\color{red}{f}}$};
		\node(a)at(-1,1){$3$};
		\node(a)at(-1,-1){$3$};
		\node(a)at(1,1){$3$};
		\node(a)at(1,-1){$3$};
		\node(a)at(-2,2){$*$};
		\node(a)at(-2,-2){$*$};
		\node(a)at(2,2){$*$};
		\node(a)at(2,-2){$*$};
		\node(a)at(-1,4){$1$};
		\node(a)at(-1,-4){$1$};
		\node(a)at(1,4){$1$};
		\node(a)at(1,-4){$1$};
		\node(a)at(-4,1){$2$};
		\node(a)at(-4,-1){$2$};
		\node(a)at(4,1){$2$};
		\node(a)at(4,-1){$2$};
		\node(o)at(0,0.4){$\{1,2,3\}$};
		\node(o)at(0,6.4){$\{1,2,3\}$};
		\node(o)at(6,0.4){$\{1,2,3\}$};
		\node(o)at(-6,0.4){$\{1,2,3\}$};
		\node(o)at(0,-6.4){$\{1,2,3\}$};
		\node(o)at(0,1.7){$\{1,3\}$};
		\node(o)at(0,4.6){$\{1,3\}$};
		\node(o)at(0,-1.7){$\{1,3\}$};
		\node(o)at(0,-4.6){$\{1,3\}$};
		\node(o)at(1.4,0){$\{2,3\}$};
		\node(o)at(-1.8,0){$\{2,3\}$};
		\node(o)at(4.2,0){$\{2,3\}$};
		\node(o)at(-4.8,0){$\{2,3\}$};
		\node(o)at(0,-6.4){$\{1,2,3\}$};
		\node(o)at(0,0.4){$\{1,2,3\}$};
		\node(o)at(0,6.4){$\{1,2,3\}$};
		\draw (0,0) circle (0.1);
		\draw (0,3) circle (0.1);
		\draw (0,6) circle (0.1);
		\draw (3,0) circle (0.1);
		\draw (6,0) circle (0.1);
		\draw (-3,0) circle (0.1);
		\draw (-6,0) circle (0.1);
		\draw (0,-3) circle (0.1);
		\draw (0,-6) circle (0.1);
		\draw (3,3) circle (0.1);
		\draw (3,-3) circle (0.1);
		\draw (-3,3) circle (0.1);
		\draw (-3,-3) circle (0.1);
	\end{tikzpicture}
	\caption{Type $D_4$, with $\frakl$ of length $l=1$.}
	\label{FigureD_41}
\end{figure}
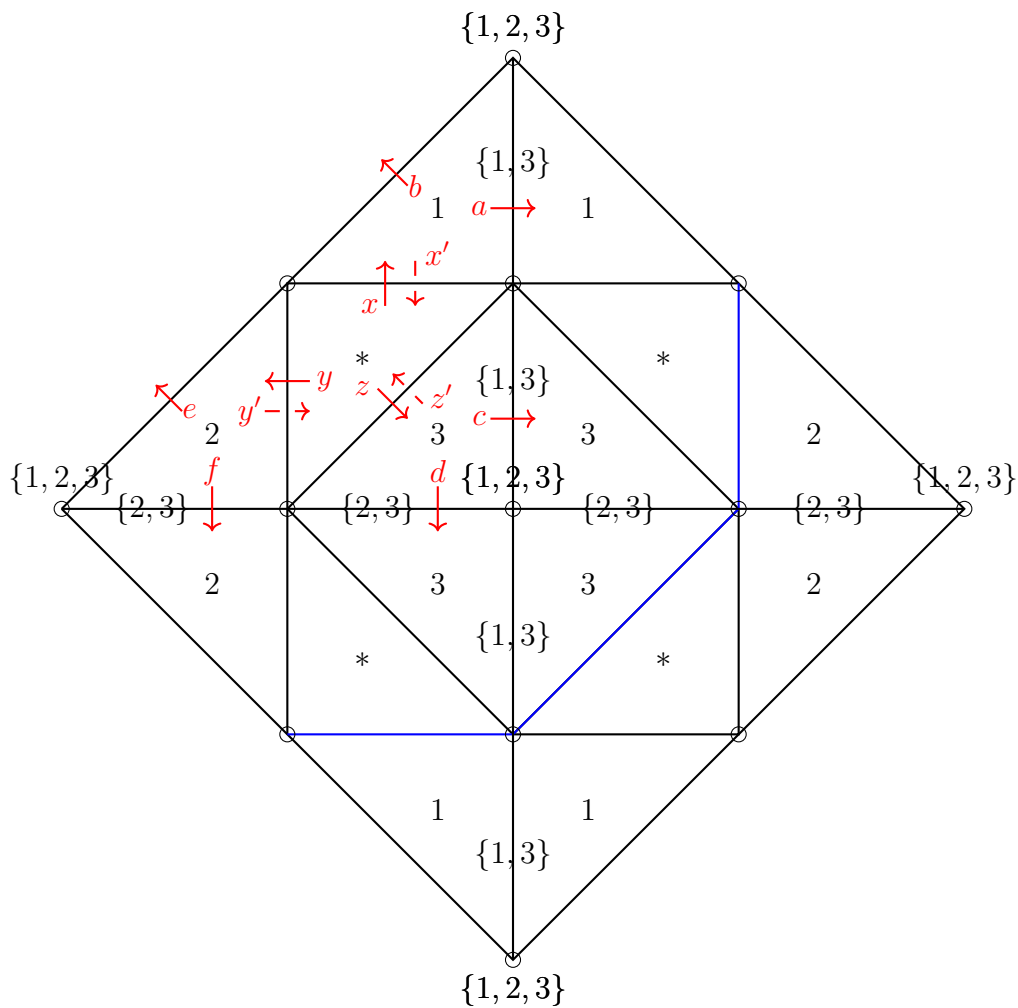

\vspace{.2in}

\subsection{{Type $D_4$ and $l=2$}} 
\begin{figure}
	\begin{tikzpicture}
		\draw[thick] (60:4)--(240:4);
		\draw[thick] (120:4)--(300:4);
		\draw[thick] (-4,0)--(4,0);
		\draw (0,0) circle (0.1);
		\node(a)at(0,2.5){$\{1,*\}$};
		\node(a)at(150:2.5){$\{3,*\}$};
		\node(a)at(30:2.5){$\{2,*\}$};
		\node(a)at(-150:2.5){$\{2,*\}$};
		\node(a)at(-30:2.5){$\{3,*\}$};
		\node(a)at(0,-2.5){$\{1,*\}$};
	\end{tikzpicture}
	\caption{Type $D_4$ with $l=2$; case of $I=\{\alpha,\beta\}$ with $(\alpha,\beta)=-1$.}
	\label{FigureD_42a}
\end{figure}
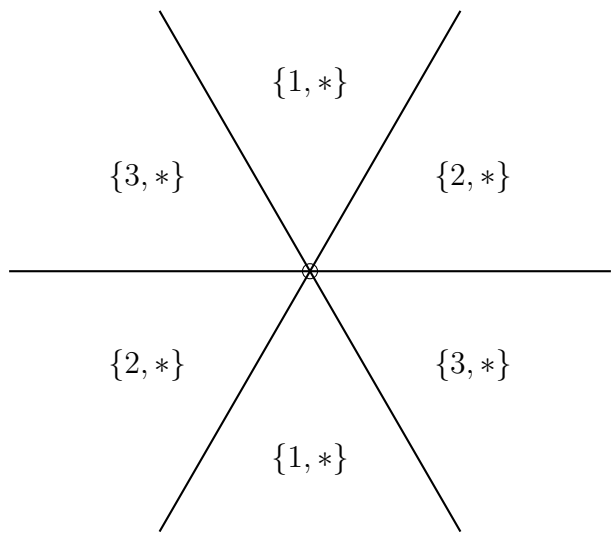

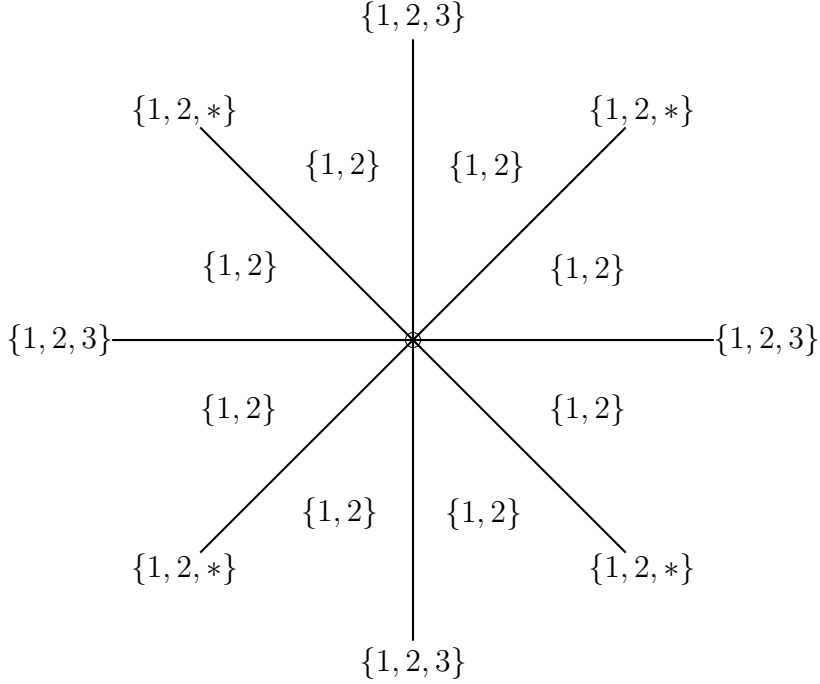
\begin{figure}[h]
	\begin{tikzpicture}
		\draw[thick] (45:4)--(45:-4);
		\draw[thick] (-45:4)--(-45:-4);
		\draw[thick] (-4,0)--(4,0);
		\draw[thick] (0,4)--(0,-4);
		\draw (0,0) circle (0.1);
		\node(a)at(67:2.5){$\{1,2\}$};
		\node(a)at(135:4.3){$\{1,2,*\}$};
		\node(a)at(135:-4.3){$\{1,2,*\}$};
		\node(a)at(-90:4.3){$\{1,2,3\}$};
		\node(a)at(45:-4.3){$\{1,2,*\}$};
		\node(a)at(22:2.5){$\{1,2\}$};
		\node(a)at(157:2.5){$\{1,2\}$};
		\node(a)at(112:2.5){$\{1,2\}$};
		\node(a)at(112:-2.5){$\{1,2\}$};
		\node(a)at(-22:2.5){$\{1,2\}$};
		\node(a)at(22:-2.5){$\{1,2\}$};
		\node(a)at(67:-2.5){$\{1,2\}$};
		\node(b1)at(4.7,0){$\{1,2,3\}$};
		\node(b1)at(-4.7,0){$\{1,2,3\}$};
		\node(c1)at(45:4.3){$\{1,2,*\}$};
		\node(e1)at(90:4.3){$\{1,2,3\}$};
	\end{tikzpicture}
	\caption{Type $D_4$ and $l=2$; case of $I=\{\alpha,\beta\}$ with $(\alpha,\beta)=0$.}
	\label{FigureD_42b}
\end{figure}

Keeping the notation of the preceding subsection, let us further intersect the complex of Figure~\ref{FigureD_41} with a root hyperplane; up to some natural isomorphisms there are two cases, corresponding to the hyperplane indicated in blue, or to one of the coordinate axes. Accordingly we respectively obtain the complexes in Figure~\ref{FigureD_42a} and Figure~\ref{FigureD_42b} (again, the types of unmarked cells is determined by rules (1) and (2) as above). In the case of Figure~\ref{FigureD_42a}, we have $W(\frakl)=\Z/2\Z$ (symmetry through the origin) and in Figure~\ref{FigureD_42b} we have $W(\frakl)=(\Z/2\Z)^2$ (reflection along the two coordinate axes).

\medskip

The combinatorial presentation of the fundamental groupoid for Figure~\ref{FigureD_42a} is as follows :

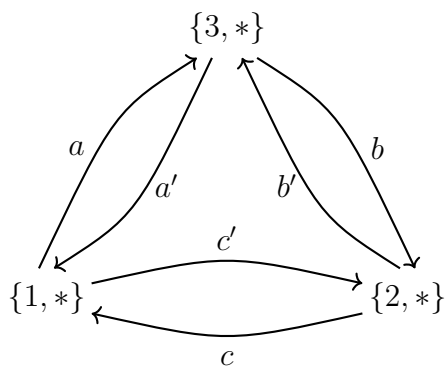
\centerline{
\begin{tikzpicture}
	\node(a)at(2.4,0){$\{2,*\}$};
	\node(a)at(-2.4,0){$\{1,*\}$};
	\node(a)at(0,3.6){$\{3,*\}$};
	\node(a)at(0,0.8){$c'$};
	\node(a)at(0,-0.8){$c$};
	\node(a)at(0.8,1.5){$b'$};
	\node(a)at(-0.8,1.5){$a'$};
	\node(a)at(2,2){$b$};
	\node(a)at(-2,2){$a$};						
	\draw[thick,->](-2.5,0.4)..controls(-1.5,2.5)..(-0.4,3.2);
	\draw[thick,<-](2.5,0.4)..controls(1.5,2.5)..(0.4,3.2);
	\draw[thick,<-](-2.3,0.4)..controls(-1.2,1.1)..(-0.2,3.2);
	\draw[thick,->](2.3,0.4)..controls(1.2,1.1)..(0.2,3.2);	
	\draw[thick,->](-1.8,0.2)..controls(0,0.6)..(1.8,0.2);
	\draw[thick,<-](-1.8,-0.2)..controls(0,-0.6)..(1.8,-0.2);		
\end{tikzpicture}
}

with relations $$(abc)^2=1, \qquad (c'b'a')^2=1.$$

The combinatorial presentation of the fundamental groupoid for Figure~\ref{FigureD_42b} is as follows :

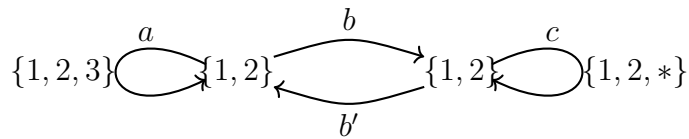
\centerline{
\begin{tikzpicture}
	\node(a)at(-1.5,0){$\{1,2\}$};
	\node(a)at(1.5,0){$\{1,2\}$};
	\node(a)at(-3.8,0){${\tiny{\{1,2,3\}}}$};
	\draw[thick,->](-1.9,0.1)..controls(-3.5,0.9) and (-3.5,-0.9)..(-1.9,-0.1);
	\node(a)at(3.8,0){${\tiny{\{1,2,*\}}}$};
	\draw[thick,->](1.9,0.1)..controls(3.5,0.9) and (3.5,-0.9)..(1.9,-0.1);
	\draw[thick,->](-1,0.2)..controls(0,.5)..(1,0.2);
	\draw[thick,<-](-1,-0.2)..controls(0,-.5)..(1,-0.2);	
	\node(a)at(-2.7,0.5){$a$};
	\node(a)at(2.7,0.5){$c$};
	\node(a)at(0,0.7){$b$};
	\node(a)at(0,-0.7){$b'$};
\end{tikzpicture}
}

with relations 
$$[a,b^{-1}cb]=0, \qquad [a,b'c(b')^{-1}]=0.$$

\medskip

\vfill\eject

\section{The perverse Coxeter (P-Coxeter) category}\label{sec:percox}

\vspace{.1in}

\subsection{ The Tits product and collinearity}  We keep the notation of \S \ref{sec:notation}. Let $\frakp, \frakq \in \calP_{G,T}$. The cells $C_{\frakp}, C_{\frakq}$ linearly generate a subspace $H_\frakl$  which is associated to a Levi subalgebra 
\[
\frakl =  \frakl(\frakp, \frakq) =   \pi(\frakp) \cap \pi(\frakq) \in \calL_{G,T}.
\]

The {\em Tits product}  $C^\circ_{\frakp}\circ C^\circ_{\frakq}$
of the (open) cells $C^\circ_{\frakp}$ and $C^\circ_{\frakq}$ is defined to be the (uniquely determined) first chamber  of $H_\frakl$  met by any line interval $[a,b]$ going  from $a\in C^\circ_{\frakp}$ to
$b\in C^\circ_{\frakq}$.  
Note that this chamber may be $C^\circ_{\frakp}$ itself (if $C^\circ_{\frakp}$ is a chamber of $H_\frakl$). 
See \cite[\S 2.7]{K-S-realhyper}  and references therein including the
original  work of Tits \cite{Tits}.

We denote by $\frakp \circ \frakq$ the parabolic subalgebra associated to the cell
$C^\circ_{\frakp}\circ C^\circ_{\frakq}$. The operation $\circ$ on $\calP_{G,T}$ is associative
but not commutative: 
$\frakp \circ \frakq$ differs in general from $\frakq \circ \frakp$ but they both are chambers of the same subspace $H_\frakl$. A direct  Lie-algebraic definition of  $\frakp \circ \frakq$  is as follows \cite[Prop. 5.5] {K-S-hW}:
\[
\frakp\circ\frakq  = (\frakp\cap\frakq)+ \frakn_\frakp. 
\] 	

We say that three parabolic subalgebras $\frakp, \frakq, \frakr\in\calP_{G,T}$ are {\em collinear}, if  there
are three points $a\in C^\circ_\frakp, b\in C^\circ_\frakq, c\in C^\circ_\frakr$ such that
$b\in[a,c]$.  In this case the cells  $ C^\circ_\frakp,  C^\circ_\frakq,  C^\circ_\frakr$ will be called
collinear as well. This concept depends on the order of the three subalgebras but is invariant
under $(\frakp,\frakq, \frakr) \mapsto (\frakr, \frakq, \frakp)$. 


\subsection{The double incidence category} \label{subsec:double-inc-cat}

Let $\calQ$ be the  $\k$-linear category defined by generators and relations as follows:

\begin{itemize}
\item[(0)] Objects of $\calQ$ are parabolic subalgebras $\frakp\in\calP_{G,T}$. 

\item[(1)] For any pair of parabolic subalgebras $\frakp,\frakp' \in \calP_{G,T}$ such that $\frakp \subset \frakp'$,  we have two generating morphisms
\[
\Ind_\frakp^{\frakp'}: \frakp\to  \frakp' , \qquad \Res_{\frakp'}^{\frakp}:  \frakp'\to  \frakp.
\]

\item[(2)] These morphisms are subject to the following relations:
\begin{itemize}
\item[(2a)] (Identity) $\Ind_\frakp^{\frakp} = \Res_{\frakp}^\frakp = \Id_\frakp$. 

\item[(2b)] (Transitivity, or (co)associativity) 
For any $\frakp \subset \frakp' \subset \frakp''$,
\[
\Ind_{\frakp'}^{\frakp''} \circ \Ind_{\frakp}^{\frakp'}=\Ind_{\frakp}^{\frakp''}, \qquad \Res_{\frakp'}^{\frakp}
\circ \Res_{\frakp''}^{\frakp'}=\Res_{\frakp''}^{\frakp}.
\]

\item[(2c)] (Idempotency) For any $\frakp\subset \frakp'$ we have 
$\Res_{\frakp'}^\frakp \circ\Ind_\frakp^{\frakp'} = \Id_\frakp$. 

The relations (2a)-(2c) allow one to define for any $\frakp, \frakp'\in \calP_{G,T}$ (not necessarily
contained one in another) a morphism $\tau_\frakp^{\frakp'}: \frakp \to\frakp'$ as
$\tau_\frakp^{\frakp'} = \Res_\frakq^{\frakp'} \circ \Ind_\frakp^\frakq$ for any
$\frakq\in\calP_{G,T}$ containing both $\frakp$ and $\frakp'$. The result is independent on $\frakq$;
one  can always take, for instance, $\frakq=\frakg$. Note that if $\frakp\subset\frakp'$, then 
(taking $\frakq=\frakp'$), we have $\tau_\frakp^{\frakp'}= \Ind_\frakp^{\frakp'}$ and
$\tau_{\frakp'}^\frakp = \Res_{\frakp'}^\frakp$. The remaning relations are formulated in terms of the
$\tau_\frakp^{\frakp'}$: 

\item[(2d)] (Collinear transitivity) If $\frakp, \frakp', \frakp''$ are collinear, then
$\tau_\frakp^{\frakp''} = \tau_{\frakp'}^{\frakp''}\circ \tau_{\frakp}^{\frakp'}$. 

\item[(2e)] (Invertibility) Let $\frakp'\subset \frakp \supset \frakp''$ be such that 
$C_{\frak p'}\supset C_\frakp\subset C_{\frakp''}$ are codimension $1$ inclusions of closed cells, that
$C_{\frakp'}, C_{\frakp''}$ have the same linear span and lie (in this span) on opposite sides of 
$C_\frakp$. Then $\tau_\frakp^{\frakp'}$  and $\tau_{\frakp'}^\frakp$ are invertible
(but not necessarily inverse to each other). 
\end{itemize}
\end{itemize}

We call $\calQ$ the {\em double incidence category} (associated to $G$ and $T$). 

\begin{figure}
	\begin{tikzpicture}
		\draw[thick] (45:2)--(45:-2);
		\draw[thick] (-45:2)--(-45:-2);
		\draw[thick] (-2,0)--(2,0);
		\draw[thick] (0,2)--(0,-2);
		\draw (0,0) circle (0.1);
		\draw[thick,red, ->, line join=round,
		decorate, decoration={
			zigzag,
			segment length=4,
			amplitude=.9,post=lineto,
			post length=2pt
		}](0.7,1.2)--(0,1.2); 
	\draw[thick,blue, ->, line join=round,
	decorate, decoration={
		zigzag,
		segment length=4,
		amplitude=.9,post=lineto,
		post length=2pt
	}](0,0)--(1.1,-0.6); 
		\node(a)at(0.85,1.3){$\frakp$};
		\node(e1)at(90:2.2){$\frakp'$};
		\node(o)at(0.2,0.5){$\frakg$};
		\node(o)at(1.3,-0.6){$\frakp''$};
	\end{tikzpicture}
	\caption{The morphisms $\text{Ind}_\frakp^{\frakp'}$ (in red) and the $\text{Res}_{\frakg}^{\frakp''}$ (in blue).}
	\label{FigureIndandRes}
\end{figure}
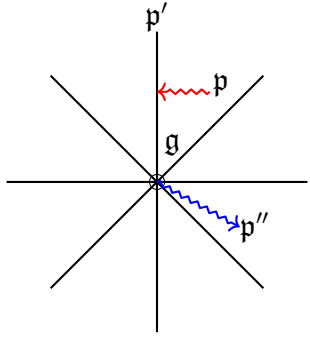

\begin{rem}
All generators and relations in $\calQ$ have a direct geometric interpretation in terms of the realization of the Coxeter complex in $\frakh$. The generating morphism $\Ind_{\frakp}^{\frakp'}$ corresponds to the operation of passing from a cell (face) $C^\circ_\frakp$ to a lower-dimensional cell $C^\circ_{\frakp'}$ in its closure; it is the other way around for the morphism $\Res^\frakp_{\frakp'}$, see Figure \ref{FigureIndandRes}. More generally, for any arrangement $\calE$ of hyperplanes in a real vector space $V$ one can define a category $\calQ(\calE)$ whose
objects are faces of $\calE$ and for any embedding $C'\subset \ol C$ of face closures there are two 
generating morphisms
$I_{C}^{C'}: C\to C'$, $R_{C'}^C: C'\to C$ satisfying relations similar to (2). 
These relations have been written in  \cite{K-S-realhyper} in the form of conditions on certain 
diagrams 
of vector spaces  (``double quivers'') labelled by the faces. Such double quivers  can be viewed as
as functors $\calQ(\calE) \to\Vect_\k$; as shown in \cite{K-S-realhyper}, they
describe perverse sheaves on the complexification $V_\C$
smooth with respect to the stratification into generic strata of the flats. See Section \ref{sec:Perv-hW} below. 
\flushright{$\triangle$} 
\end{rem}

\subsection {Analysis of the double incidence category}\label{subsec:anal-DIC}
We now highlight some features of the category $\calQ$ which will later give rise to a categorical version
of the Langlands formula.

\vskip .3cm

\paragraph{\bf {Braiding operators}} Let $\frakl\in \calL_{G,T}$ be a Levi subalgebra. The set
$\pi^{-1}(\frakl)\subset \calP_{G,T}$ of parabolics with Levi $\frakl$ is identified with the set of  chambers
(maximal, i.e., open, faces) of the arrangement in the flat $H_\frakl\subset\frakh_\R$ induced by
the root arrangement $\calH$ in $\frakh_\R$. For any $\frakp, \frakp'\in \pi^{-1}(\frakl)$ the
morphism $\tau_\frakp^{\frakp'}\in\Hom_\calQ(\frakp, \frakp')$ will be called the
{\em braiding operator}. 

Recall that $\Pi_1(H^\circ_{\frakl,\C})$ is the fundamental groupoid of the complex generic stratum
$H^\circ_{\frakl,\C}\subset H_{\frakl, \C}$ w.r.t. a set of base points, one in each real chamber.
Thus $\Ob \Pi_1(H^\circ_{\frakl,\C}) \= \pi^{-1}(\frakl)$. The relations (2d) of \S \ref{subsec:double-inc-cat}
and Proposition \ref{prop:pi_1-linear} imply: 

\begin{prop}
The morphisms $\tau_\frakp^{\frakp'}$, $\frakp, \frakp'\in \pi^{-1}(\frakl)$, give rise to a functor
$\Pi_1(H^\circ_{\frakl,\C})\to\calQ$. \qed
\end{prop}

\vskip .3cm

\paragraph{\bf Functoriality of $\Ind$ and $\Res$ w.r.t. braiding operators} \label{par:funct-Ind-Res}
Let $\frakl$ and $\frakl'$ be two adjacent Levi subalgebras, i.e. $H_{\frakl'} \subset H_\frakl$ is a hyperplane.
Let $H_\frakl^+$ be one of the open half spaces which are the connected components of 
$H_\frakl - H_{\frakl'}$  and let $\frakp_1, \frakp_2 \in \calP_{G,T}$ be parabolic subalgebras such that 
$\pi(\frakp_1)=\pi(\frakp_2)=\frakl'$. 
The (unique) chambers in $H^{+,\circ}_\frakl$ which are adjacent to $\frakp_1$ and $\frakp_2$ respectively correspond to parabolic subalgebras which we denote  $\frakp^+_1, \frakp_2^+ \in \calP_{G,T}$.
They satisfy
$\pi(\frakp_1^+)=\pi(\frakp^+_2)=\frakl$ and $C_{\frakp_i} \subset \overline{C_{\frakp^+_i}}$. 

\begin{figure}[h]
	\begin{tikzpicture}
		\draw[thick] (45:2)--(45:-2);
		\draw[thick] (-45:2)--(-45:-2);
		\draw[thick] (-2,0)--(2,0);
		\draw[thick] (0,2)--(0,-2);
		\draw[thick,red, ->, line join=round,
		decorate, decoration={
			zigzag,
			segment length=4,
			amplitude=.9,post=lineto,
			post length=2pt
		}](-1.3,1.3)--(-0.7,1.3);
		\draw[thick,blue, ->, line join=round,
		decorate, decoration={
			zigzag,
			segment length=4,
			amplitude=.9,post=lineto,
			post length=2pt
		}](1.3,-1.3)--(1.3,-0.7);
		\draw[thick, red, ->](-0.65,1.25)--(1.25,-0.65);
		\draw (0,0) circle (0.1);
		\draw[thick, blue, ->](-1.3,1.3)--(1.25,-1.25);
		\node(e1)at(135:2.2){$\frakp_1$};
		\node(a)at(-45:2.2){$\frakp_2$};
		\node(o)at(-0.6,1.8){$\frakp_1^+$};
		\node(o)at(1.8,-0.6){$\frakp_2^+$};
	\end{tikzpicture}
	\caption{Functoriality of $\text{Res}$ : the two morphisms $\frakp_1 \to \frakp_2^+$ are equal.}
	\label{Figurefunctoriality}
\end{figure}
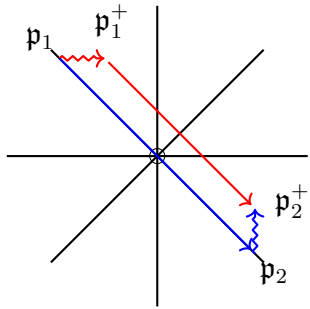

\begin{prop}\label{prop:Ind-Res-funct-Q}
We have the equalities of morphisms in $\calQ$, cf.  Fig. \ref{Figurefunctoriality}:
\[  
\tau^{\frakp^+_2}_{\frakp^+_1} \circ \Res_{\frakp_1}^{\frakp_1^+}=\Res_{\frakp_2}^{\frakp_2^+} \circ \tau^{\frakp_2}_{\frakp_1}, \qquad 
\tau^{\frakp_2}_{\frakp_1} \circ \Ind_{\frakp_1^+}^{\frakp_1}=\Ind_{\frakp_2^+}^{\frakp_2} \circ \tau^{\frakp^+_2}_{\frakp^+_1}.
\]
 \end{prop}
 
\noindent{\sl Proof:} We prove the first equality, the second one being similar.  By definition, 
$\tau_{\frakp_1}^{\frakp_2} = \Res_{\frakg}^{\frakp_2}\circ \Ind_{\frakp_1}^{\frakg}$,
so the RHS of the claimed equality is
\[
\Res_{\frakp_2}^{\frakp_2^+}\circ \Res_{\frakg}^{\frakp_2} \circ \Ind_{\frakp_1}^\frakg =
\Res_{\frakg}^{\frakp_2^+} \circ \Ind _{\frakp_1}^{\frakg} = \tau_{\frakp_1}^{\frakp_2^+}. 
\]
Now, interpreting $\Res_{\frakp_1}^{\frakp_1^+}$ as $\tau_{\frakp_1}^{\frakp_1^+}$, we write
the LHS as $\tau_{\frakp_1^+}^{\frakp_2^+}\circ\tau_{\frakp_1}^{\frakp_1^+}$ which is the same
as $\tau_{\frakp_1}^{\frakp_2^+}$ by relations (2d), since $\frakp_1, \frakp_1^+, \frakp_2^+$
are collinear. \qed

\vskip.3cm

\paragraph{\bf The proto-Langlands formula}    The following $1$-term expression of $\Res\circ\Ind$
will  later produce a full-fledged Langlands formula (with summation over double cosets). 

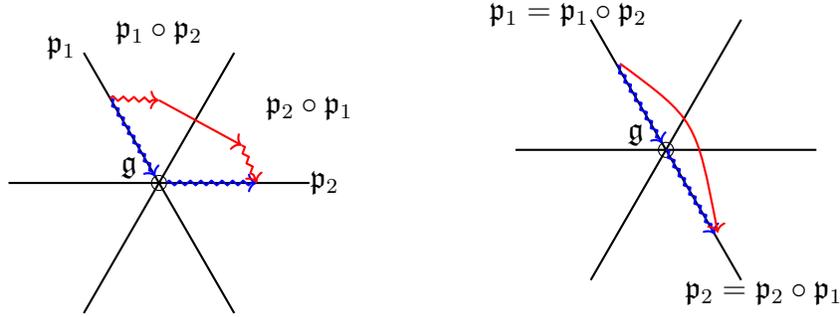
\begin{figure}[h]
	\begin{tikzpicture}
		\draw[thick] (60:2)--(240:2);
		\draw[thick] (120:2)--(300:2);
		\draw[thick] (-2,0)--(2,0);
		\draw (0,0) circle (0.1);
		\node(a)at(-1.3,1.8){$\frakp_1$};
		\node(a)at(2.2,0){$\frakp_2$};
		\node(a)at(-0.4,0.2){$\frakg$};
		\draw[thick,blue, ->, line join=round, decorate, decoration={
		zigzag, segment length=4, amplitude=.9,post=lineto, post length=2pt}](-60:-1.3)--(-60:-0.1);
		\draw[thick,blue, ->, line join=round, decorate, decoration={
		zigzag, segment length=4, amplitude=.9,post=lineto, post
		length=2pt}](0.1,0)--(1.3,0);
		\draw[thick,red, ->, line join=round, decorate, decoration={
		zigzag, segment length=4, amplitude=.9,post=lineto, post
		length=2pt}](-60:-1.3)--(0,1.1);
		\draw[thick,red, ->, line join=round, decorate, decoration={
		zigzag, segment length=4, amplitude=.9,post=lineto, post
		length=2pt}](1.1,0.5)--(1.3,0);
		\draw[thick, red, ->](0,1.1)--(1.1,0.5);
		\node(a)at(0,2){$\frakp_1 \circ \frakp_2$};
		\node(a)at(2,1){$\frakp_2 \circ \frakp_1$};
	\end{tikzpicture}
\hspace{.5in}
\begin{tikzpicture}
	\draw[thick] (60:2)--(240:2);
	\draw[thick] (120:2)--(300:2);
	\draw[thick] (-2,0)--(2,0);
	\draw (0,0) circle (0.1);
	\node(a)at(-1.3,1.8){$\frakp_1=\frakp_1 \circ \frakp_2$};
	\node(a)at(1.3,-1.9){$\frakp_2=\frakp_2 \circ \frakp_1$};
	\node(a)at(-0.4,0.2){$\frakg$};
	\draw[thick,blue, ->, line join=round, decorate, decoration={
	zigzag, segment length=4, amplitude=.9,post=lineto, post length=2pt}](-60:-1.3)--(-60:-0.1);
	\draw[thick,blue, ->, line join=round, decorate, decoration={
	zigzag, segment length=4, amplitude=.9,post=lineto, post
	length=2pt}](0,0)--(-60:1.3);
	\draw[thick, red, ->](-62:-1.3)..controls(0.4,0.4).. (-58:1.3);
\end{tikzpicture}
	\caption{Two examples of  the proto-Langlands formula.}
	\label{FigurepartialGKformula}
\end{figure}

\begin{prop}\label{prop:proto-Lang}
For any triple of parabolic subalgebras $\frakp_1, \frakp, \frakp_2$ in $\calP_{G,T}$ such that $\frakp_1 \subset \frakp \supset \frakp_2$ we have, see Fig. \ref{FigurepartialGKformula}: 
\[
\Res_{\frakp}^{\frakp_2} \circ \Ind_{\frakp_1}^{\frakp}=
\Ind^{\frakp_2}_{\frakp_2 \circ \frakp_1} \circ \tau_{\frakp_1 \circ \frakp_2}^{\frakp_2 \circ \frakp_1} \circ
\Res_{\frakp_1}^{\frakp_1 \circ \frakp_2}.
\]
\end{prop}

\noindent{\sl Proof:} Writing $\Res$ and $\Ind$ in the RHS as particular cases of $\tau$
and using the definition of $\tau$ in the LHS, we rewrite the claim as
\[
\tau_{\frakp_1}^{\frakp_2} = \tau_{\frakp_1\circ\frakp_1}^{\frakp_2} \circ
\tau_{\frakp_1\circ\frakp_2}^{\frakp_2\circ\frakp_1}\circ
\tau_{\frakp_1}^{\frakp_1\circ\frakp_2}.
\]
But this follows from the relations (2d), as the four parabolics
$\frakp_1, \frakp_1\circ\frakp_2, \frakp_2\circ\frakp_1, \frakp_2$ are collinear (i.e., any subset of three
them in this order, is). \qed


\subsection {The perverse Coxeter category}

The category $\calQ$ is defined in terms of the root arrangement $\calH$ in $\frakh_\R$. The Weyl group $W$
acts on $\frakh_\R$ preserving $\calH$.  Note, for example,  that the Tits product is $W$-equivariant, i.e. we have 
$(w \cdot \frakp) \circ (w \cdot \frakq)=w \cdot (\frakp \circ \frakq)$ for any $w , \frakp, \frakq$.

Therefore $W$ acts on $\calQ$ (in the strict sense, see \S  \ref{subsec:gr-act-cat}). 
Explicitly, $W$ acts on the set $\calP_{G,T}$ of objects of $\calQ$ (which are in bijection with faces of $\calH$)
and the action on the generating morphisms is given by
\[
w \cdot \Res_{\frakp}^{\frakq}:=\Res_{w \cdot \frakp}^{w \cdot \frakq}, \qquad 
w \cdot \Ind_{\frakq}^{\frakp}:=\Ind_{w \cdot \frakq}^{w \cdot \frakp}
\]
for any $\frakq \subset \frakp$ which implies 
\[
w \cdot \tau_{\frakp}^{\frakp'}=\tau_{w \cdot \frakp}^{w \cdot \frakp'}
\]
for any pair  $\frakp, \frakp'$. The relations are obviously preserved, so we  indeed get an  action on $\calQ$. 

\vskip .2cm

We denote by $\coxcat = \coxcat_{G,T}$ the invariant orbit category $\calQ^W$, see \S \ref{subsec:orb-cat-II},
and call it the {\em perverse Coxeter category} ({\em P-Coxeter category}, for short).
 Thus $\Ob\coxcat =  W\bs \calP_{G,T}$ is the set of $W$-orbits on
$\calP_{G,T}$ and the Hom-spaces are given by \eqref{eq:orbit-cat-inv-hom}. 

\vskip .2cm

As pointed out in \S \ref{subsec:comb-cox}
$W\bs\calP_{G,T}$ is identified with $\calP_{G,B}$, the set of standard parabolic subalgebras
$\frakp_I, I\subset \Pi$.
So we can think of objects of $\coxcat$ as being the $\frakp_I$. 

\begin{ex}
Let $G=GL_n$, so $\frakg = \gl_n$. Standard parabolics in $\gl_n$ are parametrized by ordered partitions
$\bfn = (n_1,\cdots, n_k)$ of $n$, see Example
\ref {ex:gl_n-strata}. So we can think of  objects of $\coxcat_{GL_n, T}$ as being such $\bfn$. 
\flushright{$\triangle$}
\end{ex}

\subsection{Transitivity and functoriality of $\uInd$ and $\uRes$ in the perverse Coxeter category}\label{subsec:tran-Ind-Res-C}

 Following  \S \ref{subsec:orb-cat-II},
for any $a\in\Hom_\calQ(\frakp_1,\frakp_2)$ we denote by 
$\ul a \in\Hom_\coxcat(W\frakp_1, W\frakp_2)$
the image of $a$ under the  averaging map $\xi_{\frakp_1, \frakp_2}$, see  \eqref{eq:map-xi}.
In general, averaging is not compatible with composition of morphisms, 
i.e.,  $\ul{a\circ b}\neq \ul a\circ \ul b$. However, there are instances when
equality holds. 

\begin{lemma}\label{lem:ulab}
Let $\frakp, \frakq, \frakr\in \calP_{G,T}$ and let $b\in\Hom_\calQ(\frakp,\frakq),
a\in \Hom_\calQ(\frakq,\frakr)$. Assume that:
\begin{itemize}
\item[(1)] $a$ is invariant under $\Stab(\frakq, \frakr)\subset W$ 
and $b$ under
$\Stab(\frakp,\frakq)$

\item[(2)] $\Stab(\frakp,\frakq) = \Stab(\frakp), \,\, \Stab(\frakq,\frakr) = \Stab(\frakq)$
or, alternatively,

\item[(2')] $\Stab(\frakp,\frakq) = \Stab(\frakq), \,\, \Stab(\frakq,\frakr) = \Stab(\frakr)$.
\end{itemize}
Then  $\ul{a\circ b}=  \ul a\circ \ul b$.
\end{lemma} 

\noindent{\sl Proof:}  Assume (1) and (2) for definiteness. Then, by (1), 
$\ul a$ and $\ul b$
are the collections
\be\label{eq:ula-ulb}
\ul a  = \bigl( wa\bigl|\,  w\in W/\Stab(\frakq,\frakr)\bigr), \quad \ul b = \bigl( yb\bigl| \, 
y\in
W/\Stab(\frakp, \frakq)\bigr).
\ee
Let $W(\frakp,\frakq)\subset \calP_{G,T}\times\calP_{G,T}$ be the $W$-orbit
of the pair $(\frakp, \frakq)$ and similarly for $(\frakq,\frakr)$ and $(\frakp,\frakr)$
By (2), the projections
\[
W(\frakq,\frakr) \lra W\frakq, \quad W(\frakp,\frakq) \to W\frakp
\]
are bijections and therefore
\be\label{eq:Wqr}
W(\frakq,\frakr)\times_{W\frakq} W(\frakp, \frakq)  \= W(\frakp,\frakr).  
\ee
Now \eqref{eq:ula-ulb} realizes $\ul a$ as a collection labelled by $W(\frakq,\frakr)$
and $\ul b$ as a collection labelled by $W(\frakp,\frakq)$. Further, from
(1) and (2) we see that $\ul{ab}$ is represented as a collection labelled
by $W(\frakp,\frakr)$, i.e., by the RHS of  \eqref{eq:ula-ulb}. Further, by the composition rule
in $\coxcat = \calQ^W$, the morphism $\ul a \circ\ul b$ 
is represented as a collection labelled by the LHS of \eqref{eq:ula-ulb}
and we see that the morphisms with the corresponding labels coincide, 
hence our statement. The case of (1) and (2') is similar. \qed

\vskip .2cm

We now deduce several consequences.

\begin{prop}\label{prop:uInd-tran}
Let $\frakp_1\subset \frakp_2\ \subset \frakp_3$ be three parabolics in
$\calP_{G,T}$. Then 
\[
\uInd_{\frakp_1}^{\frakp_3} =   \uInd_{\frakp_2}^{\frakp_3}\circ
\uInd_{\frakp_1}^{\frakp_2}, \quad
\uRes^{\frakp_3}_{\frakp_1} = \uRes_{\frakp_1}^{\frakp_2} \circ \uRes^{\frakp_3}_{\frakp_2}. 
\
\]
\end{prop}

\noindent{\sl Proof:}  To prove the first equality, we apply Lemma \ref{lem:ulab}
with conditions (1) and (2) to $\frakp=\frakp_1, \frakq=\frakp_2, \frakr=\frakp_3$
and $a=\Ind_{\frakp_2}^{\frakp_3}$, $b=\Ind_{\frakp_1}^{\frakp_2}$. The condition
(1) is clear, while (2) is a consequence of the following property 
of the $W$-action on the Coxeter complex.
If $w\in W$ preserves a face,
it fixes each point of this face and so preserves each subface. 
Our statement now follows from the lemma and the transitivity of
the $\Ind$-morphisms in $\calQ$. The case of $\uRes$ is similar, 
applying the lemma to $\frakp=\frakp_3, \frakq=\frakp_2, \frakr=\frakp_1$
and using  the conditions (1) and (2'). \qed.

\begin{prop}\label{prop:utau-tran}
Let $\frakp_1, \frakp_2, \frakp_3\in\calP_{G,T}$ be such that the faces
$C_{\frakp_1}, C_{\frakp_2}, C_{\frakp_3}\subset \frakh_\R$ have the same
linear span $L$, are open there and collinear. Then
$\ul\tau_{\frakp_1}^{\frakp_3} =
\ul\tau_{\frakp_2}^{\frakp_3} \circ \ul\tau_{\frakp_1}^{\frakp_2}$.
\end{prop} 

\noindent{\sl Proof:}  We apply Lemma  \ref{lem:ulab} to 
$\frakp=\frakp_1, \frakq=\frakp_2, \frakr=\frakp_3$ and $a = \tau_{\frakp_2}^{\frakp_3}$,
$b=\tau_{\frakp_1}^{\frakp_2}$. The condition (1) is again clear, while
(2) (as well as (2')) follows because 
\[
\Stab(\frakp_i,\frakp_j) = \Stab(\frakp_i) = \Stab(\frakp_j)
\]
coincides with the stabilizer of $L$.  \qed

\begin{prop}\label{prop:uRes-funct}
In the situation of  Proposition \ref{prop:Ind-Res-funct-Q}, we have equalities of morphisms in $\coxcat$: 
\[  
\ul \tau^{\frakp^+_2}_{\frakp^+_1} \circ \uRes_{\frakp_1}^{\frakp_1^+}=\uRes_{\frakp_2}^{\frakp_2^+} \circ\ul \tau^{\frakp_2}_{\frakp_1}, \qquad 
\ul\tau^{\frakp_2}_{\frakp_1} \circ \uInd_{\frakp_1^+}^{\frakp_1}=\uInd_{\frakp_2^+}^{\frakp_2} \circ\ul  \tau^{\frakp^+_2}_{\frakp^+_1}.
\]
\end{prop}

\noindent{\sl Proof:} To prove the first statement, we apply Lemma  \ref{lem:ulab}  to get the equalities
\[
\ul \tau^{\frakp^+_2}_{\frakp^+_1} \circ \uRes_{\frakp_1}^{\frakp_1^+} = \ul{
\tau^{\frakp^+_2}_{\frakp^+_1} \circ \Res_{\frakp_1}^{\frakp_1^+} } \buildrel { \on{Pr.} \ref{prop:Ind-Res-funct-Q} } \over = 
\ul{
\Res_{\frakp_2}^{\frakp_2^+} \circ\ \tau^{\frakp_2}_{\frakp_1} } 
=
\uRes_{\frakp_2}^{\frakp_2^+} \circ\ul \tau^{\frakp_2}_{\frakp_1}
\]
The conditions (1) and (2) for the individual morphisms have been already established 
in the proofs of Propositions \ref{prop:uInd-tran} and  \ref{prop:utau-tran}. The second statement is proved similarly.  \qed

	
\subsection{The Langlands formula}\label{subsec:Lan-form}
We now deduce the following.

\begin{prop}[Langlands formula] \label{prop:langlands-form}
Let $\frakp_1, \frakp_2$ be parabolic subalgebras, $O_1=W \frakp_1, O_2=W\frakp_2$. Let also $W_i=\Stab(\frakp_i)\subset W$ be the Weyl group of $\frakp_i$. Then
\[
\uRes_{\frakg}^{\frakp_2} \circ \uInd_{\frakp_1}^\frakg=
	\sum_{w \in W_1 \backslash W /W_2} \uInd_{w \cdot \frakp_2 \circ\frakp_1}^{w \cdot \frakp_2}
	 \circ \underline{\tau}_{\frakp_1 \circ w\cdot \frakp_2}^{w \cdot \frakp_2 \circ\frakp_1} \circ 
	 \uRes^{\frakp_1 \circ w\cdot \frakp_2}_{\frakp_1}
 \]
\end{prop}

\noindent{\sl Proof:}
We first note that for any $\frakp\ \frakq\in\calP_{G,T}$ with 
$\frakq \subset \frakp$ the morphisms  $\Ind_\frakq^\frakp$
and $\Res_\frakp^\frakq$ in $\calQ$ are invariant under $\Stab(\frakp, \frakq)\subset W$.
Therefore 
\[
\ul\Ind_{\frakp_1}^\frakg = \bigl( \Ind_{\frakq_1}^\frak g\bigr)_{\frakq_1\in O_1}
\in \biggl(\bigoplus_{\frakq_1\in O_1} \Hom_\calQ(\frakq_1,\frakg)\biggr)^W, \quad
\ul\Res_{\frak g}^{\frakp_2} = \bigl( \Res_\frakg^{\frakq_2}\bigr)_{\frakq_2\in O_2}
\in\biggl( \bigoplus_{\frakq_2\in O_2} \Hom_{\calQ}(\frakg, \frakq_2)\biggr)^W.
\]

By definition of composition in the orbit category and the proto-Langlands formula
(Proposition \ref{prop:proto-Lang}), the LHS of the claimed equality is the collection
\be\label{eq:LHS-Lang}
\biggl(\Ind_{\frakq_2 \circ \frakq_1}^{\frakq_2} \circ \tau_{\frakq_1 \circ \frakq_2}^{\frakq_2 \circ \frakq_1} \circ \Res_{\frakq_1}^{\frakq_1 \circ \frakq_2}   \biggr)_{\frakq_i \in O_i}  \in
 \biggl( \bigoplus_{\frakq_i\in O_i} \Hom_\calQ(\frakq_1, \frakq_2)\biggr)^W. 
\ee
Now, the $W$-invariant collection of morphisms given by   \eqref{eq:LHS-Lang} is uniquely recovered from the
subcollection given by representatives of the $W$-orbits on $O_1\times O_2$. Noticing that
\[
W\bs(O_1\times O_2) \= W_1\bs W/W_2, 
\]
we identify this subcollection with 
\[
\biggl(\Ind_{w \cdot \frakp_2 \circ\frakp_1}^{w \cdot \frakp_2}
	 \circ {\tau}_{\frakp_1 \circ w\cdot \frakp_2}^{w \cdot \frakp_2 \circ\frakp_1} \circ 
	 \Res^{\frakp_1 \circ w\cdot \frakp_2}_{\frakp_1}\biggr)_{w \in W_1 \backslash W /W_2}
\]
The morphism obtained by extending this subcollection by $W$-invariance, is precisely the RHS if the claimed equality.
\qed

\vskip .3cm

More generally, for any $\frakp_1 \subset \frakp_3 \supset \frakp_2$ 
with Weyl groups $W_1\subset W_3\supset W_2$
we have
\be
\uRes_{\frakp_3}^{\frakp_2} \circ \uInd_{\frakp_1}^{\frakp_3}=
\sum_{w \in W_1 \backslash W_3 /W_2}
\uInd_{w \cdot \frakp_2 \circ\frakp_1}^{w \cdot \frakp_2} \circ \ul{\tau}_{\frakp_1 \circ w\cdot \frakp_2}^{w \cdot \frakp_2 \circ\frakp_1} \circ \uRes^{\frakp_1 \circ w\cdot \frakp_2}_{\frakp_1}.
\ee
The proof is identical to the above. 

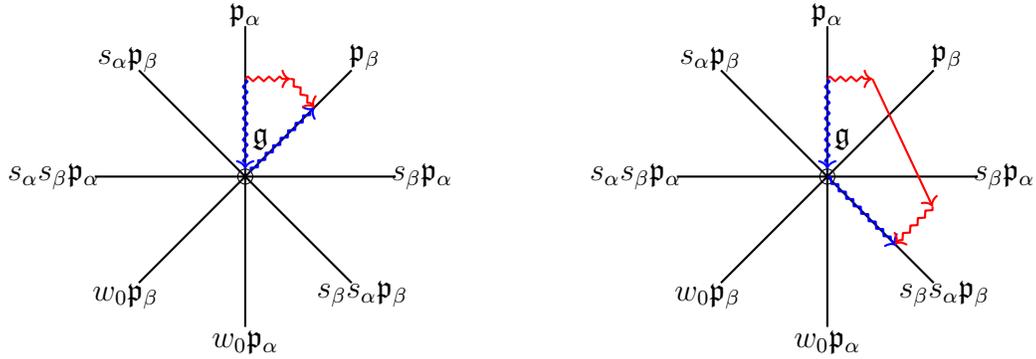
\begin{figure}[h]
	\begin{tikzpicture}
		\draw[thick] (45:2)--(45:-2);
		\draw[thick] (-45:2)--(-45:-2);
		\draw[thick] (-2,0)--(2,0);
		\draw[thick] (0,2)--(0,-2);
		\draw (0,0) circle (0.1);
		\node(a)at(135:2.2){$s_\alpha\frakp_\beta$};
		\node(a)at(135:-2.2){$s_\beta s_\alpha\frakp_\beta$};
		\node(a)at(-90:2.2){$w_0\frakp_\alpha$};
		\node(a)at(45:-2.25){$w_0\frakp_\beta$};
		\node(b1)at(2.35,0){$s_\beta\frakp_\alpha$};
		\node(b1)at(-2.55,0){$s_\alpha s_\beta\frakp_\alpha$};
		\node(c1)at(45:2.25){$\frakp_\beta$};
		\node(e1)at(90:2.15){$\frakp_\alpha$};
		\node(o)at(0.2,0.5){$\frakg$};
		\draw[thick,blue, ->, line join=round, decorate, decoration={
			zigzag, segment length=4, amplitude=.9,post=lineto, post length=2pt}](-90:-1.3)--(-90:-0.1);
		\draw[thick,blue, ->, line join=round, decorate, decoration={
			zigzag, segment length=4, amplitude=.9,post=lineto, post
			length=2pt}](45:0.1)--(45:1.3);
		\draw[thick,red, ->, line join=round, decorate, decoration={
			zigzag, segment length=4, amplitude=.9,post=lineto, post
			length=2pt}](-90:-1.3)--(0.6,1.3);
		\draw[thick,red, ->, line join=round, decorate, decoration={
			zigzag, segment length=4, amplitude=.9,post=lineto, post
			length=2pt}](0.6,1.3)--(45:1.3);
	\end{tikzpicture}
	\hspace{.5in}
	\begin{tikzpicture}
		\draw[thick] (45:2)--(45:-2);
		\draw[thick] (-45:2)--(-45:-2);
		\draw[thick] (-2,0)--(2,0);
		\draw[thick] (0,2)--(0,-2);
		\draw (0,0) circle (0.1);
		\node(a)at(135:2.2){$s_\alpha\frakp_\beta$};
		\node(a)at(135:-2.2){$s_\beta s_\alpha\frakp_\beta$};
		\node(a)at(-90:2.2){$w_0\frakp_\alpha$};
		\node(a)at(45:-2.25){$w_0\frakp_\beta$};
		\node(b1)at(2.35,0){$s_\beta\frakp_\alpha$};
		\node(b1)at(-2.55,0){$s_\alpha s_\beta\frakp_\alpha$};
		\node(c1)at(45:2.25){$\frakp_\beta$};
		\node(e1)at(90:2.15){$\frakp_\alpha$};
		\node(o)at(0.2,0.5){$\frakg$};
		\draw[thick,blue, ->, line join=round, decorate, decoration={
			zigzag, segment length=4, amplitude=.9,post=lineto, post length=2pt}](-90:-1.3)--(-90:-0.1);
		\draw[thick,blue, ->, line join=round, decorate, decoration={
			zigzag, segment length=4, amplitude=.9,post=lineto, post
			length=2pt}](0,0)--(-45:1.3);
		\draw[thick,red, ->, line join=round, decorate, decoration={
			zigzag, segment length=4, amplitude=.9,post=lineto, post
			length=2pt}](-90:-1.3)--(0.6,1.3);
		\draw[thick,red, ->, line join=round, decorate, decoration={
			zigzag, segment length=4, amplitude=.9,post=lineto, post
			length=2pt}](1.4,-0.4)--(0.9,-0.9);
		\draw[thick,red, ->](0.6,1.3)--(1.4,-0.4);
\end{tikzpicture}
	\caption{The Langlands formula for $B_2$. }
	\label{FigureGKformula}
\end{figure}

\begin{ex} Let us consider the root system $B_2$ again. There are two $W$-orbits of maximal proper parabolic subalgebras $O_1=W \cdot \frakp_\alpha$ and $O_2=W \cdot \frakp_\beta$. The composition 
$\uRes_{\frakg}^{\frakp_\beta} \circ \uInd^{\frakg}_{\frakp_\alpha}$
is given by a sum of two terms corresponding to the two  (up to $W$ conjugation)
paths going from a cell in $W \cdot \frakp_\alpha$ to the closed cell $\frakg$ to a cell in $W \cdot \frakp_\beta$, as in Fig. 
\ref{FigureGKformula}.
\flushright{$\triangle$}
\end{ex}

 \begin{rem}\label{rem:cox-gen-rel}
 This  section exhibits several morphisms (braidings, $\Ind$, $\Res$)  in $\coxcat$ and relations among them.
 We conjecture that they can be used to define $\coxcat$ abstractly, by these generators
 and relations. We plan to address this in a subsequent paper. 
 
 \end{rem}

\vfill\eject

\section{The case of $GL_n$}\label{sec:GLn}

The aim of this Section is to make explicit the structure of the P-Coxeter category for the collection of reductive groups $GL_n(\C)$ and show that this recovers the structure of the PROB of (graded) braided bialgebras, as in \cite{KS-Prob}. 

\medskip

\subsection{Combinatorial identifications}\label{subsec:comb-iden}

Let  $I$ be a finite set. We denote $GL_I = GL(\C^I)$ and  let $T_I \subset GL(I)$ be  the 
maximal torus (invertible diagonal matrices) 
corresponding to the standard basis  of $\C^I$.

We denote by $\calQ_I$ and $\calC_{I}$ the double incidence category, resp. the P-Coxeter category associated to the pair 
$(GL_I,T_I)$.   
When $I=\{1,\cdots, n\}$ we  will simply write $T_n$ for $T_{\{1,\cdots, n\}}$ and use
the notations
$\calQ_n$ and $\coxcat_n$ for the above categories. Thus if $|I|=n$, then there are 
isomorphisms $\calQ_{I} \simeq \calQ_n$ and $\calC_{I} \= \coxcat_n$.

Recall that $\Ob(\calQ_n) = \calP_{GL_n, T_n}$ is the set of parabolic subalgebras of $\gl_n$ containing
$\frakh_n = \Lie(T_n)$ (the set of diagonal matrices). Such subalgebras are in bijection
with  {\em ordered set partitions},
i.e.,  sequences of subsets 
\[
\bfI=(I_1, \ldots, I_k), \quad \emptyset\neq I_\nu\subset \{1,\cdots, n\}, 
\quad I_1 \sqcup I_2 \sqcup \ldots \sqcup I_k=\{1, \ldots, n\}.
\]
For such $\bfI$ and $i\in \{1,\cdots, n\}$ we denote $\nu_\bfI(i)\in \{1,\cdots, k\}$ to be the unique $\nu$
such that $i\in I_\nu$. 
The parabolic subalgebra $\frakp_\bfI$ corresponding to an ordered set partition
$\bfI$ is  the stabilizer of the flag 
\[
F_{\bfI}= V_1 \subset V_2 \subset \ldots \subset V_k=\C^n,
\text{ where }  V_\nu=\bigoplus_{i \in I_1 \sqcup \ldots \sqcup I_\nu } \C e_i.
\]
Recall further that $\calP_{GL_n, T_n}$ is identified with the set of faces (cells) of $\calH_n$,
the arrangement of the diagonal hyperplaces in $\frakh_{n,\R} = \R^n$. If a parabolic subalgebra
$\frakp_\bfI$ is associated to an ordered set partition $\bfI$, then its  corresponding face is
\be\label{eq:cell-set-part}
C_\bfI \,=\, \bigl\{ (t_1,\cdots t_n)\in \R^n\bigl| \, t_i =  t_j \text{ when } \nu_\bfI(i) = \nu_\bfI(j) 
\text{ and } t_i < t_j \text{ when } \nu_\bfI(i) < \nu_\bfI(j)
\bigr\}. 
\ee
In plain words $(t_1,\cdots, t_n)\in C_\bfI$ iff the pattern of equalities and inequalities among the $t_i$
is given by $\bfI$.

Any ordered set partition $\bfI=(U_1,\cdots, I_k)$ defines an unordered set partition
$\{I_1,\cdots, I_k\}$ by forgetting the order of the $I_\nu$. The following is straightforward.

\begin{prop}\label{prop:faces-GL}
Let $\bfI, \bfJ$ be two ordered set partitions of $\{1,\cdots, n\}$. Then:

\vskip .2cm 

(a)    $C_\bfI \subset \ol C_{\bfJ}$ if and only if $\bfJ$ is a refinement of   $\bfI$. 

\vskip .2cm

(b)  $C_\bfI, C_\bfJ$  are chambers of the same flat of $\calH_n$  if and only if the corresponding
unordered partitions are equal.

\vskip .2cm

(c) Further,  $C_\bfI, C_{\bfJ}$  are chambers of the same flat  that lie on opposite sides of a wall if and only if 
$\bfI$ is obtained from $\bfJ$ by exchanging two adjacent parts. \qed
\end{prop}

In fact, the hyperplane arrangement induced on the flat $H_\frakl$ associated to an unordered partition, 
$\{I_1, \ldots, I_k\}$ of $\{1, \ldots, n\}$ is isomorphic to that of type $A_{k-1}$. In particular, the group $W(\frakl)$ is isomorphic to the symmetric group $\mathfrak{S}_k$ and acts on chambers by permuting the parts.
Further, the Tits product on $\Ob(\calQ_n)$ can be combinatorially described as follows,
cf. \cite[\S 1.6] {aguiar:species-Hopf} and references therein. 

\begin{prop}\label{prop:Tits-GL} 
Let $\bfI=(I_1, \ldots, I_k)$ and $\bfJ=(J_1, \ldots, J_l)$ be two ordered set partitions of $\{1, \ldots, n\}$. 
Set $K_{i,j}=I_i \cap J_j$, so that we have a disjoint decomposition
$\bigsqcup_{i,j} K_{i,j}=\{1,\ldots, n\}$. Endowing the set of pairs $(i,j)$ with $1 \leq i \leq k, 1 \leq j \leq l$ with the lexicographic order and removing any empty entry, we obtain an ordered set partition $\bfK$ (with at most $kl$ sets).
Then $C_\bfI \circ C_\bfJ =C_\bfK$.
\end{prop}

\noindent{\sl Proof:}  
Let $(t_1,\cdots, t_n)\in C_\bfI$ and $(s_1,\cdots, s_n)\in C_\bfJ$. 
By  definition of $\circ$ and by \eqref{eq:cell-set-part}, $C_\bfI \circ C_\bfJ = C_\bfL$, 
where $\bfL$ expresses the pattern of equalities and inequalities among the numbers $t_i+\eps s_i$,
where $\eps>0$ is very small. Applying  \eqref{eq:cell-set-part} to $C_\bfI$ and $C_\bfJ$,
we see that the equivalence relation on $\{1,\cdots, n\}$ given by the 
pattern of equalities among the $t_i+\eps s_i$   has, as its equivalence classes, precisely the
nonempty $K_{i,j}$. After this we see that the pattern of strict inequalities is given by the lexicographic
order on the $K_{i,j}$, i.e., $\bfL=\bfK$. \qed   

\vskip .2cm

Further, $\Ob(\coxcat_n)$, i.e., the set of $\frakS_n$-orbits on the collection of ordered set partitions
of $\{1,\cdots, n\}$, is 
the collection of ordered partitions $\bfn=(n_1, \ldots, n_k)$ of the number  $n$. 
The projection map $\Ob(\calQ_n) \to \Ob(\coxcat_n)$ sends $(I_1, \ldots, I_k)$ to $(|I_1|, \ldots, |I_k|)$.


\subsection{Graded bialgebras }\label{subsec:grad-set-bi}
In this paper we will use the term ``braided category' to signify, for short,  a $\k$-linear braided monoidal category. 
Let $(\calV, \otimes)$ be a  braided  category  with the unit object and the 
braiding  denoted by $\1$ and $R= (R_{VW}: V\otimes W\to W\otimes V)$ respectively. 
Recall that a {\em bialgebra} in $\calV$ is an object $A\in\calV$ together with an associative
multiplication $\mu: A\otimes A\to A$, coassociative comultiplication $\Delta: A\to A\otimes A$,
unit $\eps: \1\to A$ and a counit $\eta: A\to\1$ satisfying:
\begin{itemize}
\item [(B1)]  (Unitality) The two compositions $\mu\circ(\eps\otimes\Id), \mu\circ(\Id\otimes\eps): A\to A$ are equal to $\Id_A$.

\item [(B2)] (Counitality) The two compositions $(\eta\otimes\Id)\circ\Delta, (\Id\otimes\eta)\circ\Delta: A\to A$
are equal to $\Id_A$.
\item  [(B3)] (Exchange) The diagram commutes:
\[
\xymatrix{
A\otimes A \ar[d]_{\Delta\otimes\Delta}
\ar[r]^\mu  &  A \ar[r]^\Delta  & A\otimes A
\\
A\otimes A\otimes A\otimes A \ar[rr]_{\Id\otimes R_{AA}\otimes\Id} && A\otimes A\otimes A\otimes A.
\ar[u]_{\mu\otimes \mu}
}
\]

\end{itemize}
We will need two versions of this concept, one in this and one in the next subsection. 

\begin{ex} (Graded bialgebras). Assume first that $\calV$ has direct sums
(which by $\k$-linearity of $(\calV, \otimes)$ are then automatically distributive with respect to $\otimes$).
Let $\calV^{\N}$ be the category of $\N$-graded objects $V=\bigoplus_{n=0}^\infty V_n$ of $\calV$.
The category $\calV^{\N}$ inherits from $\calV$
the braided monoidal structure. 
An object  of   $\calV^{\N}$    can also be thought as a collection $V=(V_n)_{n\geq 0}$ of its graded components,
so that
\be\label{eq:graded-otimes}
(V\otimes W)_n = \bigoplus_{n_1+n_2=n} V_{n_1} \otimes W_{n_2}. 
\ee
By a {\em graded
bialgebra} in $\calV$ we will mean a bialgebra  $A$ in    $\calV^{\N}$ with  
$A_0=\1$ and  such that   $\eps$
is the embedding $\1=A_0\to A$ and $\eta$ is the projecton $A\to A_0=\1$. 

More explicitly, a graded bialgebra in $\calV$ can be defined componentwise, as a datum of objects $A_n\in\calV$, $A_0=\1$
and morphisms
\[
\mu_{n_1, n_2}: A_{n_1}\otimes A_{n_2} \to A_{n_1+n_2}, \quad \Delta_{n_1, n_2}: A_{n_1+n_2}
\to A_{n_1}\otimes A_{n_2}, \quad n_1, n_2\geq 0,
\]
satisfying

\begin{itemize}

\item [(GB1)]  (Unitality) The morphisms $\mu_{0,n}, \mu_{n,0}$ are equal to  the unit isomorphisms.

\item [(GB2)]  (Counitality) The morphisms $\Delta_{0,n}, \Delta_{n,0}$ are equal to the inverses
of the unit isomorphisms.

\item [(GB3)]  (Exchange) Let $p_1+p_2 = q_1+q_2=n$. Then the composition

\[
A_{p_1}\otimes A_{p_2} \buildrel \mu_{p_1,p_2}\over\lra A_n \buildrel \Delta_{q_1, q_2} \over\lra A_{q_1}\otimes
A_{q_2}
\]
is equal to the sum of compositions
\[
\begin{gathered}
\sum_{\begin{bmatrix} n_{11}& n_{12} \\ n_{21} & n_{22}\end{bmatrix} }
\biggl\{ A_{p_1}\otimes A_{p_2}\buildrel \Delta_{n_{11}, n_{22}}\otimes\Delta_{n_{21}, n_{22}}\over\lra
A_{n_{11}}\otimes A_{n_{12}}\otimes A_{n_{21}} \otimes A_{n_{22}} 
\buildrel \Id\otimes R_{A_{n_{12}}, A_{n_{21}}}\otimes\Id\over\lra 
\\
\lra
A_{n_{11}}\otimes A_{n_{21}} \otimes A_{n_{12}}\otimes A_{n_{22}}
\buildrel \mu_{n_{11}, n_{21}}\otimes\mu_{n_{12}, n_{22}}\over\lra A_{q_1}\otimes A_{q_2} \biggr\},
\end{gathered}
\]
the summation being over matrices
\[
\begin{bmatrix} n_{11}& n_{12} \\ n_{21} & n_{22}\end{bmatrix}\in {\Mat}_2(\Z_+),\quad \text{s.t.} \quad 
n_{i1}+n_{i2}=p_i, \quad n_{1j}+n_{2j}=q_j. 
\]
\end{itemize}
This explicit  form allows us to speak about graded bialgebras in $\calV$ without assuming the existence
of direct sums. 
\end{ex}

\subsection{The braided category of $\frakS$-modules and set bialgebras}

Assuming again the existence of   
direct sums in $\calV$, we can construct a new
braided   category $\calV^\frakS$ of {\em $\frakS$-modules}, or {\em species} in $\calV$,
cf.   \cite{aguiar-species-book}  
\cite {aguiar:species-Hopf}.
By definition, an object of $\calV^\frakS$ is a collection $X=(X_n)_{n\geq 0}$ of objects of $\calV$
together with the action of the symmetric group $\frakS_n$ on  $X_n$ for each $n$. The category
$\calV^\frakS$ inherits from $\calV$ the
braided monoidal structure, with the monoidal operation $\odot$ given by
\be\label{eq:odot}
(X\odot Y)_n = \bigoplus_{n_1+n_2=n} \Ind_{\frakS_{n_1}\times\frakS_{n_2}}^{\frakS_n} 
(X_{n_1}\otimes Y_{n_2}).
\ee
The unit object of $\calV^\frakS$ is the collection $\1_{\calV^\frakS}$ whose $0$th component is $\1\in\calV$
and other components are $0$. 

Alternatively, we can view an object $V$ of $\calV^\frakS$ as a $\calV$-valued  functor 
\[
X: I\mapsto V_I, \quad I\mapsto X_I
\]
from the groupoid of finite sets and their isomorphisms, so that $X_n$ in the earlier 
sense is recovered as $X_{\{1,\cdots, n\}}$. 
In this language the formula \eqref{eq:odot} for $\odot$ appears as
\be
(X\odot Y)_I = \bigoplus_{I= I_1\sqcup I_2} X_{I_1}\otimes Y_{I_2},
\ee
the direct sum over all disjoint decompositions ${I= I_1\sqcup I_2}$. 

If $\calV$ is Karoubi-complete (possesses the images of all projectors), then  for any object $E\in\calV$
with an action of a finite group $G$ we have a well defined subobject of invariants $E^G$,
see \S \ref{subsec:karoubi}.  In this
case we have a  functor 
\be\label{eq:I-br-mon}
\Inv: \calV^\frakS\lra \calV^{\N}, \quad (X_n)_{n\geq0}  \mapsto (X_n^{\frakS_n})_{n\geq 0}.
\ee
This functor is braided monoidal  because  \eqref{eq:odot} reduces, after taking invariants, to
\eqref{eq:graded-otimes}. It was considered in \cite[Def. 15.5] {aguiar-species-book}  
for the case $\calV=\Vect_\k$.

\begin{ex}\label{ex:set-bialg} (Set bialgebras)
By a {\em set bialgebra} in $\calV$ we will mean a  bialgebra $U=(U_n)_{n\geq 0}$ in $\calV^\frakS$
with $U_0=\1$ whose unit (resp. counit) morphism is the embedding of (resp. projection to) $U_0$. 

More explicitly, a set bialgebra $U$ can be seen as a following set of data:

\begin{itemize}
\item[(i)] A functor  $I\mapsto U_I$ from the groupoid of finite sets to $\calV$, with $U_\emptyset=\1$.

\item[(ii)] For any (ordered) pair $(I,J)$ of finite sets, morphisms
\[
{\mu}_{I,J}: U_{I} \otimes U_{J} \lra U_{I \sqcup J}, \qquad {\Delta}_{I,J}: U_{I \sqcup J} \
\lra U_I \otimes U_J,
\]
\end{itemize}

natural in (isomorphisms of) $I$ and $J$ and 
satisfying the following conditions:

\begin{itemize}

\item[(SB0)]  The maps $\mu$, resp. $\Delta$ are associative, resp. coassociative in the obvious sense.

\item[(SB1-2)] ((Co)unitality) The maps $\mu_{I,\emptyset}, m_{\emptyset, I}, \Delta_{I, \emptyset}, \Delta_{\emptyset,I}$ are all given by unit morphisms or their inverses.

\item[(SB3)] (Exchange) Suppose a finite set $I$  has two disjoint decompositions $I=J_1\sqcup J_2 = L_1\sqcup L_2$.
Put $K_{ij} = J_i\cap L_j$. Then 
the composition 
\[
\xymatrix{U_{J_1} \otimes U_{J_2} \ar[r]^-{\mu_{J_1,J_2}}& U_{I} \ar[r]^-{\Delta_{L_1,L_2}}& U_{L_1} \otimes U_{L_2}}
\]
is equal to the composed morphism
\begin{equation}\label{Eq:setbialg}
\begin{split}
U_{J_1} \otimes U_{J_2} 
&\stackrel{\Delta_{K_{11}, K_{12}} \otimes \Delta_{K_{21}, K_{22}} }{\longrightarrow} 
U_{K_{11}} \otimes U_{K_{12}}
\otimes U_{K_{21}} \otimes U_{K_{22}} \to   \\
&\stackrel{R_{U_{K_{12}},U_{K_{21}}}}{\longrightarrow}
U_{K_{11}} \otimes U_{K_{21}} \otimes U_{K_{12}} \otimes
U_{K_{22}}	
\stackrel{\mu_{K_{11}, K_{21}} \otimes \mu_{K_{12}, K_{22}} } {\longrightarrow} U_{L_1} \otimes U_{L_2}. 
\end{split}
\end{equation}

\end{itemize}

\end{ex}

\begin{rem} As before, this explicit form allows us to speak about set bialgebras in $\calV$ without assuming 
the existence of direct sums. In fact, since the relations above do not involve summation,
one can speak about set bialgebras in not necessary additive categories.
Such structures have been discussed and studied in
\cite{aguiar-species-book} \cite{aguiar:species-Hopf} under the
name {\em Hopf monoids in the category of species}. Our considerations are nominally more general
as we allow species taking values in arbitrary ($\k$-linear) braided categories $\calV$ (the case of the braiding
on $\calV=\Vect_\k^\N$
corresponding to $q$-deformation of commutativity was already considered in 
\cite{aguiar-species-book}  \cite{aguiar:species-Hopf}). 
\end{rem}
 
  Axioms (SB1-2) imply that for any ordered set partition $\bfI=(I_1,\cdots, I_k)$
of $\{1,\cdots, n\}$ we have  well defined iterated multiplication and comultiplication morphisms
\[
\mu_\bfI = \mu_{I_1,\cdots, I_k}: U_{I_1}\otimes\cdots\otimes U_{I_k} \to U_n, 
\quad
\Delta_\bfI = \Delta_{I_1,\cdots, I-k}: U_n\lra U_{I_1}\otimes\cdots\otimes U_{I_k}. 
\]
Further, (SB3) implies, by iteration, the following general exchange rule.

\begin{prop}\label{prop:gen-exch}
In the situation and notation of Proposition \ref{prop:Tits-GL}, the composition
\[
U_{I_1}\otimes\cdots\otimes U_{I_k} \buildrel \mu_\bfI \over\lra U_n\buildrel\Delta_\bfJ \over\lra U_{J_1}\otimes\cdots\otimes
U_{J_l}
\]
is equal to the composition
\[
\bigotimes_i U_{I_i} \buildrel \bigotimes_i \Delta_{K_{i1},\cdots, K_{il}}\over\lra \bigotimes_i \bigotimes_j U_{K_{ij}}
\buildrel R\over\lra  \bigotimes_j \bigotimes_i U_{K_{ij}}
\buildrel \bigotimes_j \mu_{K_{1j},\cdots, K_{kj}}\over\lra \bigotimes_j U_{J_j}, 
\]
where $R$ is the braiding isomorphism in $\calV$ connecting the horizontal (along the rows) and the vertical (along the
columns) orderings of the tensor product $\bigotimes_{i,j} U_{K_{ij}}$. \qed
\end{prop}

Recall the braided monoidal functor $\Inv$  from \eqref{eq:I-br-mon}.

\begin{prop}\label{prop:Inv(U)}
 Assume that $\calV$ is Karoubi-complete. If $U$ is a set bialgebra in $\calV$, then $\Inv(U) = (U_n^{\frakS_n})_{n\geq 0}$
is a graded bialgebra.

\end{prop}

\noindent{\sl Proof:}   If $\calV$ has direct sums, this follows from the fact that the functor $\Inv$ in
\eqref{eq:I-br-mon} is braided monoidal. The general case is reduced to this by passing to $\calV^\oplus$,
the direct sum completion of $\calV$.  
\qed


\subsection{The PROBs of graded bialgebras and of set bialgebras}
Let $\calV, \calW$ be   braided categories. We denote by $\BFun(\calV, \calW)$
the category of $\k$-linear braided monoidal functors $\calV\to\calW$. 
By $\GBialg(\calV)$ and $\SBialg(\calV)$ we denote the categories of graded bialgebras and set bialgebras
in $\calV$ respectively. 

As in  \cite{KS-Prob}, we denote by $\frakB$ the PROB governing graded bialgebras, i.e.,
a braided   category generated by the components of the {\em universal graded bialgebra} $\bfa$. 
It is characterized uniquely up to braided monoidal equivalence by the  (2-)natural equivalence of categories
\be\label{eq:gr-bialg-B}
\GBialg(\calV) \,\= \, \BFun(\frakB, \calV)
\ee
for any   braided   category $\calV$. 

Explicitly, $\frakB$ is generated as a  braided   category, by formal objects $\bfa_n$, $n\geq 0$,
$\bfa_0=\1$ and generating morphisms
\be\label{eq:universal-mu-delta}
\mu_{n_1, n_2}: \bfa_{n_1}\otimes \bfa_{n_2}\lra \bfa_{n_1+n_2}, \quad
\Delta_{n_1, n_2}: \bfa_{n_1+n_2}\lra \bfa_{n_1}\otimes\bfa_{n_2}, 
\ee
which are subject only to  associativity, coassociativity and
the relations (GB1)-(GB3) forming the axioms of a graded bialgebra. 
The identification \eqref{eq:gr-bialg-B} takes $F\in \BFun(\frakB, \calV)$ to the graded bialgebra
$A=(A_n = F(\bfa_n))_{n\geq 0}$ whose multiplication and comultiplication morphisms are
the images under $F$  of \eqref{eq:universal-mu-delta}. 

Thus  objects of $\frakB$ are tensor products $\bfa_\bfn = \bfa_{n_1}\otimes\cdots\otimes \bfa_{n_k}$
for all sequences $\bfn = (n_1,\cdots, n_k)$ of positive integers, as well as $\bfa_0 = \1$. 
Denote $|\bfn|=\sum n_i$. 
Let $\frakB_n\subset\frakB$ be the full subcategory on objects $\bfa_\bfn$ with $|\bfn|=n$.
Then we have a decomposition into blocks
\[
\frakB = \coprod_{n\geq 0} \frakB_n, \text{ i.e., } \Hom_\frakB (\frakB_m, \frakB_n) = 0, \,\, m\neq n.
\]
Note that two sequences $\bfn, \bfn'$ differing by permutation of the $n_i$, represent isomorphic
(by braiding) objects. So isomorphism classes of objects of $\frakB_n$ are labelled by unordered
partitions of $n$. 

Similarly, we denote by $\frakSB$ and call the {\em PROB of set bialgebras} the braided  
category characterized by the (2-)natural equivalence
\be\label{eq:SB_n-charact}
\SBialg(\calV) \,\=\, \BFun(\frakSB, \calV)
\ee
for any braided category $\calV$. Explicitly, $\frakSB$ is generated, as a braided category,
by the components of the {\em universal set bialgebra} $\bfu$, i.e., by formal objects $\bfu_I$,
one for each finite set $I$, with $u_\emptyset=\1$ and by generating morphisms
\[
\begin{gathered}
\gamma_\sigma: \bfu_I \to \bfu_J \text{ for any bijection } \sigma: I\to J, 
\\
\mu_{I,J}: \bfu_I \otimes \bfu_J \lra \bfu_{I\sqcup J}, \qquad \Delta_{I,J}: \bfu_{I \sqcup J} \lra
\bfu_I \otimes \bfu_J
\end{gathered}
\]
subject only to the relations of functoriality ($\gamma_{\sigma\sigma'} = \gamma_\sigma\circ\gamma_{\sigma'}$,
$\gamma_\Id=\Id$), of  naturality of the $\mu_{I,J}$ and $\Delta_{I,J}$  and to  the relations (SB0-3) defining set bialgebras.

Thus objects of $\frakSB$ other than $\1 = \bfu_\emptyset$ have the form 
$\bfu_\bfI = \bfu_{I_1}\otimes\cdots\otimes\bfu_{I_k}$ for all sequences $\bfI = (I_1,\cdots, I_k)$ 
of nonempty finite sets.  As before, we have a decomposition into blocks
\[
\frakSB = \coprod_{n\geq 0} \frakSB_n, 
\]
where $\frakSB_n\subset\frakSB$ is the full subcategory on objects $\bfu_\bfI$
with $\sum |I_\nu|=n$.  

We will also write $\bfu_n$ for $\bfu_{\{1,\cdots, n\}}$. 
Here is the main result of this section.

\begin{theorem}\label{thm:GL}
We have an  equivalences of categories
 $\frakSB_n \= \calQ_n\rtimes\frakS_n$.

\end{theorem}

Here $\calQ_n\rtimes\frakS_n$ is the semidirect product, see \S \ref{subsec:semidir}.

The proof will occupy the remainder of this section. 


\subsection{Braided monoidal structure on $\coprod_n (\calQ_n\rtimes\frakS_n)$} \label{subsec:braid-QS}
To prove Theorem \ref{thm:GL} we need, first of all, construct functors 
$\frakSB_n\to \calQ_n\rtimes\frakS_n$ for all $n$. 
For this it is enough, in view
of \eqref{eq:SB_n-charact}, to make $\coprod_n (\calQ_n\rtimes\frakS_n)$ into a braided category
and exhibit a set bialgebra $U=(U_n)$ in it such that each  $U_n$ lies in the $n$th block. 
To do this, we note the following.

\vskip .2cm

The collection of complex reductive groups of the type 
\[
GL_{n_1} \times \cdots \times GL_{n_r}, \qquad n_1, \ldots, n_r \geq 1
\]
has the remarkable (and often used) property of being stable under the operation of passing to a parabolic subgroup. In particular, the parabolic subgroup 
\[
P_{n,m}=\begin{pmatrix}
GL_n& \star \\ 0 & GL_m
\end{pmatrix}
\]
gives rise to a canonical parabolic correspondence
\[
\xymatrix{GL_n \times GL_m & P_{n,m} \ar[l]_-{p} \ar[r]^-{i} & GL_{m+n}}.
\]
We use this to define a  functor 
\be
\otimes = \otimes_{n,m}:\calQ_n \times \calQ_m \to \calQ_{n+m}
\ee
 as follows. Denote by the same symbols $p,i$ the maps of Lie algebras. This provides an identification of Cartan subalgebras $\frakh_{GL_n} \oplus \frakh_{GL_m} \simeq \frakh_{GL_{n+m}}$. Under this identification, the Coxeter complex of $GL_{n+m}$ is a refinement of that of $GL_n \times GL_m$.
At the level of objects, we define $\otimes_{n,m}$ to be given by the map
\be\label{eq:kappa-nm}
\kappa_{n,m} : \Ob (\calQ_n) \times \Ob(\calQ_m) \to \Ob(\calQ_{n+m}), \qquad  (\frakp_1, \frakp_2) \mapsto \frakp_1\otimes\frakp_2 := p^{-1}(\frakp_1 \oplus \frakp_2).
\ee
In the language of \S \ref{subsec:comb-iden}, objects $\frakp_\bfI$ of $\calQ_n$ are represented by
set partitions $\bfI = (I_1,\cdots, I_k)$ of $\{1,\cdots n\}$, and those of $\calQ_n$ by
set parititons $\bfJ=(J_1,\cdots, J_l)$ of $\{1,\cdots, m\}$. If 
$(\frakp_\bfI,\frakp_\bfJ) \in \Ob(\calQ_n) \times \Ob(\calQ_m)$,
then  $\frakp_\bfI \otimes \frakp_\bfJ = \frakp_{\bfI * \bfJ}$, where
\be\label{eq:p1otimesp2}
\bfI * \bfJ  \, : = \, \bfI \sqcup (\bfJ+n) \, := \, (I_1,\cdots, I_k, J_1+n,\cdots, J_l+n)
\ee
is the concatenation of $\bfI$ and $\bfJ$. 
Here $J_\nu+n$, resp $(\bfJ+n)$ stands for the result of adding $n$ to all integers occuring in  
$J_\nu$ resp. in a part of $\bfJ$.

Let us now define $\otimes_{n,m}$ on morphisms. 
Observe  that  $\kappa_{n,m}$ preserves incidence of parabolic subalgebras in each variable,
i.e.,  the operation $(\bfI, \bfJ) \mapsto \bfI * \bfJ$ preserves refinement of
set partitions in each variable:  if $\bfI$ refines $\bfI'$, then $\bfI * \bfJ$ refines
$\bfI' * \bfJ$ and similarly in the other variable. So we define $\otimes = \otimes_{n,m}$ on
generating morphisms (in each variable) by putting
\be\label{eq:Ind-tensor}
\Ind_{\frakp_\bfI}^{\frakp_{\bfI'}} \otimes \frakp_{\bfJ} \,=\, \Ind_{\frakp_{\bfI * \bfJ}}^{\frakp_{\bfI' * \bfJ}},\quad
\Res_{\frakp_\bfI}^{\frakp_{\bfI'}} \otimes \frakp_{\bfJ} \,=\, \Res_{\frakp_{\bfI * \bfJ}}^{\frakp_{\bfI' * \bfJ}}
\ee
and similarly in the other variable. 

\begin{prop} \label{prop:kappa-nm}
The values of $\otimes_{n,m}$ on generating morphisms defined above,
together with the map  \eqref{eq:kappa-nm} on objects,   define 
a functor $\otimes_{n,m}:\calQ_n\times\calQ_m\to\calQ_{n+m}$. The collection of functors
$\otimes_{n,m}$
endows $\coprod_n \calQ_n$
with the structure of a monoidal category. 
\end{prop}

\noindent{\sl Proof:}  
Notice that $\kappa_{n,m}$ also preserves
collinearity in each variable: if $(\bfI, \bfI', \bfI'')$ are collinear, then so are
$(\bfI *\bfJ, \bfI' * \bfJ, \bfI'' * \bfJ)$ and similarly in the other variable. This implies that for each $\bfJ$
the rule \eqref{eq:Ind-tensor} extends to 
  a functor  $-\otimes \frakp_{\bfJ}: \calQ_n \to \calQ_{n+m}$, and for each $\bfI$ 
 we simillarly get   a functor  $\frakp_\bfI\otimes - : \calQ_m\to \calQ_{n+m}$.  To see that they unite
 into a functor $\otimes_{n,m}: \calQ_n\times\calQ_m\to\calQ_{n+m}$, we need to check the
 commutativity of these single variable tensor products. That is, for each $f\in\Hom_{\calQ_n}(\frakp_{\bfI}, \frakp_{\bfI'})$
 and $g\in\Hom_{\calQ_m}(\frakp_{\bfJ}, \frakp_{\bfJ'})$ the diagram in $\calQ_{n+m}$
 \be\label{eq:f-otimes-g}
 \xymatrix{
 \frakp_\bfI \otimes\frakp_{\bfJ} \ar[r]^{\frakp_\bfI \otimes g}
 \ar[d]_{f\otimes \frakp_\bfJ}
  & \frakp_\bfI\otimes\frakp_{\bfJ'}
  \ar[d]^{f\otimes \frakp_{\bfJ'}}
 \\
 \frakp_{\bfI'}\otimes\frakp_{\bfJ} \ar[r]_{\frakp_{\bfI'}\otimes g} & \frakp_{\bfI'} \otimes\frakp_{\bfJ'} 
 }
 \ee
is commutative, so both composition represent a well-defined morphism $f\otimes g: \frakp_\bfI\otimes\frakp_\bfJ
\to\frakp_{\bfI'}\otimes\frakp_{\bfJ'}$. 

It is sufficient to verify this commutativity when $f$ and $g$ are generating morphisms, i.e., of the form 
$\Ind$ or $\Res$. If they are both of the type $\Ind$, this follows from transitivity of the  morphisms $\Ind$
in $\calQ_{n+m}$ and similarly in the case when they are both of type $\Res$. 
Let us consider the mixed case when
$f=\Ind_{\frakp_\bfI}^{\frakp_{\bfI'}}$, $g= \Res_{\frakp_{\bfJ'}}^{\frakp_\bfJ}$ where $\bfI$ refines $\bfI'$ and $\bfJ'$
refines $\bfJ$. In this case the lower path in the diagram \eqref{eq:f-otimes-g} represents the morphism
$\tau_{\frakp_\bfI \otimes\frakp_\bfJ}^{\frakp_{\bfI'}\otimes\frakp_{\bfJ'}}$ (as defined in the axiom (2c) of
 \S \ref{subsec:double-inc-cat}), since it goes through a common subface of the faces $C_{\bfI * \bfJ}$ and
 $C_{\bfI' * \bfJ'}$. The arrows in the upper path are identified with the corresponding $\tau$-morphisms
 since they go through a common over-face of $C_{\bfI * \bfJ}$ and
 $C_{\bfI' * \bfJ'}$. So the equality of these morphisms follows from collinear transitivity, since 
 $(\frakp_\bfI\otimes \frakp_\bfJ, \frakp_\bfI\otimes\frakp_{\bfJ'}, \frakp_{\bfI'}\otimes\frakp_{\bfJ'})$
 are collinear (the parabolic in the middle corresponds to the common over-face). 
 The other mixed case when $f$ is of type $\Res$ and $g$ is of type $\Ind$, is considered similarly.
 This proves the well-definedness of the functors $\otimes_{n,m}$. The fact that they define a monoidal
 structure (i.e., associativity) is straightforward. \qed

\medskip

Note that the functor $\otimes _{n,m}$ is equivariant with respect  to the symmetric group actions
on the source  and target and the embedding $\frakS_n\times\frakS_m\to\frakS_{n+m}$. 
So it induces a functor which we denote by the same symbol
\[
\otimes_{n,m}: (\calQ_n\rtimes \frakS_n) \times (\calQ_m\rtimes\frakS_m) \lra
\calQ_{n+m}\rtimes \frakS_{n+m}.
\]
Proposition \ref{prop:kappa-nm} implies that it makes $\coprod_{n\geq 0} (\calQ_n\rtimes \frakS_n)$
into a monoidal category.  We further introduce a braiding $R$ on this category as follows. 

Let $(\frakp_\bfI,\frakp_\bfJ) \in \Ob(\calQ_n) \times \Ob(\calQ_m)$,
  By the above,  $\frakp_\bfI \otimes \frakp_\bfJ$ and $\frakp_\bfJ \otimes \frakp_\bfI$
correspond to the ordered set partitions
$\bfI \sqcup (\bfJ+n)$ and $\bfJ \sqcup (\bfI+m)$ of $\{1, \ldots, n+m\}$ respectively.  
Recall \eqref{eq:semidir-Hom} that morphisms in $\calQ_{n+m}\rtimes\frakS_{n+m}$ from $\frakp_\bfI\otimes\frakp_\bfJ$
to $\frakp_\bfJ\otimes\frakp_J$ are formal sums of pairs $(\xi, g)$ where $g\in\frakS_{n+m}$ and 
$\xi$ is a morphism in $\calQ_{n+m}$ from $g(\frakp_\bfI\otimes\frakp_\bfJ)$ to $\frakp_\bfJ\otimes\frakp_\bfI$. 
Let us take $g$ to be the maximal $(n,m)$-shuffle
\be\label{eq:shuffle-sigma}
g=\sigma_{n,m}:= \begin{pmatrix}
1 & \cdots & n & n+1 & \cdots & n+m
\\
m+1&\cdots & m+n& 1&\cdots & m
\end{pmatrix} 
\ee
Then $\frakp' = \sigma_{n,m}(\frakp_\bfI\otimes\frakp_\bfJ)$ is the parabolic corresponding to the ordered set partition
$(\bfI+m) \sqcup \bfJ$ and has the same Levi subalgebra as $\frakp_\bfJ\otimes\frakp_\bfJ$. 
Therefore we have the morphism in the braid groupoid
\[
\xi = \tau_{\frakp'}^{\frakp_\bfJ \otimes \frakp_\bfI}\in\Hom_{\calQ_{n+m}}(\frakp', \frakp_\bfJ\otimes\frakp_\bfI) = 
\Hom_{\calQ_{n+m}}(\sigma_{n,m}(\frakp_\bfI\otimes\frakp_\bfJ), \frakp_\bfJ \otimes\frakp_\bfI)
\]
and we define the braiding operator to be a single pair 
\be\label{eq:braiding-GL}
R_{\frakp_\bfI, \frakp_\bfJ} = (\xi, g) = (\tau_{\frakp'}^{\frakp_\bfJ \otimes \frakp_\bfI}, \sigma_{n,m}) \in
\Hom_{\calQ_{n+m}\rtimes\frakS_{n+m}}(\frakp_\bfI\otimes\frakp_\bfJ, \frakp_\bfJ\otimes\frakp_\bfI). 
\ee

\begin{prop}
This makes $\coprod_n (\calQ_n\rtimes\frakS_n)$ into a braided category. 
\end{prop}

\noindent{\sl Proof:} The naturality of the braiding isomorphisms follows from the functoriality of $\Ind$ and $\Res$
w.r.t. the braiding operators (\S \ref{subsec:anal-DIC}\ref{par:funct-Ind-Res}). The commutativity of the braiding triangles,
namely
\[
\xymatrix{
& \frakp_\bfI\otimes\frakp_\bfK\otimes\frakp_\bfJ 
\ar[dr]^{R_{\frakp_\bfI, \frakp_\bfK} \otimes \frakp_\bfJ}
&
\\
\frakp_\bfI\otimes\frakp_\bfJ\otimes\frakp_\bfK
\ar[ur]^{\frakp_\bfI\otimes R_{\frakp_\bfJ, \frakp_\bfK}}
 \ar[rr]_{R_{\frakp_\bfI\otimes\frakp_\bfJ, \frakp_\bfK}} && 
\frakp_\bfK\otimes\frakp_\bfI\otimes\frakp_\bfJ
}
\]
and the other similar triangle where  $\frakp_\bfI$   is moved past $\frakp_\bfJ \otimes \frakp_\bfK$,
follows from the corresponding standard identities for shuffles \eqref{eq:shuffle-sigma} which  lift into
 identities in the braid group(oid). \qed


\subsection{Proof of Theorem \ref{thm:GL}:  Set bialgebras and the double incidence category}
Let $\QS =\coprod_n(\calQ_n\rtimes\frakS_n)$. 
We will prove a stronger statement implying  Theorem \ref{thm:GL}(a):

\begin{theorem}\label{thm:SB=QS}
$\frakSB$ and $\QS$ are equivalent as braided categories in a way compatible with
decompositions into blocks. 
\end{theorem}

To prove the theorem, we first construct a set bialgebra $U = (U_n)_{n\geq 0}$ in $\QS$.
Recall that $\Ob(\calQ_n\rtimes\frakS_n) = \Ob(\calQ_n) = \calP_{GL_n, T_n}$, the set of
parabolic subalgebras in $\gl_n$ containing $\frakh_n=\Lie(T_n)$. We put $U_n= \gl_n$.
In the language of ordered set partitions of $\{1,\cdots, n\}$, the object $U_n$ corresponds to
the $1$-part set partition consisting of $\{1,\cdots, n\}$ itself. The object $U_n=\gl_n$ is fixed under the $\frakS_n$-action
so  $\frakS_n$ acts on it in $\calQ_n\rtimes\frakS_n$. 
Next, by \eqref{eq:p1otimesp2},
\[
U_n\otimes U_m = \frakp_{(\{1,\cdots, n\}, \{n+1, \cdots, n+m\})}\in \Ob \calQ_{n+m}.
\]
It follows that as an object of $\QS^\oplus$,
\[
U_n\odot U_m = \Ind_{\frakS_n\times\frakS_m}^{\frakS_{n+m}} (U_n\otimes U_m) = \bigoplus_{(J_1, J_2)}
\frakp_{(J_1, J_2)}
\]
where the sum is over all $2$-part ordered set parititons $(J_1, J_2)$ of $\{1,\cdots, n+m\}$ with $|J_1|=n$
and $|J_2|=m$. We now define the mutiplication and comultiplication morphisms in $\QS^\oplus$
\[
\mu_{n,m}: U_n\odot U_m\lra U_{n+m}, \quad \Delta_{n,m}: U_{n+m}\lra U_n\odot U_m
\]
via their components in $\QS$ as follows:
\[
\begin{gathered}
\mu_{n,m} = \bigl( \mu_{J_1, J_2} := \Ind_{\frakp_{(J_1, J_2)}}^{\gl_{n+m}}: \frakp_{(J_1, J_2)} \lra \gl_{n+m}\bigr), 
\\
\Delta_{n,m} = \bigr(\Delta_{J_1, J_2}: = \Res _{\frakp_{(J_1, J_2)}}^{\gl_{n+m}}:  \gl_{n+m}
\lra  \frakp_{(J_1, J_2)}\bigr). 
\end{gathered}
\]

\begin{prop}
These morphisms make $U=(\gl_n)_{n\geq 0}$ into a set bialgebra in $\QS$.
\end{prop}

\noindent{\sl Proof:} The (co)unitality properties are obvious.
Associativity of $\mu$ follows from transitivity of $\Ind$. Coassociativity of $\Delta$
follows from transitivity of $\Res$. To see the exchange property, let $(J_1, J_2)$ and
$(L_1, L_2)$ be two $2$-part ordered partitions of $\{1,\cdots, n+m\}$.  The equality \eqref{Eq:setbialg}
expressing the exchange property in terms of the components, follows directly from
the proto-Langlands formula (Proposition \ref{prop:proto-Lang}), once we take into account the
interpretation  of the Tits product in Proposition \ref{prop:Tits-GL} and the definition  \eqref{eq:braiding-GL} of the
braiding  in $\QS$. \qed

\vskip .2cm

By \eqref{eq:SB_n-charact}, the set bialgebra $U$ defines a braided monoidal functor
\[
F: \frakSB\lra \QS, \quad \bfu_n = \bfu_{\{1,\cdots, n\}}\mapsto U_n=\gl_n,
\]
so that $F(\frakSB_n)\subset \calQ_n\rtimes\frakS_n$ for each $n$. We prove that $F$ is an equivalence by constructing a
quasi-inverse functor $G: \QS \to\frakSB$. An object of $\calQ_n\rtimes\frakS_n$ is a parabolic
$\frakp_\bfI = \frakp_{(I_1,\cdots, I_k)}$  for an ordered set partition $\bfI=(I_1,\cdots, I_k)$
of $\{1,\cdots, n\}$. We define $G$  on objects by
\[
G(\frakp_\bfI) = \bfu_\bfI := \bfu_{I_1}\otimes\cdots\otimes\bfu_{I_k}. 
\]
Next, morphisms in $\calQ_n\rtimes\frakS_n$ are generated by:
\begin{itemize}

\item  Generating morphsims of $\calQ_n$ which are $\Ind_\bfI^{\bfI'}: \frakp_\bfI\to\frakp_\bfI'$
and $\Res_\bfI^{\bfI'}: \frakp_\bfI'\to\frakp_\bfI$ for any pair $(\bfI,\bfI')$ of ordered set partitions of $\{1,\cdots, n\}$
such that $\bfI$ is a refinement of $\bfI'$. 

\item Isomorphisms $\xi_{\sigma,\bfI}: \frakp_I \to\frakp_{\sigma( \bfI)}$ for any
$\sigma\in\frakS_n$ and any $\bfI = (I_1\cdots, I_k)$.
Here $\sigma(\bfI) = (\sigma(I_1),\cdots, \sigma(I_k))$. 

\end{itemize}

We define $G$ on generating morphisms by:
\begin{itemize}
\item  $G(\Ind_\bfI^{\bfI'}) = \mu_{\bfI}^{\bfI'}: \bfu_\bfI\to\bfu_\bfI'$,  the morphism of (iterated) multiplication.

\item $G(\Res^{\bfI'}_\bfI) = \Delta^{\bfI'}_\bfI: \bfu_{\bfI'}\to\bfu_\bfI$, the morphism of (iterated) comultiplication. 

\item $G(\xi_{\sigma, \bfI}) : \bfu_\bfI = \bigotimes_\nu \bfu_{I_\nu} \lra \bigotimes_\nu \bfu_{\sigma(I_\nu)} = 
\bfu_{\sigma(\bfI)}$ is defined as the tensor product of isomorphisms $(\sigma|_{I_\nu})_*: \bfu_{I_\nu}\to\bfu_{\sigma(I_\nu)}$
which $\bfu$, as an object of $\frakSB^\frakS$, associates to the bijection $\sigma|_{I_\nu}$. 
\end{itemize} 
Next, we show that $G$ preserves the relations among the generating morphisms and so is well defined as a
functor $\QS\to\frakSB$. We focus on  relations (listed in \S \ref{subsec:double-inc-cat}, from which
we use their numeration) in $\calQ_n$, as preservation of relations involving the $\xi_{\sigma, \bfI}$ amounts
to $\frakS_n$-equivariance of the rest and is clear from the construction. 

The relation (2a) (identity) is obvious. Transitivity (2b) follows from associativity and coassociativity in $\bfu$. 
Next, preservation of (2c) (idempotency) means that $\Delta_\bfI^{\bf I'} \circ \mu_{\bfI}^{\bfI'} = \Id$. 
By concatenation of the corresponding statements for individual parts of $\bfI'$ this reduces to the cases (for all $n$)
when $\bfI' = (\{1,\cdots,  n\})$ consists of one part.  In such a case the statement follows at once from the exchange
property of
Proposition \ref{prop:gen-exch} with $\bfI = \bfJ$. 

The rest of the relations involve morphisms
\[
T_\bfI^\bfJ = G(\tau_\bfI^\bfJ) = \Delta_\bfJ^\bfK \mu_\bfI^\bfJ
\]
for any two ordered set partitions $\bfI, \bfJ$ of $\{1,\cdots, n\}$ with $\bfK$  being arbitrary such that both $\bfI, \bfJ$
are its refinements. We can and will assume $\bfK=(\{1,\cdots, n\})$  and so $T_\bfI^\bfJ$ is the composition
\[
\bfu_\bfI \buildrel \mu_I \over\lra \bfu_n\buildrel \Delta_J\over\lra \bfu_\bfJ. 
\]
We start with (2e) (invertibility). Let $\bfI,\bfJ',\bfJ''$ be such that $C_{\bfJ'}\supset C_{\bfI} \subset C_{\bfJ''}$ are
codimension $1$ inclusions, with $C_{\bfJ'}, C_{\bfJ''}$ having the same linear span and lying on the
opposite sides of $C_\bfI$. By Proposition  \ref{prop:faces-GL}(c) this means that $\bfJ'=(J'_1,\cdots, J'_l)$
and $\bfJ''=(J''_1,\cdots, J''_l)$ differ by exchange of two adjacent parts, i.e., for some $\nu=1,\cdots, l-1$ we have
\[
\begin{gathered}
J'_a=J''_a, \,\,a\neq \nu,\nu+1, \quad J'_\nu = J''_{\nu+1},\,\, J'_{\nu+1} =J''_\nu, \quad\text{and so}
\\
\bfI = (J'_1, \cdots, J'_{\nu-1}, J'_\nu \cup J'_{\nu+1}, J'_{\nu+2},\cdots, J'_l). 
\end{gathered}
\]
We need to prove that $T_{\bfJ'}^{\bfJ''} = \Delta_{\bfJ''} \mu_{\bfJ'}: \bfu_{\bfJ'}\to \bfu_{\bfJ''}$ is invertible. 
By    Proposition \ref{prop:gen-exch} we see that this composition equals the isomorphism 
$R_{\bfu_{J'_\nu}, \bfu_{J'_{\nu+1}}}$ tensored with the identity operators on the left and on the right and so it
is invertible. 

We now verify (2d) (collinear transitivity) in several stages. First, we consider two ordered set partitions
$\bfJ=(J_1,\cdots, J_l)$ and $\bfJ'=(J'_1,\cdots, J'_l)$ differing by renumbering of parts, i.e., $J'_\nu = J_{s(\nu)}$
for some $s\in \frakS_l$. Recall that the projection $\Br_l\to\frakS_l$ from the braid group to the symmetric one
has a canonical section $s\mapsto \lambda(s)$ characterized by the partial homomorphism property
\be\label{eq:lambda-part-hom}
\lambda(s_1s_2) = \lambda (s_1) \lambda(s_2),\quad \text{if} \quad l(s_1s_2) = l(s_1)+l(s_2)
\ee
and the normalization $\lambda( (\nu, \nu+1))= e_\nu$, the standard generator of $\Br_l$. 
The exchange property (Proposition  \ref{prop:gen-exch}) implies:

\begin{prop}
The morphism 
\[
T_\bfJ^{\bfJ'}:\bfu_{J_1}\otimes\cdots\otimes \bfu_{J_l} \lra \bfu_{J'_1}\otimes\cdots\otimes \bfu_{J'_l} =
\bfu_{J_{s(1)}} \otimes\cdots\otimes \bfu_{J_{s(l)}}
\]
is equal to the braiding isomorphism $R_{\lambda(s)}$. \qed
\end{prop}

\begin{cor}
Let $\bfJ, \bfJ',\bfJ''$ differ from each other by renumbering of parts. Assume that $(C_\bfJ, C_{\bfJ'}, C_{\bfJ''})$
are collinear. Then $T_{\bfJ}^{\bfJ''} = T_{\bfJ'}^{\bfJ''} T_\bfJ^{\bfJ'}$. 
\end{cor}

\noindent {\sl Proof:} This follows from  the description of the braid groupoid by collinearity relations in 
Proposition \ref{prop:pi_1-linear}. Alternatively, let $s_1, s_2\in\frakS_l$ be defined by $J''_\nu = J'_{s_1(\nu)}$,
$J'_\nu = J_{s_2(\nu)}$. Then  collinearity of  $(C_\bfJ, C_{\bfJ'}, C_{\bfJ''})$
is equivalent to the condition $l(s_1s_2) = l(s_1) l(s_2)$, so our statement follows from \eqref{eq:lambda-part-hom}. 
\qed

\vskip .2cm

Finally, in the presense of the relations already established,  general collinear transitivity (2d) is
a consequence of the proto-Langlands formula (Proposition \ref{prop:proto-Lang}) which
expresses $\tau_\bfI^\bfJ\in\Hom_{\calQ_n}(\frakp_\bfI, \frakp_\bfJ)$ as a composition of
``local'' morphisms associated to successive segments of the straight interval
$[x,y]$, $x\in C_\bfI$, $y\in C_\bfJ$. So to prove preservation of (2d) under $G$, 
it suffices to establish a counterpart of the proto-Langlands formula for the
$T_\bfI^\bfJ\in\Hom_\frakSB(\bfu_\bfI, \bfu_\bfJ)$. But such a counterpart is
provided by Proposition \ref{prop:gen-exch}. Theorems \ref{thm:SB=QS} 
and \ref{thm:GL} are proved. 


\subsection{ Graded bialgebras and the P-Coxeter category}

For any ordered partition $\bfn=(n_1,\cdots, n_k)$ of $n\geq 0$ we write
\[
\bfa_\bfn = \bfa_{n_1}\otimes\cdots \otimes \bfa_{n_k}, \quad 
\bfu_\bfn = \bfu_{n_1}\otimes\cdots\otimes\bfu_{n_k}, \quad
\frakS_\bfn = \frakS_{n_1}\times
\cdots\times\frakS_{n_k}. 
\]
Thus $\frakS_\bfn$ acts on $\bfu_\bfn$.

\begin{prop}\label{prop:Inv*}
(a) The object $(\bfu_n^{\frakS_n})_{n\geq 0} \in (\frakSB^\Kar)^\N$
is a graded bialgebra in $\frakSB^\Kar$ and so defines a braided monoidal functor
\[
\Inv^*: \frakB \lra \frakSB^\Kar, \quad \bfa_{\bfn}  \mapsto
\bfu_{\bfn} ^{\frakS_{\bfn}}. 
\]

(b) The image of $\frakB_n$ under $\Inv^*$ is equivalent to $\coxcat_{GL_n}$. 
 \end{prop}
 
 \noindent {\sl Proof} (a) follows from  Proposition \ref{prop:Inv(U)}.
 Part (b) is a consequence of Proposition \ref{prop:inv=kar}(b), since
 $\coxcat_{GL_n}$ is defined as a subcategory of  $\frakS_n$-invariants in $\mathcal{Q}_{GL_n}$,
 while $\frakSB_n$ is identified with $\mathcal{Q}_{GL_n}\rtimes\frakS_n$ by Theorem 
 \ref{thm:GL}. \qed

\begin{rems}
(a) We conjecture that $\Inv^*$ defines an equivalence $\frakB_n\to \coxcat_{GL_n}$. 
This is essentially a particular case of the conjecture of Remark \ref{rem:cox-gen-rel}. 

\vskip .2cm 
(b) Independently of the validity of the conjecture in (a), we have the following
supply of functors from $\coxcat_{GL_n}$. Namely,   each set bialgebra $U$
 in a $\k$-linear braided monoidal abelian category $\calV$
gives a functor $f_n: \frakSB_n\to\calV$ for each $n$ and therefore
a functor $g_n = \Inv^*\circ f_n: \coxcat_{GL_n}\to\calV$. 
 
\end{rems}

\vfill\eject

\section {The P-Coxeter category and perverse sheaves on $\hW$}\label{sec:Perv-hW}

\subsection{Statement of the result.}\label{subsec:stat-result} 

Recall the complex stratification $\calS$ of $\hW$ labelled by $W\bs \calL_{G,T}$, the set of $W$-orbits on
the set of Levi subalgebras. It gives rise to the category $\Perv(\hW)$ of $\calS$-smooth perverse sheaves of $\k$-vector
spaces on $\hW$, where $\k$ is our fixed field of characteristic zero. To recall the definition, for each stratum $S\in\calS$ denote $i_S: S\to\hW$ the embedding.
By definition, $\Perv(\hW)$ is the full subcategory in the derived category $D^b\Sh_\hW$, see \S \ref{subsec:desc-sh},
formed by complexes of sheaves $\calF$ such that for each $S\in\calS$ the three following conditions are satisfied:
\begin{itemize}
\item[($\Perv^0$)] All $\ul H^p (i_S^*\calF)$, the cohomology sheaves of the complex $i_S^*\calF$,
are locally constant.
\end{itemize}

This condition (for any $S\in\calS$) is referred to as {\em $\calS$-smoothness} of $\calF$. It implies that
 all $\ul H^p (i_S^!\calF)$ are locally constant as well. The other two conditions (perversity)
 are imposed on these two types
 of sheaves:

  \begin{itemize}

\item[($\Perv^-$)] $\ul H^p (i_S^*\calF)=0$ for $p>-\dim_\C S$. 

\item[($\Perv^+$)] $\ul H^p (i_S^!\calF)=0$ for $p< -\dim_\C S$. 
\end{itemize}

It is known that $\Perv(\hW)$ is an abelian category, being the heart of a t-structure in the category of
$\calS$-smooth complexes. 

\begin{rems}\label{rems:def-perv}
(a) Unlike some standard treatments, we do not require that the stalks of the $\ul H^p (i_S^*\calF)$ are finite-dimensional.
In fact, one can consider $\calS$-smooth perverse sheaves with values in any abelian category $\calV$ by realizing
them insider the category of complexes of cellular sheaves for an appropriate cell decomposition refining $\calS$,
e.g., as in \cite{K-S-hW} and understanding cellular sheaves as functors from the poset of cells to $\calV$,
cf.   \cite{Vybornov} or \cite[Prop. 1.8]{K-S-realhyper}. 
In such an approach, Verdier duality involves passing to the opposite category of $\calV$
which need not be identified with $\calV$. 

\vskip .2cm

(b) The  normalization of the perversity conditions used here, differs by shift from that of
\cite{K-S-realhyper}.  
\flushright{$\triangle$}
\end{rems}

Here is the main result of this section.

\begin{theorem}\label{thm:Perv-hW}
We have an equivalence of categories $\Perv(\hW)\= \Fun(\coxcat, \Vect_\k)$. 
\end{theorem}

The proof will be given in the next two subsections. 

\begin{rem}
A different description of $\Perv(\hW)$, involving more redundant data, was given in 
\cite{K-S-hW}. The advantage of  Theorem \ref{thm:Perv-hW}  is its direct parallelism
with the algebraic properties of induction and restriction operations in representation theory,
in particular, with the Langlands formula. 
\flushright{$\triangle$}
\end{rem}

\subsection{Comparison of descent for modules and perverse sheaves} We start with $\calS_\frakh$,
the stratification of $\frakh$ into generic strata of complex flats, see \S \ref{subsec:rest-root}. It gives rise
to the category $\Perv(\frakh)$ of $\calS_\frakh$-smooth perverse sheaves of $\k$-vector spaces of 
$\frakh$, defined similarly to \S \ref{subsec:stat-result}. The main result of \cite[Th. 8.1]{K-S-realhyper}
can be viewed as an equivalence
\be\label{eq:Perv-h}
\Perv(\frakh) \= \Fun(\calQ,\Vect_\k).
\ee
Our proof of Theorem  \ref{thm:Perv-hW} will proceed from this by an appropriate
descent procedure. 

More precisely, both categories in \eqref{eq:Perv-h} are equipped with natural $W$-actions; the
equivalence is, by construction, compatible with the actions  and so gives an equivalence of
equivariant categories
\be\label{eq:Perv-h-W}
\Perv[\hW] \buildrel {\on{def}}\over = \Perv(\frakh)^{[W]} \= \Fun(\calQ, \Vect_\k)^{[W]}. 
\ee
Here $[\hW]$ is the orbifold (stack) quotient so the first equality in \eqref{eq:Perv-h-W} is just a
different notation for the same category. 

Let $A=A(\calQ)$ be the total algebra of the category $\calQ$, see \eqref{eq:A(V)}. The $W$-action on $\calQ$
gives an action on $A$ and we have
\[
\Fun(\calQ,\Vect_\k)^{[W]} \= _{A[W]}\Mod, \quad \Fun(\coxcat, \Vect_\k) \= _{A^W}\Mod,
\]
as $\coxcat = \calQ^W$ and so $A(\coxcat) = A^W$, see \S \ref{subsec:orb-cat-II}. We now consider the
diagram
\be\label{eq:diag-compar}
\xymatrix{
\Perv(\hW) \ar@{-->}[d]
\ar@<0.1cm>[r]^{p^{-1}} & \Perv[\hW] \ar@<0.1cm>[l]^{p_*} \ar[d]_\mu
& \calK_\top \ar@{^{(}->}[l]
\\
{A^W}\Mod  \ar@<0.1cm>[r]^{I} & _{A[W]}\Mod  \ar@<0.1cm>[l]^{R}   & \calK_\alg \ar@{^{(}->}[l],
} 
\ee
where:
\begin{itemize}
\item $(p^{-1}, p_*)$ is the adjoint pair of the pullback and pushforward functors induced by the projection
$p:[\hW]\to \hW$ of the stacky quotient to its coarse moduli space, which is in this case the topological and GIT quotient. Explicitly, these functors
are given by \eqref{eq:p^{-1}}, \eqref{eq:p_*} for single sheaves and are exact (on the abelian categories
of sheaves) since $q: \frakh\to\hW$ is a finite  flat morphism. 
They extend componentwise to the derived categories and further preserve the subcategories of perverse sheaves
(for the same reason that $q$ is finite). Since short exact sequences in the abelian categories of perverse
sheaves come from exact triangles in the ambient triangulated categories, both $p^{-1}$ and $p_*$ in
\eqref{eq:diag-compar} are exact.

\item $(I,R)$ is the adjoint pair of the induction and restriction functors \eqref{eq:I,R-general}. Since $R$
consists of taking $W$-invariants \eqref{eq:I,R-explicit}, it is an exact functor (whereas $I$ is, a priori, only
right exact).

\item $\calK_\top$ resp. $\calK_\alg$ is the full subcategory in $\Perv[\hW]$ resp. $_{A[W]}\Mod$
consisting of objects annihilated by $p_*$ resp. $R$.

\item $\mu$ is the equivalence \eqref{eq:Perv-h-W}. 
\end{itemize}

We want to construct an equivalence depicted by the dotted arrow on the left in a way compatible
with the other data in the diagram. Corollary \ref{cor:der-unit} implies (by restricting to $\Perv(\hW)$) that
the unit of the adjunction $(p^{-1}, p^*)$ in \eqref{eq:diag-compar} is an isomorphism, and 
Proposition \ref{prop:unit-RI} gives a similar statement for the unit of $(I,R)$. 
Therefore, by Proposition \ref{prop:W/K=V}, $\calK_\top$ and $\calK_\alg$ are Serre subcategories with quotients
\[
\Perv[\hW]/\calK_\top \= \Perv(\hW), \quad _{A[W]}\Mod/\calK_\alg \= _{A^W}\Mod. 
\]
So our task will be accomplished once we prove the following.

\begin{prop}\label{prop:mu-Ktop-Kalg}
The equivalence $\mu$ of  \eqref{eq:Perv-h-W} and \eqref{eq:diag-compar} takes
$\calK_\top$ into $\calK_\alg$.
\end{prop}

This will be done in the next subsection. 

\subsection{Comparison of kernels} Let $\calF\in\Perv(\frakh)$ and $F: \calQ\to\Vect_\k$ be the functor
(``double quiver'' in the terminology of \cite{K-S-realhyper}) associated to $\calF$ by the equivalence
\eqref{eq:Perv-h}. Explicitly  \cite{K-S-realhyper}, for an object $\frakp\in\Ob\calQ = \calP_{G,T}$ the value
$F(\frakp)$ is the stalk at the face $C^\circ_\frakp\subset \frakh_\R$ of the sheaf $\ul\H^0_{\frakh_\R}(\calF)$.
The appearance of $\ul\H^0$ and not $\ul\H^{\dim \frakh_\R}$, as in \cite{K-S-realhyper} is due to the
different normalization of the perversity conditions, see Remark \ref{rems:def-perv}(b). 

Further, let $i:\frakh_\R\to\frakh$ be the embedding. It was shown in  \cite{K-S-realhyper}  that the
complex $i^!\calF$ of which $\ul\H^0_{\frakh_\R}(\calF)$ is the $0$th cohomology sheaf,  has no other
cohomology sheaves, i.e., is quasi-isomorphic to  $\ul\H^0_{\frakh_\R}(\calF)$. 

Thuss, if $\calF$ is $W$-equivariant, i.e., $\calF\in\Perv[\hW]$, then the $A[W]$-module $\mu(\calF)$
is, by \eqref{eq:mod-functor}, the direct sum
\[
\mu(\calF) = \bigoplus_{\frakp\in\calP_{G,T}} F(\frakp) \,=\, \bigoplus_{\frakp\in\calP_{G,T}} (i^!\calF)_{C^\circ_\frakp}.
\]
So 
\be\label{eq:FinK-top}
\calF\in\calK_\top \quad \Leftrightarrow \quad (q_*\calF)^W=0, 
\ee
while
\be\label {eq:F-in-Kalg}
\mu(\calF)\in\calK_\alg \quad \Leftrightarrow \quad 
\biggl(\bigoplus_{\frakp\in\calP_{G,T}} (i^!\calF)_{C^\circ_\frakp}\biggr)^W =0,
\ee
and Proposition \ref{prop:mu-Ktop-Kalg} means that the conditions in the RHS of 
\eqref{eq:FinK-top} and \eqref{eq:F-in-Kalg} are equivalent. To prove this, consider the Cartesian diagram
\be\label{eq:diagr-h-hR}
\xymatrix{
\frakh_\R \ar[d]_{q_\R}
\ar[r]^i &  \frakh \ar[d]^q
\\
\hW_\R \ar[r]_j &\hW
}
\ee
with $i,j$ being closed embeddings and $q_\R$ being induced by $q$. The group $W$ acts on the top row.
By base change and exactness of taking $W$-invariants, 
\[
(q_{\R*} i^!\calF)^W \= (j^!q_*\calF)^W \= j^!((q_*\calF)^W) = j^!p_*\calF. 
\]
Now, $\hW_\R$ is identified with the closure of the positive Weyl chamber. Its faces are in bijection with
$W$-orbits on  faces of the arrangement $\calH$ in $\frakh$, i.e., with objects of the P-Coxeter category
$\coxcat = \calQ^W$. So
\[
M^W = \bigoplus_{C\subset \hW_\R} (q_{\R*}i^!\calF)^W_C \=
\bigoplus_{C\subset \hW_\R} (j^!p_*\calF)_C,
\]
the direct sum of the stalks of $j^!p_*\calF$ over all the faces in $\hW_\R$. We conclude that
\[
M^W=0 \quad \Leftrightarrow\quad j^! p_*\calF = 0. 
\]
This automatically implies that $\mu(\calK_\top)\subset\calK_\alg$, i.e., that $M^W=0$ whenever $p_*\calF=0$.
The reverse inclusion comes from the following

\begin{prop}\label{prop:i^G=0}
Let $\calG\in\Perv(\hW)$ be such that  $j^!\calG=0$. Then $\calG=0$. 
\end{prop}

\noindent{\sl Proof:} For any $m\geq 0$ let $\hW^{\geq m}\subset\hW$ be the union of the complex
strata of complex codimension $\geq m$ and 
\[
\hW^{\geq m}_\R = (\hW_\R)\cap (\hW^{\geq m}) \,\,\subset \,\, \hW_\R
\]
be the union of the real
strata (faces) of real codimension $\geq m$. Denote
$j_{\geq m}: \hW^{\geq m}_\R \to    \hW^{\geq m} $ the embedding. Let
\[
j_m: \hW^{=m}_\R \,\,:=\,\,  (\hW_\R^{\geq m}) - (\hW_\R^{\geq m+1}) \lra \hW^{=m} := \,\,
(\hW^{\geq m}) - (\hW^{\geq m+1})
\]
be the emdedding of the union of real strata of real odimension exactly $m$ into the union of complex
strata of complex codimension exactly $m$. Note that $\hW^{=m}$ is a complex manifold and
$j_m$ is an embedding of a totally real submanifold of middle codimension which has nonempty
intersection with each connected component of $\hW^{=m}$. This implies:

\begin{lemma}
If $\calN$ is a (shifted) local system on $\hW^{=m}$ such that $j_m^!\calN=0$, then $\calN=0$.
\qed
\end{lemma}

We now prove by increasing induction on $m$ that  the restriction of $\calG$ on
the open set $(\hW)-(\hW)^{\geq m}$ is zero. Indeed, for $m=1$  we are talking about
the  restriction  $\calG|_{\hW^{=0}}$. This is a shifted local system, denote it $\calN$.
Since $j^!\calG=0$, we have $j_0^!\calN=0$ and so $\calN=0$. So $\calG$ is supported
on $\hW^{\geq 1}$. Therefore $\calN_1 = \calG|_{\hW^{=1}}$ is a shifted local system
satisfying $j_1^!\calN_1=0$, so $\calN_1=0$ and $\calG$ is supported on $\hW^{\geq 2}$.
Continuing like this, we find that $\calG=0$.

\vskip .2cm

This finishes the proof of Proposition \ref{prop:i^G=0} and thus of Proposition \ref{prop:mu-Ktop-Kalg} and Theorem
\ref{thm:Perv-hW}.

\vfill\eject

 \appendix
 
 \section{Orbit categories and descent}
 
 \subsection{Coreflexive subcategories and localization}
 Let $\calV, \calW$ be categories, $i: \calV\to\calW$ a functor, having a right adjoint
 $\rho: \calW\to\calV$. Recall the well known
 
 \begin{prop}\label{prop:coreflect}
 The following are equivalent:
 \begin{itemize}
 \item[(i)] $i$ is fully faithful, i.e., embeds $\calV$ as a full subcategory in $\calW$.
 
 \item[(ii)] The unit of the adjunction $\eps: \Id_\calV\to \rho\circ i$ is an isomorphism of functors.
 \end{itemize}
 
 \end{prop}
 
 \noindent{\sl Proof:} This is a part of (the dual version of) \cite[Prop. I.1.3]{GZ}. For example,
 (ii)$\Rightarrow$(i) follows by writing
 \[
 \Hom_\calW(i(V), i(V')) \= \Hom_\calV (V,\rho i(V')) \= \Hom_\calV (V, V'). \qed
 \]
 
 \vskip .2cm
 
 A full subcategory whose embedding functor has a right adjoint, is called {\em coreflective},
 see \cite[Def.4.5.12]{Riehl} \cite[\S I.3.5]{Borceux}. The following is also standard
 \cite[Prop.I.1.3(iii)]{GZ}.
 
 \begin{prop}\label{prop:localiz}
 In the situation of Proposition \ref{prop:coreflect}, let $\Sigma_\rho\subset\Mor(\calW)$
 be the class of morphisms $u$ such that $\rho(u)$ is invertible and denote,
 as usual, by $\calW[\Sigma_\rho^{-1}]$
  the localization of $\calW$  along $\Sigma_\rho$. Then $\rho$ induces an equivalence
 $\calW[\Sigma_\rho^{-1}]\to\calV$. \qed
 \end{prop}
 
 Now assume, in the addition to the situation of Propositions \ref{prop:coreflect} and
 \ref{prop:localiz} that $\calV,\calW$ are abelian and $\rho$ is exact ($i$ is always right exact
 and $\rho$ is left exact as they form an adjoint pair). 
 Let $\calK\subset\calW$ be the ``kernel'' of $\rho$, i.e., the full subcategory formed by objects
 $K\in\calW$ such that $\rho(K)=0$. 
 
 \begin{prop}\label{prop:W/K=V}
 Under the above assumptions, $\calK$ is a Serre subcategory in $\calW$ and the quotient category
 $\calW/\calK$ is canonically equivalent to $\calV$.
  \end{prop}
  
  \noindent{\sl Proof:} Being a Serre category means closure under 
  passing to subobjects, quotient objects and
  extensions which is true since $\rho$ is exact. Next, by definition $\calW/\calK = \calW[\Sigma_\calK^{-1}]$,
  where $\Sigma_\calK$ consists of morphisms $u$ suct that $\Ker(u)$ and $\Coker(u)$ are in $\calK$. 
  But because $\rho$ is exact, $\rho(\Ker(u)) = \Ker(\rho(u))$ and $\rho(\Coker(u))=\Coker(\rho(u))$,
  so $\Sigma_\calK = \Sigma_\rho$. The identification $\calW/\calK \=\calV$ now follows from
  Proposition \ref{prop:localiz}.\qed

 
 \subsection{Group actions on categories}\label{subsec:gr-act-cat}
 Let $G$ be a group and $\calV$ a category. Recall that a {\em  (weak) $G$-action}   of $G$ on $\calV$ is a datum $\gamma$ of:
 \begin{itemize}
 \item [(0)] A functor $\gamma_g: \calV\to\calV$ given for each $g\in G$ such that $\gamma_1=\Id$.
 
 \item[(1)] A natural isomorphism of functor $\gamma_{g,h}: \gamma_g \gamma_h \to \gamma_{gh}$
  given for each $g,h\in G$ such that:
 
 \item[(2)] The $2$-cocycle condition holds, i.e., for  any $g,h,k\in G$ the diagram of natural transformations
 \[
 \xymatrix{
 \gamma_g \gamma_h \gamma_k \ar[r]^{\gamma_{g,h}\circ \gamma_k}
 \ar[d]_{\gamma_g\circ \gamma_{h,k}}  & \gamma_g \gamma_{hk}
 \ar[d]^{\gamma_{g, hk}}
 \\
 \gamma_g \gamma_{hk} \ar[r]_{\gamma_{g, hk}} & \gamma_{ghk}
 }
 \]
 is commutative. 
  \end{itemize}
  In this case we say that $\calV$ is a $G$-{\em category}. 
  
  By a {\em $G$-equivariant object} of a $G$-category $\calV$ we mean an object $V\in\calV$
  together with isomorphisms $\phi_g=\phi_{g,V}: V\to \gamma_g(V)$ given for each $g\in G$ and such that
  for any $g,h\in G$ the diagram
  \[
  \xymatrix{
  V\ar[r]^{\phi_g} 
  \ar[d]_{\phi_{gh}}& \gamma_g(V)\ar[d]^{\gamma_g(\phi_h)}
  \\
  \rho_{gh}(V) \ar[r]_{\gamma_{g,h,V}} &\gamma_g\gamma_h(V)
  }
  \]
  commutes. $G$-equivariant objects in $\calV$ form a category which we denote $\calV^{[G]}$. 

A $G$-action $\gamma$ is called {\em strict} (and we say that $\calV$ is a {\em strict $G$-category}),
if all $\gamma_{g,h}$ are identity transformations, in particular,
$\gamma_g \gamma_h = \gamma_{gh}$. This concept is used mostly in the case when $\calV$ is a small category.

\begin{exas}\label{ex:G-equiv-sh-mod}
(a) Let $X$ be a topological space (see \S \ref{sec:notation} for conventions) By $\Sh_X$ we denote the
category of sheaves of $\k$-vector spaces on $X$. If $G$ acts on $X$ in the usual sense
then $\Sh_X$ becomes a $G$-category. The category $\Sh_X^{[G]}$ consists of $G$-equivariant
sheaves in the standard sense. 

\vskip .2cm

(b) Let $A$ be an associative $\k$-algebra with unit. By $_A\Mod$ we denote the category
of unital  left $A$-modules. For a non necessarily unital $A$-module $M$ we denote
\be\label{eq:M^u}
M^u = \{ m\in M\, | \, 1\cdot m=m\}
\ee
the maximal unital submodule of $M$.
If $G$ acts on $A$ in the usual sense
(notation $(g, a) \mapsto g(a)$),
 then $_A\Mod$ becomes a $G$-category. 
The category $_A\Mod^{[G]}$ consists of $G$-equivariant $A$-modules. Alternatively, let $A[G]$
be the {\em twisted group algebra} of $G$ acting on $A$. By definition, it consists of  formal sums
$\sum_{g\in G} a_g[g]$ where $a_g\in A$ are zero for almost all $g$, and the formal symbol $[g]$
satisfies the commutation rule $[g] a = g(a) [g]$. Then  a $G$-equivariant $A$-module is the same
as an $A[G]$-module, i.e., $_A\Mod^{[G]}\,  \simeq\,  _{A[G]}\Mod$. 

\flushright{$\triangle$}
\end{exas}


\subsection{Orbit categories I. Coinvariants}\label{subsec:coinv-orbits}
 Let $\calV$ be a small category with a strict action of $G$.
Thus $G$ acts on the set $\Ob\calV$. Let us assume that this action is free. In this case we can define
the {\em coinvariant orbit category} $\calV_G$ as follows. We put $\Ob \calV_G = G\backslash \Ob\calV$,
the set of $G$-orbits $O\subset \Ob\calV$ and for two such orbits $O_1, O_2$ put
\be\label{eq:Hom-V[G]} 
\Hom_{\calV_G}(O_1, O_2) = G\big\backslash \coprod_{x_1\in O_1, x_2\in O_2} \Hom_\calV (x_1, x_2),
\ee
the set of $G$-orbits on the disjoint union.  

If $O_1, O_2, O_3\subset\Ob\calV$ are three orbits, 
$u\in\Hom_{\calV_G}(O_1, O_2)$ is represented by $\wt u \in\Hom_{\calV}(x_1, x_2)$, $x_i\in O_i$
and $v\in\Hom_{\calV_G}(O_2, O_3)$ is represented by $\wt v\in \Hom_{\calV}(x'_2, x'_3)$, $x'_i\in O_i$,
then there is a unique $g\in G$ such that $g(x_2) = x'_2$ and $vu\in\Hom_{\calV_G}(O_1, O_3)$
is represented by $\wt v \circ g(\wt u)\in \Hom_\calV (g(x_1), x_3)$. Here it is essential that the $G$-action
on $\Ob\calV$ is free. See \cite[\S2]{ciblis} for a linear version of this construction. 

Note that we have a canonical functor $\xi: \calV \to\calV_G$. 

\begin{ex}\label{ex:fund-grpd-coinv}
 $p: X\to Y$ be an unramified Galois covering  of topological spaces
 with Galois group $G$, so $G$ acts freely on $X$ and $Y=G\backslash X$. 
Let $B\subset Y$ be a subset and $\Pi_1(Y,B)$ be the fundamental groupoid  of $Y$ with the set  of base points $B$, so
 $\Ob \Pi_1(Y,B)=B$. Let $\wt B=p^{-1}(B)$ and $\calV=\Pi_1(X,\wt B)$ be the fundamental groupoid of  $X$ with
  set of base points $\wt B$. Then $G$ acts strictly on $\calV$  the action on $\Ob\calV = \wt B$ is free with
  set of orbits $B$. In this case we have an isomorphism of groupoids $\Pi_1(X,\wt B)_G \simeq \Pi_1(Y,B)$, and the functor $\xi$ is identified with $p_*$, the morphism
  of fundamental groupoids induced by $p$. 
  \flushright{$\triangle$}
\end{ex}


\subsection{Karoubi completions and objects of invariants}\label{subsec:karoubi}
Let $\calV$ be a $\k$-linear category.  We denote by:

\begin{itemize}
\item $\calV^\oplus$ the finite direct sum completion of $\calV$. Its objects are formal finite direct
sums $\bigoplus_i V_i$ of objects of $\calV$ and morphisms are given by matrices.
If $\calV$ already has finite direct sums, then $\calV^\oplus\=\calV$. 

\item $\calV^\Kar$ the {\em Karoubi completion} of $\calV$ obtained by adjoining the images of projectors.
Thus objects of $\calV^\Kar$ are formal symbols $\Im(p)$ for all  projectors, i.e., endomorphisms 
$p\in\End(V)$, $V\in\Ob(\calV)$ such $p^2=p$. If $q\in\End(W)$  is another projector, then
\[
\Hom_{\calV^\Kar}(\Im(p), \Im(q)) \,=\,  q \Hom_\calV(V,W) p \, \subset\,  \Hom_\calV(V,W). 
\]
We say that $\calV$ is {\em Karoubi complete}, if $\calV\hookrightarrow \calV^\Kar$, $V\mapsto \Im(\Id_V)$,
is an equivalence, i.e., each projector $p: V\to V$ is realized as the projection on a direct summand. 

\item $\calV^{\oplus/2} = (\calV^\oplus)^\Kar$. Note that it is both Karoubi complete and closed
under direct sums as $\Im(p)\oplus\Im(q) = \Im(p\oplus q)$. 
\end{itemize}
If $\calV$ has a $\k$-bilinear (braided) monoidal structure $\otimes$, then $\otimes$ is distributive
with respect to direct sums in $\calV$ whenever they exist. The (braided) monoidal structure on $\calV$
is extended to $\calV^\oplus$ by distributivity and to $\calV^\Kar$ by 
$\Im(p)\otimes\Im(q) = \Im(p\otimes q)$, and therefore extends to $\calV^{\oplus/2}$ as well. 

\vskip .2cm

Let $G$ be a finite group. Then $\eps = {1\over|G|} \sum_{g\in G} [g]\in\k[G]$ is an idempotent:
$\eps^2=\eps$. If $G$ acts on an object $V\in \calV$, then the action of $\eps$ on $V$ is
a projector $\eps_V: V\to V$. We denote by $V^G$ and call
the {\em object of invariants} the object $\Im(\eps_V)\in\calV^\Kar$. 

Assume that $(\calV,\otimes)$ is monoidal. If $G_i$ acts on $V_i$, $i=1,2$, then
$(V_1\otimes V_2)^{G_1\times G_2} \= V_1^{G_1} \otimes V_2^{G_2}$. 


\subsection{The semidirect product category}\label{subsec:semidir}
Let $\calV$ be a $\k$-linear cateory  with finitely many objects, 
  equipped with a strict action of a finite group $G$. 
 The {\em semi-direct product}  $\calV \rtimes G$ is the category with the same set of objects as $\calV$ and 
 morphisms given by
\be\label{eq:semidir-Hom}
\Hom_{\calV \rtimes G}(V,W):=\biggl\{\sum_{g\in G}(\xi,g)\;|\; \xi \in \Hom_\calV(\gamma_g(V),W) 
\biggr\}, \quad V,W\in \Ob(\calV). 
\ee
  For each given $g$, the symbols $(\xi,g)$ are assumed to be $\k$-linear in $\xi$. 
The composition of morphisms is given by
\begin{equation}\label{eq:semidirectgroupoid}
	(\xi,g) \circ (\zeta,h)=( \xi \circ \gamma_g(\zeta), gh)
\end{equation}
whenever this makes sense. 
Intuitively, $\calV\rtimes G$ is the category obtained from $\calV$ by adding isomorphisms 
$(\Id,g): V \to \gamma_g(V)$. 

\begin{exas}\label{exas:semidir}
 (a) Suppose $\calV$ has one object $V$, i.e., reduces to an associative algebra $A=\End_\calV(V)$
with $G$-action. 
Then $\calV\rtimes G$ also has one object and its endomorphism algebra is
$\End_{\calV\rtimes G}(V) = A[G]$, the twisted group algebra of $G$, see Example 
\ref{ex:G-equiv-sh-mod}(b).  

\vskip .2cm

(b) Generalizing the  property cited in Example 
\ref{ex:G-equiv-sh-mod}(b), we have  an identification
\[
\Fun(\calV, \calW)^{[G]} \simeq \Fun(\calV\rtimes G,\calW)
\]
for any category $\calW$ equipped with a trivial $G$-action.
\end{exas}

The  relation to $\calV\rtimes G$ with  twisted group algebras can be generalized as follows.  Denote by  
$\calV^{\oplus}$ the category obtained from $\calV$ by adding formal direct sums. 
We can assign to $\calV$ an associative $\k$-algebra called the {\em total algebra} of $\calV$,  defined as
\be\label{eq:A(V)}
A(\calV) = \bigoplus_{V_1,V_2\in\Ob\calV} \Hom_{\calV}(V_1, V_2) = \End_{\calV^\oplus} \biggl(\bigoplus_{V\in\Ob\calV}
V\biggr). 
\ee
This algebra has a unit which is decomposed into a sum of orthogonal idempotents 
$1=\sum_{V\in\Ob\calV}\eps_V$, where the idempotent $\eps_V$ corresponds to $\Id_V$. We have an equivalence of
categories
\[
_{A(\calV)}\Mod \= \Fun(\calV, \Vect_\k).
\]
Expilcitly, given a functor $F: \calV\to\Vect_\k$, we associate to it the $A(\calV)$-module 
\be \label{eq:mod-functor}
M = \bigoplus_{V\in\Ob\calV} F(V)
\ee
and  $F$ is recovered from $M$ by $F(V) = \eps_V M$. 

If now  $G$ acts on $\calV$ in a strict sense then it also acts on $A(\calV)$  
and so $A(\calV)[G]$ is defined. Then we have an isomorphism of algebras and
by applying Example   \ref{exas:semidir}(b), an equivalence of categories
\[
A(\calV\rtimes G) \= A(\calV)[G], \quad 
_{A(\calV)[G]}\Mod \= \Fun(\calV\rtimes G, \Vect_\k)
\]

Further, let $V\in\calV$ be a $G$-equivariant object with structure isomorphisms 
$\phi_g: V\buildrel\sim\over\to \gamma_g(V)$. Then $\lambda_g = (\phi_g^{-1}, g)$ is an automorphism of $V$ in
$\calV\rtimes G$. In this way we get an action of $G$ on $V$ by automorphisms of $\calV\rtimes G$.
Note that an object with $G$-action is the same as a $G$-equivariant object with respect to the trivial
action of $G$ on the category. 

\begin{prop}\label{V^G=VXG}
The above construction defines an equivalence of categories 
\[
\calV^{[G]} \lra (\calV\rtimes G)^{[G]},
\]
where the $G$-action of $G$ on $\calV\rtimes G$ is trivial.  \qed
\end{prop}

\vfill\eject

\subsection{Orbit categories II. Invariants} \label{subsec:orb-cat-II}
Let $\calV$ be a $\k$-linear category with finitely many objects and let $G$ be a finite group acting on $\calV$ in a strict sense.
We define a new category $\calV^G$ called the {\em invariant orbit category} of $\calV$ under the $G$-action.
By definition, $\Ob \calV^G = G\backslash \Ob\calV$ is the set of $G$-orbits on $\Ob\calV$. For two such orbits $O_1, O_2$ we put
 \be\label{eq:orbit-cat-inv-hom}
 \Hom_{\calV^G}(O_1, O_2) = \biggl(\bigoplus_{V_1\in O_1, V_2\in O_2} \Hom_\calV(V_1, V_2)\biggr)^G,
 \ee
the space of $G$-invariants in the direct sum, as opposed to coinvariants in \eqref{eq:Hom-V[G]}. 
The composition of morphisms in $\calV^G$ is defined as follows :
if $O_1, O_2, O_3$ are $G$-orbits in $\Ob\calV$ and if $a=(a_{V_1, V_2})_{V_1,V_2} \in \Hom_{\calV^G}(O_1, O_2)$, $b=(b_{V_2, V_3})_{V_2, V_3} \in \Hom_{\calV^G}(O_2, O_3)$ then $b \circ a$ is defined as
\[
(b \circ a)_{V_1,V_3}=\sum_{V_2 \in O_2} b_{V_2, V_3} a_{V_1, V_2} \in \Hom_{\calV }(V_1, V_3).
\]

The construction of the invariant orbit category also has an interpretation in terms of total algebra $A(\calV)$. We can form the invariant subalgebra $A(\calV)^G$. For any $G$-orbit $O\subset \Ob\calV$ the
 element $\eps_O=  \sum_{V\in O} \eps_V$ is an idempotent in $A(\calV)^G$; the unit $1 \in A(\calV)$ 
 lies in $A(\calV)^G$ and 
 is the sum $\sum_{O\in G\backslash\Ob\calV} \eps_O$. 
 It follows that 
 \be
 A(\calV)^G = A(\calV^G).
 \ee

 For a fixed pair $(V_1, V_2)$ of objects of $\calV$ there is a canonical
 {\em averaging  map} 
\be\label{eq:map-xi}
\xi_{V_1, V_2}: \Hom_{\calV}(V_1, V_2) \to  \Hom_{\calV^G}(G V_1, G V_2)
\ee
 which assigns to a morphism $a \in \Hom_{\calV}(V_1, V_2)$  the element
 \be\label{eq:def-average}
\ul a =\xi_{V_1, V_2}(a) := {1\over |\Stab(V_1, V_2)|}   \sum_{g\in G} \bigl\{ ga: gV_1\to gV_2\bigr\}\quad 
 \in \bigoplus_{V'_1\in GV_1, V'_2\in GV_2} \Hom_\calV(V'_1, V'_2). 
 \ee
 Here:
 \begin{itemize}
 
 \item $\Stab(V_1, V_2) \subset G$ is the stabilizer of the pair 
 $(V_1, V_2)\subset \Ob(\calV)\times\Ob(\calV)$. 
 
 \item  $ga: gV_1\to gV_2$ is considered as an element of the direct sum
 on the right in \eqref{eq:def-average} lying in a single summand labelled by 
 $V'_1=gV_1,$, $V'_2=gV_2$. 
 
 \end{itemize}   
 
 \begin{ex}
 Suppose $\Ob(\calV)=\{V\}$ consists of one object, so $\calV$ is given by
 an algebra $A=\End_\calV(V)$ with $G$-action. Then 
 \[
 \xi = \xi_{V,V}: A\lra A^G, \quad a\mapsto \ul a = {1\over |G|} \sum_{g\in G} ga
 \]
 is the usual averaging map from $A$ to its invariant subalgebra. 
 \end{ex}

This shows that unlike the case of $\calV_G$ from \S \ref{subsec:coinv-orbits}, the maps $\xi_{V_1, V_2}$ do
not, in general,  extend to a functor $\calV\to\calV^G$.  

\medskip

Assume now that  the $G$-action on $\Ob\calV$ is free. 
In this situation one can use the maps $\xi_{V_1, V_2}$ to compare the construction 
of $\calV^G$ with that of  \S \ref{subsec:coinv-orbits}
 in the case when both are defined.
 That is, for any category $\calE$ let $\k[\calE]$ denote its $\k$-linear envelope.
  It is a new category with the same objects as $\calE$
 and morphisms being $\k$-vector spaces spanned by the sets of morphisms in $\calE$.  
 
 \begin{prop}
(a)  Let $\calE$ be a category with finitely many objects and a strict action of a finite group $G$ which is free on
  $\Ob\calE$. Then we have an isomorphism of categories
 $
 \k[\calE_G] \= \k[\calE]^G.  
 $
 \vskip .2cm
 
 (b) Let $\calV$ be a $\k$-linear category with finitely many objects and a strict action of a finite group $G$ which is free on
  $\Ob\calV$. Then we have an equivalence of categories
 $
 \calV^G\= \calV\rtimes G. 
 $
 \end{prop} 
 
 \noindent{\sl Proof:} (a) This boils down to the following observation : if $G$ is a finite group acting freely on a finite set $S$,  then
 the two $\k$-vector spaces $(\bigoplus_{s \in G} \k s)^G$ and $\bigoplus_{\overline{s} \in G \backslash S} \k \overline{s}$ are canonically isomorphic.
 
 \vskip .2cm
 
 (b) Let $V,W\in\Ob(\calV)$, so the  orbits $GV, GW$ are objects of $\calV^G$. By definition and
 by freeness of the action on $\Ob(\calV)$, 
 \[
 \Hom_{\calV\rtimes G}(V,W) = \bigoplus_{g\in G} \Hom_\calV(gV, W) = \biggl(
 \bigoplus_{g,h\in G} \Hom(\calV(gV, hW)\biggr)^G
 =
 \Hom_{\calV^G}(GV, GW). 
 \]
 These identifications for all $V,W$ are compatible with composition of morphisms and so define
 a functor $\calV\rtimes G\to \calV^G$. This functor is surjective on objects and bijective on
 Hom-spaces, i.e., 
 is an equivalence of cateories. \qed 
 
 \vskip .2cm
 
 We now relate the construction of $\calV^G$ with other constructions in the Appendix.  
 First, for any $G$-orbit $O\subset\Ob(\calV)$ we form the object 
 \[
 V_O = \bigoplus_{V\in O} V\, \in \, \calV^\oplus.
 \]
 The object $V_O$ is naturally $G$-equivariant, and \eqref{eq:orbit-cat-inv-hom} is the same as 
  $\Hom_{(V^\oplus)^{[G]}} (V_O, V_{O'})$. Further, as in  Proposition \ref{V^G=VXG}, 
  $V_O$ considered as an object of $\calV^\oplus \rtimes G$ has a natural $G$-action, so the object of $G$-invariants
  $V_O^G \in (V^\oplus\rtimes G)^\Kar$ is defined. If $V\in O$ is any object and $G_V\subset G$ is the
  stabilizer of $V$ (with respect to the $G$-action on $\Ob(\calV)$), then $V_O^G \= V^{G_V}$.  
  Next, the total algebra $A(\calV\rtimes G)$ contains the group algebra $\k[G]$, in particular, 
  the idempotent $\eps = {1\over |G|} \sum_{g\in G} [g]$. 
  
  \begin{prop}\label{prop:inv=kar}
  (a) The correspondence $O\mapsto V_O$ embeds $\calV^G$ as a full subcategory in $(\calV^\oplus)^{[G]}$.
  
  \vskip .2cm
  
  (b) The correspondence $O\mapsto V_O^G$ embeds $\calV^G$ as a full subcategory in $(\calV^\oplus \rtimes G)^\Kar$.
  In fact, it takes values (up to isomorphism) in $(\calV\rtimes G)^\Kar$ since $V_O^G\= V^{G_V}$. 
  
  \vskip .2cm
  
  (c) The total algebra $A(\calV^G) = A(\calV)^G$ is identified with the subalgebra $\eps A(\calV\rtimes G)\eps\subset
  A(\calV\rtimes G)$. 
  \end{prop}

\noindent{\sl Proof:} (a)  Follows because \eqref{eq:orbit-cat-inv-hom} is the same as
$\Hom_{(\calV^\oplus)^{[G]}}(V_O, V_{O'})$.

\vskip .2cm

 (b)  By definition of Karoubi completion and Proposition \ref{V^G=VXG},
\[
\Hom_{(\calV^\oplus\rtimes G)^\Kar } (V_O^G, V_{O'}^G) = \eps \Hom_{\calV^\oplus\rtimes G}(V_O, V_{O'})\eps
= \Hom_{\calV^\oplus\rtimes G}(V_O, V_{O'})^G = \Hom_{\calV^G}(O, O'). 
\]

(c) This is a consequence of (b). \qed
 
 \vfill\eject

 \subsection{Descent for sheaves} \label{subsec:desc-sh}
 Let $G$ be a finite group acting on a topological space $X$ (conventions of \S \ref{sec:notation}), so we have the quotient space $G\bs X$ and the
 projection $q: X\to G\bs X$.  It is convenient to  also consider the orbifold (or stack) quotient $[G\bs X]$
 so that sheaves on $[G\bs X]$ are $G$-equivariant sheaves on $X$, i.e., $\Sh_{[G\bs X]} = \Sh_X^{[G]}$. 
 The map $q$ yields a projection $p: [G\bs X]\to G\bs X$. In particular, we have the 
   pullback functor 
   \be\label{eq:p^{-1}}
   p^{-1}: \Sh_{G\backslash X} \lra \Sh_{[G\bs X]} = \Sh_X^{[G]}
   \ee
 which takes a sheaf $\calF\in\Sh_{G\bs X}$ to $q^{-1}(\calF)\in \Sh_X$ with its obvious $G$-equivariant structure. This functor
 has a right adjoint
 \be\label{eq:p_*}
 p_*: \Sh_{[G\bs X]} = \Sh_X^{[G]} \lra \Sh_{G\bs X}, \quad \calG \mapsto q_*(\calG)^G,
 \ee
 the subsheaf of invariants in the usual direct image. Note that both $p^{-1}, p_*$ are exact
 functors of abelian categories. 
 
 \begin{prop}\label{prop:p^-1emb}
 The unit $\eps: \Id_{\Sh_{G\bs X}}\to p_* p^{-1}$ is an isomorphism of functors.

 \end{prop}
 
 \noindent{\sl Proof:} 
  Let $\calG\in\Sh_{G\bs X}$. It is enough to show that for any $y\in G\bs X$ the induced map of stalks
 $\eps_{\calG, y}: \calG_y \to (p_*p^{-1}\calG)_y$ is an isomorphism. 
 Now, for any $x\in X$ we have $(q^{-1}\calG)_x = \calG_{q(x)}$. The preimage $q^{-1}(x)$ forms a single $G$-orbit
 in $X$, so 
 $(q_* q^{-1}\calG)_y$ is the direct sum of  identical copies of $\calG_y$, labelled by
 $q^{-1}(y)$ and transitively permuted by $G$.
 Therefore 
 \[
 (p_*p^{-1}\calG)_y = \biggl(\bigoplus\nolimits_{x\in q^{-1}(y)} \calG_y\biggr)^G = \calG_y,
 \]
 with $\eps_{\calG, y}$ giving the identification. 
   \qed
 
\vskip .2cm

Further, let $D^b\Sh_{G\bs X}$, $D^b\Sh_{[G\bs X]}$ be the bounded derived categories of $\Sh_{G\bs X}$ and  
$\Sh_{[G\bs X]}$ respectively. Since $G$ is a finite group and $\ch(\k)=0$, any $\k$-linear representation
of $G$ has trivial higher cohomology. So our definition of $D^b\Sh_{[G\bs X]}$ agrees with other,
more sophisticated definitions of the equivariant derived categories such as in
\cite{Bernstein-Lunts}. 

Since $p^{-1}, p_*$ are exact functors of abelian categories, they extend componentwise to
an adjoint pair of derived functors which we denote by the same symbols:
\[
\xymatrix{
D^b\Sh_{G\bs X} \ar@<0.1cm>[r]^{p^{-1}} 
&  D^b\Sh_{[G\bs X]} \ar@<0.1cm>[l]^{p_*}.
}
\]

\begin{corollary}\label{cor:der-unit}
The unit of the derived adjunction $\eps: \Id_{D^b\Sh_{G\bs X}}\to p_* p^{-1}$ is an isomorphism
of functors. \qed
\end{corollary}

In particular, Proposition \ref{prop:coreflect} gives that both versions of pullback:
\[
p^{-1}: \Sh_{G\bs X}\lra \Sh_{[G\bs X]}, \quad p^{-1}: D^b \Sh_{G\bs X}\lra D^b\Sh_{[G\bs X]}
\]
are embeddings of full (coreflexive) subcategories.

 \vfill\eject
  
\subsection {Descent for modules} \label{subsec:desc-alg} 
Let $A$ be an associative algebra with unit and
$G$ a finite group acting on $A$, see Example \ref{ex:G-equiv-sh-mod}.
We have an embedding of algebras
\[
\iota: A^G \hookrightarrow A[G], \quad a\mapsto a\cdot \eps, \quad \eps := {1\over |G|}
\sum_{g\in G} \,\, [g].
\]
  Note that unless $|G|=1$, the embedding $\iota$ does not
preserve units: $\iota(1)=\eps$ is the idempotent corresponding, under the isomorphism
$\k[G]\= \bigoplus_{V\in \Irr(G)} \End_\k(V)$, to the identity of the trivial representation. 

The embedding $\iota$ gives rise to the adjoint pair $(I,R)$, where
\be\label{eq:I,R-general}
\begin{gathered}
I: \, {_{A^G}\Mod} \lra {_{A[G]}\Mod}, \quad N\mapsto A[G]\otimes_{A^G}N,
\\
R:\,  {_{A[G]}\Mod} \lra{ _{A^G}\Mod}, \quad M\mapsto (\iota^*M)^u
\end{gathered}
\ee
 are the standard induction and restriction functors. Here $\iota^*M$ is $M$ considered
 as an $A^G$-module (which may not be unital) and ``$u$'' means the maximal
 unital submodule \eqref{eq:M^u}. 
 
 Viewing $A[G]$-modules as $G$-equivariant $A$-modules, we can describe the functors
 $I$ and $R$ more directly:
  \be\label{eq:I,R-explicit}
  \begin{gathered}
 I(N) = A\otimes_{A^G}N, \quad g(a\otimes n) := g(a) \otimes n, \quad N\in {_{A^G}\Mod},
 \\
 R(M) = M^G, \quad M\in {_{A[G]}\Mod} = {_A\Mod^{[G]}}. 
 \end{gathered}
 \ee
 In particular, $R$ is always exact (not just left exact) while $I$ is exact (not just right exact)
 iff $A$ is flat as an $A^G$-module. 
 
 \begin{prop}\label{prop:unit-RI}
 The unit $\eps: \Id_{_{A^G}\Mod}\to RI$ of the adjunction $(I,R)$ is an isomorphism of functors.
 \end{prop}
 
 In particular (Proposition \ref{prop:coreflect}), 
the functor $I$ embeds $_{A^G}\Mod$ as a full subcategory
 in $_{A[G]}\Mod$. 
 
 \vskip .2cm
 
 \noindent{\sl Proof:} The statement means that for any (unital) $A^G$-module $N$ the natural morphism
 \[
 \eps_N: N\lra (A\otimes_{A^G} N)^G, \quad n\mapsto 1_A\otimes n,
 \]
is an isomorphism.  This is manifestly true for $N=A^G$ and therefore for any free $A^G$-module
$F$. Now, any $N$ has a free resolution, as in the top row of the diagram below:
\[
\xymatrix{
\cdots \ar[r]& F^{-1}
\ar[d]_{\eps_{F^{-1}}} \ar[r] & F^0 \ar[r] \ar[d]_{\eps_{F^0}} & N
\ar[d]^{\eps_N} \ar[r]&0
\\
\cdots \ar[r]&(A\otimes_{A^G} F^{-1})^G  \ar[r] & (A\otimes_{A^G} F^0)^G \ar[r] & 
(A\otimes_{A^G}N)^G \ar[r]&0.
}
\]
 Since the functor $\otimes$ is right exact and taking $G$-invariants is exact as $\ch \k=0$, the bottom row of
 the diagram is exact as well. Since $\eps_{F^i}$ are isomorphisms for all $i$,
 so is $\eps_N$.
 \qed

\vfill\eject

\vskip 0.5cm

\small{

M.K.:  Kavli IPMU (WPI), UTIAS, The University of Tokyo, Kashiwa, Chiba 277-8583, Japan. Email: {\tt mikhail.kapranov@protonmail.com}

\vskip .2cm

V. S.: Institut Math\'ematique de Toulouse, Universit\'e Paul Sabatier, 118 rue de Narbonne, 31062 Toulouse, France, and
 Kavli IPMU (WPI), UTIAS, The University of Tokyo, Kashiwa, Chiba 277-8583, Japan.
Email: {\tt schechtman@math.ups-toulouse.fr}

\vskip .2cm

O.S.: Laboratoire de Math\'ematique d'Orsay, Universit\'e Paris Saclay, 91405 Orsay, France.
Email: {\tt olivier.schiffmann@universite-paris-saclay.fr}

\vskip .2cm

J.Y.: Laboratoire de Math\'ematique d'Orsay, Universit\'e Paris Saclay, 91405 Orsay, France.
Email: {\tt jiangfan.yuan@universite-paris-saclay.fr}
}

 \end{document}